\documentclass[10pt]{article}
\usepackage{amsmath, amssymb}
\usepackage{amsthm} 
\usepackage[all]{xy}
\usepackage[pdftex]{graphicx}
\usepackage{tikz}

\theoremstyle{plain} 

\theoremstyle{definition}

\def\Sch{\text{\rm Sch}}
\def\F{\mathcal F}
\def\L{{\mathcal L}}

\def\bgn{\begin}

\def\CL{\text{\rm CL}}

\def\E{\bold E}
\def\J{{\mathcal J}}
\def\L{{\mathcal L}}

\def\1{{[1]}}
\def\2{{[2]}}
\def\3{{[3]}}
\def\({\left(}
\def\){\right)}
\def\s-circ{\,{\scriptstyle{\circ}}\,}
\def\<<{<\negthinspace \negthinspace<}
\def\Ad{\text{\rm Ad}}
\def\ad{\text{\rm ad}}

\def\even{\text{\rm even}}
\def\bgn{\begin}

\def\endaln{\end{align}}
\def\cou{\ss\text{\rm co}}

\def\<{<\negthinspace \negthinspace <}

\def\t{\theta}

\def\({\left(}
\def\){\right)}

\def\Re{\text{\rm Re}}
\def\[{\big[\neg\big[}
\def\]{\big]\neg\big]}
\def\al{\al}

\def\M{{\mathcal M}}

\def\Vol{\text{\rm Vol}}
\def\tr{\text{\rm tr}}
\def\a{\alpha}
\def\b{\beta}
\def\bs{\backslash}
\def\e{\varepsilon}

\def\Gam{\Gamma}

\def\del{\delta}
\def\lam{\lambda}
\def\ome{\omega}
\def\Ome{\Omega}
\def\sig{\sigma}

\def\A{{\mathcal A} }
\def\E{{\mathcal E}}

\def\D{\Bbb D}

\def\R{\Bbb R}
\def\C{\Bbb C}

\def\O{\Bbb O}

\def\M{\frak M}

\def\w{\wedge}
\def\({\left(}
\def\){\right)}

\def\nab{\nabla}
\def\neg{\negthinspace}
\def\h{\hat}

\def\til{\tilde}

\def\l{\left}
\def\ol{\overline}
\def\pa{\partial}
\def\olpa{\ol{\partial}}
\def\r{\right}
\def\ran{\rangle} 
\def\lan{\langle}

\def\ss{\scriptscriptstyle}
\def\trian{\triangle}

\def\arrow{\longrightarrow}

\def\bsh{\backslash}

\def\:{\, :\,}
\def\CL{\text{\rm CL}}
\def\TT{T\oplus T^*}
\def\complex{generalized complex }
\def\K\"ahler{generalized K\"ahler}
\def\vol{\text{\rm vol}}

\def\10{\displaystyle L^{10}}
\def\2{\displaystyle L^2}
\def\c0{\displaystyle C^0}
\def\dstyle{\displaystyle}
\def\10{\displaystyle L^{10}}
\def\2{\displaystyle L^2}
\def\del{\delta}
\def\del2{\displaystyle L^2_{0,\delta}}
\def\c0{\displaystyle C^0}
\def\dstyle{\displaystyle}

\def\del{\delta}

\def\K{{\mathcal K}}
\def\M-A{\text{\rm Monge-Amp\`ere}}
\def\O{{\mathcal O}}

\def\M-A{\text{\rm Monge-Amp\`ere}}
\def\[{\big[\,}
\def\]{\,\big]}
\def\End{\text{\rm End\,}}
\def\Hom{\text{\rm Hom}}

\def\id{\text{\rm id}}

\def\top{\text{\rm top}}
\def\rank{\text{\rm rank}}
\def\even{\text{\rm even}}
\def\odd{\text{\rm odd}}
\def\GK{\text{\rm generalized K\"ahler }}
\def\L{\text{{\mathcal L}}}
\def\TT{TM\oplus T^*M}

\def\M{{\mathcal M}}

\def\ol{\overline}
\def\part{\partial}

\def\tt{\ss T\oplus T^*}
\def\D{{\mathcal D}}
\def\F{{\mathcal F}}
\def\L{{\mathcal L}}
\def\Lam{\Lambda}
\def\Herm{\text{\rm Herm}}

\def\ad{\text{\rm ad}}
\def\rk{\text{\rm rk}}

\title{The Kobayashi-Hitchin correspondence of generalized holomorphic vector bundles over generalized K\"ahler manifolds\\
of symplectic type
}%
\author{Ryushi Goto}
\date{} 
\begin{document}
\maketitle

\begin{abstract}\noindent
In the previous paper \cite{Goto_2017},
the notion of an Einstein-Hermitian metric of a generalized holomorphic vector bundle over a generalized K\"ahler manifold of symplectic type was introduced from the moment map framework. 
In this paper we establish the Kobayashi-Hitchin correspondence, that is, the equivalence of
the existence of an Einstein-Hermitian metric and $\psi$-polystability of a generalized holomorphic vector bundle over a compact generalized K\"ahler manifold of symplectic type.
Poisson modules provide intriguing generalized holomorphic vector bundles and 
we obtain $\psi$-stable Poisson modules over complex surfaces which are not stable in the ordinary sense.
\end{abstract}
\tableofcontents
\numberwithin{equation}{section}

\noindent
\section{Introduction}
\footnote{MSC classes:	53D18, 53D17, 53C07}\\
\noindent
S.~Kobayashi \cite{S.Kob_1987} and L\"ubke \cite{Lu_1983} showed that an Einstein-Hermitian vector bundle over a smooth compact K\"ahler manifold is polystable. Conversely, Donaldson 
\cite{Do_1983}, \cite{Do_1985}, \cite{Do_1987} and Uhlenbeck-Yau \cite{UY_1986}, \cite{UY_1989} proved that a polystable holomorphic vector bundle has an Einstein-Hermitian metric. The equivalence of
the existence of an Einstein-Hermitian metric and polystability of a holomorphic vector bundle
 is called the Kobayashi-Hitchin correspondence (See also \cite{Hi_2015}
 for more details).

Generalized complex geometry and generalized K\"ahler geometry are introduced by Hitchin \cite{Hi_2003} and Gualtieri which are successful 
generalizations of the ordinary complex geometry and K\"ahler geometry.
Holomorphic Poisson geometry is closely related to generalized K\"ahler geometry.
In fact, it turns out that holomorphic Poisson structures on a K\"ahler manifold produce  
remarkable  generalized K\"ahler structures of symplectic type \cite{Hi_2007}, \cite{Lin_2006},  \cite{Goto_2009}, \cite{Goto_2010}, \cite{Goto_2012}, \cite{ABD_2015},\cite{Gua_2010}.
Poisson modules also provide typical and prominent generalized holomorphic vector bundles \cite{Hi_2011-1}, \cite{Hi_2011-2}, \cite{Gua_2010}.\par
Then an inevitable and important problem is to extend the Kobayashi-Hitchin correspondence to the cases of generalized holomorphic vector bundles over a generalized K\"ahler manifold.
In the paper, 
from a view point of moment map framework, we focus on  the class of generalized K\"ahler manifolds of symplectic type{$^{*1}$}.
\footnote{{${}^{*1}$} Our results can be applied to more general cases that $\psi$ is a $d$-closed, nondegenerate, pure spinor.
In fact, the existence result as in Section \ref{Construction of Einstein-Hermitian metrics on stable generalized holomorphic bundles}
does hold. We also extend our result to the cases $H$-twisted generalized K\"ahler structures.} A generalized K\"ahler manifold $(M, \J_1, \J_2)$ is {\it of symplectic type} if one of generalized complex structures $\J_2$ is induced from 
 a symplectic form $\ome$ and a real $d$-closed $2$-form $b$ on $M,$
 which is denoted by $\J_\psi$,where $\psi=e^{b-\sqrt{-1}\ome},$
is a $d$-closed, nondegenerate, pure spinor of
type $0$. Let $E$ be a generalized holomorphic vector bundle and $h$ an Hermitian metric over 
 a generalized K\"ahler manifold of symplectic type
$(M, \J, \J_\psi)$.
The Einstein-Hermitian condition in \cite{Goto_2017}
 is defined by the following equation:
\bgn{equation}\label{Einstein-Hermitian equation}
\K_h(\psi)=\lam \psi\,\id_E,
\end{equation}
where $\K_h(\psi)$ denotes the mean curvature of the canonical generalized connection
(see Definition \ref{Einstein-Hermitian condition} and Definition \ref{Einstein-Hermitian metric}). 
The degree of $E$ is defined in terms of the first Chern class of a vector bundle $E$ together with
the class $[\psi]$ and then the slope is also given by 
deg $E/\rank\,E$ .
By using the slope inequality,  we have the notions of $\psi$-stability and $\psi$-polystability of a generalized holomorphic vector bundle. 
Then the following is our main theorem which is proved in Subsection \ref{Proof of main theorem}.
\bgn{theorem}\label{kobayashi-Hitchin correspondence}\text{\rm[Kobayashi-Hitchin correspondence]}
There exists an Einstein-Hermitian metric on a $\psi$-polystable generalized holomorphic vector bundle.
Conversely, a generalized holomorphic vector bundle admitting an Einstein-Hermitian metric is $\psi$-polystable.
\end{theorem}
 
It is known that a generalized K\"ahler structure $(\J_1, \J_2)$ gives rise to
a bihermitian structure $(I_+, I_-, g, b)$, where 
$I_\pm$ are two ordinary complex structures and a single $g$ is an Hermitian metric with respect to both $I_+$ and $I_-$, and $b$ is a real $2$-form.
Then a generalized holomorphic vector bundle provides a holomorphic vector bundle with respect to both $I_+$ and $I_-.$
From the view point of bihermitian geometry, Hitchin \cite{Hi_2011-1} expected in general the following equation describes a stability condition of a generalized holomorphic vector bundle: 
\bgn{equation}\label{Hitchin equation}
\frac12\(F_+\w\ome_+^{n-1}+ F_-\w\ome_-^{n-1}\)=\lam\, \id_E \vol_M
\end{equation}
In \cite{HMS_2016},  Hu, Moraru and Reza  showed that the equation (\ref{Hitchin equation})
is related with
$\a$-stability which is introduced by using Hermitian forms $\ome_\pm$ of the bihermitian structure if both $\ome_\pm $ are Gauduchon metrics.
However, a generalized K\"ahler manifold $(M,\J_1, \J_2)$ does not give 
Gauduchon metrics $\ome_\pm$ if the dimension of $M$ is greater than $4.$
The $\psi$-stability in this paper is topologically defined as in the ordinary K\"ahler cases, which is different from the $\a$-stability. 
Our Einstein-Hermitian equation 
(\ref{Einstein-Hermitian equation}) is introduced from the moment map framework, which is
also different from the equation (\ref{Hitchin equation}) in general${}^{*2}$.
{\footnote {${}^{*2}$We also remark that in the cases of co-Higgs bundles over an ordinary K\"ahler manifold, 
our equation (\ref{Einstein-Hermitian equation}) coincides with the one as (\ref{Hitchin equation}).}}

This paper is organized as follows.
In Section \ref{Generalized complex structures and generalized Kahler structures}, 
we shall give a brief review of generalized complex structures and generalized K\"ahler structures focusing on nondegenerate, pure spinors. 
In Subsection \ref{The stability theorem of generalized Kahler manifolds}, 
we recall the stability theorem of generalized K\"ahler structures which is applied to construct nontrivial examples of generalized K\"ahler manifolds from holomorphic Poisson structures. 
In Section \ref{Einstein-Hermitian generalized connections} and \ref{Generalized holomorphic vector bundles},
we recall definitions of Einstein-Hermitian metrics and generalized holomorphic vector bundles and the canonical generalized connections. \par
In the cases of ordinary holomorphic vector bundles over a K\"ahler manifold, 
a weak holomorphic subbundle plays a crucial role to construct an Einstein-Hermitian metric.
In Section \ref{Weak generalized holomorphic subbundles}, 
we also introduce the notion of a weak 
generalized holomorphic subbundle.
Then the notions of $\psi$-stability and $\psi$-polystability are introduced 
in Section \ref{Definition of stability and polystability}.
In Section \ref{The second fundamental forms}, 
we obtain the formula of the second fundamental form which measures the difference between 
the first Chern form of $E$ and the one of a subbundle of $E$.
The formula is used to 
show $\psi$-polystability of a generalized holomorphic vector bundle admitting an Einstein-Hermitian metric in Section \ref{From Einstein-Hermitian metrics to stability}.
In Section \ref{Variation formula of mean curvature}, we give variation formulas of the curvature and the mean curvature of 
canonical generalized Hermitian connections under deformations of Hermitian metrics on a generalized holomorphic vector bundle. 
In Section \ref{Construction of Einstein-Hermitian metrics on stable generalized holomorphic bundles}, 
we construct an Einstein-Hermitian metric on a stable generalized holomorphic vector bundle. 
We use the continuity method. 
In Section \ref{Einstein-Hermitian metrics and stable Poisson modules}, 
we construct Poisson modules by using the Serre construction over complex surfaces and discuss the $\psi$-stability of them. 
On $\C P^2,$ by using a certain configuration of points on a line, we obtain a $\psi$-stabile Poisson 
module which is not stable in the ordinary sense. 
Thus such a Poisson module does have an Einstein-Hermitian metric as a generalized holomorphic vector bundle, however which does not have any ordinary Einstein-Hermitian metric. 
In Section \ref{Vanishing theorems of generalized holomorphic vector bundles}, 
we obtain vanishing theorems of a generalized holomorphic vector bundle with an Einstein-Hermitian metric over a generalized K\"ahler manifold of symplectic type, which give another proof of $\psi$-polystability of an Einstein-Hermitian generalized holomorphic vector bundle.
\section{Generalized complex structures and generalized K\"ahler structures}\label{Generalized complex structures and generalized Kahler structures}
\subsection{Generalized complex structures and nondegenerate, pure spinors}
Let $M$ be a differentiable manifold of real dimension $2n$.
The bilinear form $\lan\,\,,\,\,\ran_{\tt}$ on 
the direct sum $T_M \oplus T^*_M$ over a differentiable manifold $M$ of dim$=2n$ is defined by 
$$\lan v+\xi, u+\eta \ran_{\tt}=\frac12\(\xi(u)+\eta(v)\),\quad  u, v\in T_M, \xi, \eta\in T^*_M .$$
Let SO$(\TT)$ be the fibre bundle over $M$ with fibre SO$(2n, 2n)$ which is 
a subbundle of End$(\TT)$  preserving the bilinear form $\lan\,\,\,,\,\,\,\ran_{\tt}$ 
 An almost \complex structure $\J$ is a section of SO$(\TT)$ satisfying $\J^2=-\id.$ Then as in the case of almost complex structures, an almost \complex structure $\J$ yields the eigenspace decomposition :
\bgn{equation}\label{eigenspace decomposition}
(T_M \oplus T^*_M)^\C =\L_\J \oplus \ol \L_\J,
\end{equation} where 
$\L_\J$ is $-\sqrt{-1}$-eigenspace and  $\ol{\L}_\J$ is the complex conjugate of $\L_\J$. 
The Courant bracket of $\TT$ is defined by 
$$
 [u+\xi, v+\eta]_{\cou}=[u,v]+{\mathcal L}_u\eta-{\mathcal L}_v\xi-\frac12(di_u\eta-di_v\xi),
 $$
 where $u, v\in TM$ and $\xi, \eta$ is $T^*M$.
If $\L_\J$ is involutive with respect to the Courant bracket, then $\J$ is a generalized complex structure, that is, $[e_1, e_2]_{\cou}\in \Gam(\L_\J)$  for any two elements 
 $e_1=u+\xi,\,\, e_2=v+\eta\in \Gam(\L_\J)$.
Let $\CL(T_M \oplus T^*_M)$ be the Clifford algebra bundle which is 
a fibre bundle with fibre the Clifford algebra $\CL(2n, 2n)$ with respect to $\lan\,,\,\ran_{\tt}$ on $M$.
Then a vector $v$ acts on the space of differential forms $\oplus_{p=0}^{2n}\w^pT^*M$ by 
the interior product $i_v$ and a $1$-form $\t$ acts on $\oplus_{p=0}^{2n}\w^pT^*M$ by the exterior product $\t\w$, respectively.
Then the space of differential forms gives a representation of the Clifford algebra $\CL(\TT)$ which is 
the spin representation of $\CL(\TT)$. 
Thus
the spin representation of the Clifford algebra arises as the space of differential forms $$\w^\bullet T^*_M=\oplus_p\w^pT^*_M=\w^{\even}T^*_M\oplus\w^{\odd}T^*_M.$$ 
The inner product $\lan\,,\,\ran_s$ of the spin representation is given by 
$$
\lan \a, \,\,\,\b\ran_s:=(\a\w\sig\b)_{[2n]},
$$
where $(\a\w\sig\b)_{[2n]}$ is the component of degree $2n$ of $\a\w\sig\b\in\oplus_p \w^pT^*M$ and 
$\sig$ denotes the Clifford involution which is given by 
$$
\sig\b =\bgn{cases}&+\b\qquad \deg\b \equiv 0, 1\,\,\mod 4 \\ 
&-\b\qquad \deg\b\equiv 2,3\,\,\mod 4\end{cases}
$$
We define $\ker\Phi:=\{ e\in (T_M\oplus T^*_M)^\C\, |\, e\cdot\Phi=0\, \}$ for a differential form $\Phi
\in \w^{\even/\odd}T^*_M.$
If $\ker\Phi$ is maximal isotropic, i.e., $\dim_\C\ker\Phi=2n$, then $\Phi$ is called {\it a pure spinor} of even/odd type.
A pure spinor $\Phi$ is {\it nondegenerate} if $\ker\Phi\cap\ol{\ker\Phi}=\{0\}$, i.e., 
$(T_M\oplus T^*_M)^\C=\ker\Phi\oplus\ol{\ker\Phi}$.
Then a nondegenerate, pure spinor $\Phi\in \w^\bullet T^*_M$ gives an almost generalized complex structure $\J_{\Phi}$ which satisfies 
$$
\J_\Phi e =
\bgn{cases}
&-\sqrt{-1}e, \quad e\in \ker\Phi\\
&+\sqrt{-1}e, \quad e\in \ol{\ker\Phi}
\end{cases}
$$
Conversely, an almost \complex structure $\J$ locally arises as $\J_\Phi$ for a nondegenerate, pure spinor $\Phi$ which is unique up to multiplication by
non-zero functions.  Thus an almost \complex structure yields the canonical line bundle $K_{\J}:=\C\lan \Phi\ran$ which is a complex line bundle locally generated by a nondegenerate, pure spinor $\Phi$ satisfying 
$\J=\J_\Phi$.
An \complex structure 
$\J_\Phi$ is integrable if and only if $d\Phi=\eta\cdot\Phi$ for a section $\eta\in T_M\oplus T^*_M$. 
The {\it type number} of $\J=\J_\Phi$ is defined as the minimal degree of the differential form $\Phi$. Note that type number Type $\J$ is a function on a manifold which is not a constant in general.
\bgn{example}
Let $J$ be a complex structure on a manifold $M$ and $J^*$ the complex structure on the dual  bundle $T^*M$ which is given by $J^*\xi(v)=\xi (Jv)$ for $v\in TM$ and $\xi\in T^*M$.
Then a \complex structure $\J_J$ is given by the following matrix
$$\J_J=\bgn{pmatrix}J&0\\0&-J^*
\end{pmatrix},$$
Then the canonical line bundle is the ordinary one which is generated by complex forms of type $(n,0)$.
Thus we have  Type $\J_J =n.$
\end{example}
\bgn{example}
Let $\ome$ be a symplectic structure on $M$ and $\h\ome$ the isomorphism from $TM$ to $T^*M$ given by $\h\ome(v):=i_v\ome$. We denote by $\h\ome^{-1}$ the inverse map from $T^*M$ to $TM$.
Then a \complex structure $\J_\psi$ is given by the following
$$\J_\psi=\bgn{pmatrix}0&-\h\ome^{-1}\\
\h\ome&0
\end{pmatrix},\quad\text{\rm Type $\J_\psi =0$}$$
Then the canonical line bundle is given by the differential form $\psi=e^{-\sqrt{-1}\ome}$. 
Thus Type $\J_\psi=0.$
\end{example}
\bgn{example}[$b$-field action]
A real $d$-closed $2$-form $b$ acts on a \complex structure by the adjoint action of Spin group $e^b$ which provides
a \complex structure $\Ad_{e^b}\J=e^b\circ \J\circ e^{-b}$. 
\end{example}
\bgn{example}[Poisson deformations]\label{Poisson deformations}
Let $\b$ be a holomorphic Poisson structure on a complex manifold. Then the adjoint action of Spin group $e^\b$ gives deformations of new \complex structures by 
$\J_{\b t}:=\Ad_{\b^{Re} t}\J_J$.  Then Type ${\J_{\b t}}_x=n-2$ (rank of $\b_x$) at $x\in M$,
which is called the Jumping phenomena of type number.
\end{example}
Let $(M, \J)$ be a generalized complex manifold and $\ol \L_\J$ the eigenspace of eigenvalue $\sqrt{-1}$.
Then we have the Lie algebroid complex $\w^\bullet\ol{\L}_\J$:
$$
0\arrow\w^0\ol \L_\J\overset{\ol\pa_\J}\arrow\w^1\ol \L_\J\overset{\ol\pa_\J}\arrow\w^2\ol \L_\J\overset{\ol\pa_\J}\arrow\w^3\ol \L_\J\arrow\cdots
$$
The Lie algebrid complex is the deformation complex of \complex structures. 
In fact, $\e\in \w^2\ol \L_\J$ gives deformed isotropic subbundle 
$E_\e:=\{ e+[\e, e]\, |\, e\in \L_\J\}$. 
Then $E_\e$ yields deformations of \complex structures if and only if $\e$ satisfies Generalized Mauer-Cartan equation
$$
\ol{\pa}_\J\e+\frac12[\e, \e]_{\Sch}=0,
$$
where $[\e, \e]_{\Sch}$ denotes the Schouten bracket. 
The Kuranishi space of generalized complex structures is constructed.
Then the second cohomology group $H^2(\w^\bullet\ol \L_\J)$ of the Lie algebraic complex gives the infinitesimal deformations of \complex structures and the third one 
$H^3(\w^\bullet\ol \L_\J)$ is the obstruction space to deformations of \complex structures.
Let $\{e_i\}_{i=1}^n$ be a local basis of $\L_\J$ for an almost \complex structure $\J$, 
where $\lan e_i, \ol e_j\ran_{\tt}=\del_{i,j}$.
The the almost \complex structure $\J$ is written as an element of Clifford algebra,
$$
\J=\frac{\sqrt{-1}}2\sum_i e_i\cdot\ol {e}_i,
$$
where $\J$ acts on $\TT$ by the adjoint action $[\J, \,]$. 
Thus we have $[\J, e_i]=-\sqrt{-1}e_i$ and $[\J, \ol e_i]=\sqrt{-1}e_i$.
An almost \complex structure $\J$ acts on differential forms by the Spin representation which gives the decomposition into eigenspaces:
\bgn{equation}\label{eq:bunkai}
\w^\bullet T^*_M=U^{-n}\oplus U^{-n+1}\oplus\cdots U^{n},
\end{equation}
where $U^{i}(=U^i_\J)$ denotes the $i$-eigenspace. 
\subsection{Generalized K\"ahler structures}
\bgn{definition}
{\it A generalized K\"ahler structure} is a pair $(\J_1, \J_2)$ consisting of two commuting \complex structures 
$\J_1, \J_2$ such that $\h G:=-\J_1\circ\J_2=-\J_2\circ \J_1$ gives a positive definite symmetric form 
$G:=\lan \h G\,\,,  \,\,\ran$ on $T_M\oplus T_M^*$, 
We call $G$ {\it a generalized metric}.
{\it A \GK structure of symplectic type} is a \GK structure $(\J_1, \J_2)$ such that  $\J_2$ is the \complex structure $\J_\psi$ which is  induced from  a $d$-closed, nondegenerate, pure spinor $\psi:=e^{b-\sqrt{-1}\ome}.$
\end{definition}
Each $\J_i$ gives the decomposition $(\TT)^\C=\L_{\J_i}\oplus\ol \L_{\J_i}$ for $i=1,2$.
Since $\J_1$ and $\J_2$ are commutative, we have the simultaneous eigenspace decomposition 
$$
(\TT)^\C=(\L_{\J_1}\cap \L_{\J_2})\oplus (\ol \L_{\J_1}\cap \ol \L_{\J_2})\oplus (\L_{\J_1}\cap \ol \L_{\J_2})\oplus
(\ol \L_{\J_1}\cap \L_{\J_2}).
$$
Since $\h G^2=+\id$,
The generalized metric $\h G$ also gives the eigenspace decomposition: $\TT=C_+\oplus C_-$, 
where $C_\pm$ denote the eigenspaces of $\h G$ of eigenvalues $\pm1$. 
We denote by $\L_{\J_1}^\pm$ the intersection $\L_{\J_1}\cap C^\C_\pm$. 
Then it follows 
\bgn{align*}
&\L_{\J_1}\cap \L_{\J_2}=\L_{\J_1}^+,  \quad \ol \L_{\J_1}\cap \ol \L_{\J_2}=\ol \L_{\J_1}^+\\
&\L_{\J_1}\cap \ol \L_{\J_2}=\L_{\J_1}^-,\quad \ol \L_{\J_1}\cap \L_{\J_2}=\ol \L_{\J_1}^-
\end{align*}
\bgn{example}
Let $X=(M, J,\ome)$ be a K\"ahler manifold. Then the pair $(\J_J, \J_\psi)$ is a generalized K\"ahler where 
$\psi=\exp(\sqrt{-1}\ome)$. 
\end{example}
\bgn{example}
Let $(\J_1, \J_2)$ be a generalized K\"ahler structure. 
Then the action of $b$-fields gives a generalized K\"ahler structure 
$(\Ad_{e^b}\J_1, \Ad_{e^b}\J_2)$ for a real $d$-closed $2$-form $b.$
\end{example}
\subsection{The deformation-stability theorem of generalized K\"ahler manifolds}
\label{The stability theorem of generalized Kahler manifolds}
It is known that the deformation-stability theorem of ordinary K\"ahler manifolds holds
\bgn{theorem}[Kodaira-Spencer]
Let $X=(M,J)$ be a compact K\"ahler manifold and $X_t$ small deformations of $X=X_0$ as complex manifolds.
Then $X_t$ inherits a K\"ahler structure. 
\end{theorem}
The following deformation-stability theorem of generalized K\"ahler structures provides many interesting examples of generalized K\"ahler manifolds of symplectic type.
\bgn{theorem}{\rm \cite{Goto_2010}}
Let $X=(M,J,\ome)$ be a compact K\"ahler manifold and $(\J_J, \J_\psi)$ the induced generalized K\"ahler structure, 
where $\psi=e^{-\sqrt{-1}\ome}$. 
If there are analytic deformations $\{\J_t\}$ of $\J_0=\J_J$ as \complex structures, then there are deformations of $d$-closed nondegenerate, pure spinors $\{\psi_t\}$ such that 
pairs $(\J_t, \J_{\psi_t})$ are generalized K\"ahler structures, where $\psi_0=\psi$
\end{theorem}
Then we have the following:
\bgn{corollary}
Let $X=(M,J, \ome)$ be a compact K\"ahler manifold with a nontrivial holomorphic Poisson structure $\b.$
Then there exists a nontrivial deformations of generalized K\"ahler structures $(\J_{\b t}, \, \J_{\psi_t})$ such that 
$\{\J_{\b t}\}$ is the Poisson deformations given by Example \ref{Poisson deformations} and $\{\psi_t\}$ is a family of $d$-closed nondegenerate, pure spinors and $\psi_0=e^{-\sqrt{-1}\ome}.$
\end{corollary}
\section{Einstein-Hermitian generalized connections}
\label{Einstein-Hermitian generalized connections}
\subsection{Generalized connections over symplectic manifolds}
Let $M$ be a compact real manifold of dimension $2n$ and $\ome$ a real symplectic structure on $M$. 
We denote by $\psi$ the exponential of $b-\sqrt{-1}\ome$, that is, 
$$
\psi:=e^{b-\sqrt{-1}\ome}=1+(b-\sqrt{-1}\ome)+\frac1{2!}(b-\sqrt{-1}\ome)^2+\cdots+\frac1{n!}(b-\sqrt{-1}\ome)^n,
$$
where $b$ denotes a real $d$-closed $2$-form on $M.$
Then $\psi$ induces a generalized complex structure $\J_\psi$ which gives a decomposition 
$$
(\TT)^{\Bbb C}=\L_{\J_\psi}\oplus \ol {\L_{\J_\psi}}
$$
$\TT$ acts on differential forms by the spin representation which is given by the interior product and the exterior product of $\TT$ on differential forms. 
Then as in (\ref{eq:bunkai}), we have a decomposition of differential forms on $M,$
$$
\oplus_{i=0}^{2n}\w^iT^*_M=\oplus_{j=0}^{2n} U_{\J_\psi}^{-n+j}, 
$$
where $U_{\J_\psi}^{-n}$ is a complex line bundle generated by $\psi$  and 
$U_{\J_\psi}^{-n+i}$ is constructed by the spin action of $\w^i\ol{\L_{\J_\psi}}$ on $U^{-n}_{\J_\psi}.$
Let $E\to M $ be a complex vector bundle of rank $r$ over $M$ and $\Gam(E)$ a set of smooth sections of $E$.
We denote by $\Gam(E\otimes(\TT)^\C)$ the set of smooth sections of $E\otimes(\TT)^\C$. 
{\it A generalized connection}
$\D^\A$
is a map from $\Gam(E)$ to $\Gam(E\otimes(\TT)^\C)$ such that 
$$
\D^\A(fs)=s\otimes df+f\D^\A(s),\qquad \text{\rm for } s\in\Gam(E),\,\, f\in C^\infty(M).
$$
Let $h$ be an Hermitian metric of $E$. {\it A generalized Hermitian connection } is a generalized connection $\D^\A$ satisfying 
$$
dh(s, s')=h(\D^\A s, s')+h(s, \D^\A s'), \qquad \text{\rm for } s, s'\in \Gam(E).
$$
We denote by $u(E,h)(:=u(E))$ the set of skew-symmetric endomorphisms of $E$ with respect to $h.$ Then $\End(E)$ is decomposed as
\bgn{equation}\label{eq:End(E)}
\End(E)=u(E,h)\oplus \Herm(E,h),
\end{equation}
where $\Herm(E,h)$ denotes the set of Hermitian endmorphisms of $E.$
Let $\{U_\a\}$ be an open covering of $M$ which gives local trivializations of $E$. 
We take  $s_\a:=(s_{\a,1}, \cdots, s_{\a,r})$ as  a local unitary frame of $E$ over  $U_\a$. 
The set of transition functions is denoted by $\{g_{\a\b}\}$.
Then an Hermitian generalized  connection $\D^\A$ is written as 
$$\D^\A(s_{\a,p})=\sum_{q=1}^r s_{\a,q}\A_{p,\a}^q,
$$
where $\A_\a:=\(\A_{p, \a}^q\)$ is called a connection form of a generalized connection which is a section of $u(E)\otimes_\R(\TT)|_{U_\a}$. Then $\A_\a$ is denoted as 
$$
\A_\a=\sum_i \A_{\a, i} e_i,
$$
where $e_i\in (\TT)$ and $\A_\a^i\in u(E).$
Note that each $e_i$ is a real element of $\TT$.
Then $\A_\a$ is also decomposed into  
$$\A_{\a}=A_{\a}+V_{\a}.$$  
Then it turns out that $A_\a:=\(A_{p,\a}^q\)$ is an ordinary connection form and 
$V_\a:=\(V_{p,\a}^q\)$ gives a section of $u(E)\otimes T_M$
In fact, 
by using local trivializations $s_\a$, given a generalized connection $\D^\A$ is written as 
$$
\D^\A= d_\a+  \A_\a,
$$
Since 
$\D^\A :\Gam (E)\to \Gam (E\otimes (\TT))$ is globally defined as a differential operator, we have 
\bgn{align}
A_\a=&-(dg_{\a\b})g^{-1}_{\a\b}+g_{\a\b}A_\b g^{-1}_{\a\b}\\
V_\a=&g_{\a\b}V_\b g^{-1}_{\a\b}
\end{align}
Thus it follows that $A_\a$ is a connection form of a connection $D^A$ and  $V_\a$ is a $u(E)$-valued vector field. 
As shown in the ordinary connections, a generalized connection $\D^\A $ is extended to be an operator $\Gam(\End(E))\to \Gam(\End(E)\otimes(\TT))$
by the following: 
$$
(\D^\A \xi) s:=\D^\A(\xi s)-\xi(\D^\A s), \qquad \text{\rm for }\xi\in \Gam(\End(E)), \,\, s\in\Gam(E).
$$
Then the extended operator $\D^\A$  is also written as follows in terms of  $\A_\a=\sum_i\A_{\a,i}e_i$　
\bgn{align}\label{ali:3.4}
\D^\A\xi=d\xi+\sum_i[\A_{\a,i}, \xi]e_i.
\end{align}
 Each $e_i\in \TT$ acts on $\psi$ which is denoted by $e_i\cdot\psi\in U^{-n+1}_{\J_\psi}$ and then each element $\xi\otimes e_i$ of 
$\End(E)\otimes(\TT)$ acts on $\psi$ which is denoted by $\xi\otimes (e_i\cdot\psi).$ 
Then the action on $\psi$ gives the map 
$\End(E)\otimes(\TT)\to \End(E)\otimes U^{-n+1}_{\J_\psi}$. 
Thus we obtain $\D^\A\xi \cdot\psi \in \End(E)\otimes U^{-n+1}_{\J_\psi}$ for $\xi\in \End(E)$ and a generalized connection $\D^\A.$
Thus 
 an operator $d^\A$ is defined by  
 \bgn{align}
 d^\A: \Gam(\End&(E))\to\Gam(\End(E)\otimes U^{-n+1}_{\J_\psi})\\
& \xi \mapsto \D^\A\xi\cdot\psi 
 \end{align}
We also extend $d^\A$ to be an operator 
$$d^\A:\Gam(\End(E)\otimes U^{-n+1}_{\J_\psi})\to \Gam(\End(E)\otimes( U^{-n}_{\J_\psi}\oplus U^{-n+2}_{\J_\psi}))$$
by setting :
$$
d^\A(\xi\otimes e_i\cdot\psi)=(\D^\A\xi)\cdot e_i\cdot\psi+\xi\otimes d(e_i\cdot\psi)
$$
Let $a=\sum_i a_i e_i$ be a section of $\End(E)\otimes(\TT)$, where 
$a_i\in \End(E)$ and $e_i\in \TT$.
Then the extended operator $d^\A$ is written in the following form:
\bgn{align}\label{dDA(acdotpsi)}
d^\A (a\cdot\psi)=&\sum_i d^\A(a_ie_i\cdot\psi)=\sum_i (d^\A a_i)e_i\cdot\psi+
a_id(e_i\cdot\psi)\notag\\
=&\sum_{i,j}da_i\cdot e_i\cdot \psi+a_id(e_i\cdot\psi)+[\A_{\a, j}, \, a_i]e_j\cdot e_i\cdot\psi\notag\\
=&d(a \cdot\psi)+[\A\, \cdot a]\cdot\psi
\end{align}
where we are using the following notation: 
$[\A\, \cdot a]:=\sum_{i,j}[A_{\a, j}, \, a_i]e_j\cdot e_i.$
By using the local trivialization and the decomposition $\A_\a=A_\a+V_\a$, the operator $d^\A$ is described as the following:
$$
d^\A (a\cdot\psi) =(d_\a+A_\a+ V_\a) a\cdot\psi =d(a\cdot\psi)+ [V_\a\cdot a]\cdot\psi+ [A_\a\cdot a]\cdot\psi,
$$
where $[\,\cdot\,]$ is the product of 
$\End(E)\otimes \CL$ which is defined as 
$$[V_\a\cdot a]=\sum_{i,j}[V_{\a, i}, a_j]e_i\cdot e_j$$
$$[A_\a\cdot a]=\sum_{i,j}[A_{\a, i}, a_j]e_i\cdot e_j$$

Note that  our new product $[\, \cdot\, ]$ is the combinations of 
 the bracket of 
Lie algebra $\End(E)$ and the Clifford multiplications of Clifford algebra
$\CL$, that is, 
$$
[(A\otimes s)\cdot (A'\otimes s')] =[A, A']\otimes s\cdot s'\in \End(E)\otimes\CL
$$
for $ A \otimes e, A'\otimes e'$, where $A, A'\in \End(E), 
s, s'\in \CL$.
\subsection{Curvature of generalized connections}
Let $\D^\A$ be a generalized Hermitian connection of an Hermitian vector bundle $(E,h)$ over a manifold $M$ which consists an ordinary Hermitian connection $D^A$ and 
a section $V\in \Gam(u(E)\otimes T_M)$. 
Then the ordinary curvature $F_A$ is a section $\End(E)$ valued $2$-from which acts on $\psi$ by the spin
representation to obtain $F_A\cdot\psi \in \End(E)\otimes_\R (U^{-n}\oplus U^{-n+2}).$
The ordinary connection is extended to be an operator as before
$$
d^A: \End(E)\otimes_\R U^{-n+1}\to \End(E)\otimes_\R (U^{-n}\oplus U^{-n+2}).
$$
(For simplicity, we denote by $U^p$ the eigenspace $U^p_{\J_\psi}.$)
By applying $d^A$ to $V\cdot\psi \in \End(E)\otimes_\R U^{-n+1}$, we have $d^A(V\cdot\psi)\in \End(E)\otimes_\R (U^{-n}\oplus U^{-n+2}).$ 
By the spin representation, $[V\cdot V]\in u(E)\otimes_\R \w^2 T_M$ acts on $\psi$ to obtain 
$[V\cdot V]\cdot\psi \in \End(E)\otimes_\R  (U^{-n}\oplus U^{-n+2}).$
Then we have the definition of curvature of generalized connection $\D^\A$:
\bgn{definition}\label{def:Fa(psi)}[curvature of generalized connections]
 $\F_\A(\psi)$ of a generalized connection $\D^\A$ is defined by
\bgn{equation}{\F_\A}(\psi):=F_A\cdot\psi+ d^A(V\cdot\psi)+\frac12[V\, \cdot V]\cdot\psi
\end{equation}
\end{definition} 
$\F(\psi)$ is a globally defined section of $\End(E)\otimes (U^{-n}\oplus U^{-n+2})$ which is called 
the curvature of $\D^\A$.
\bgn{remark}
The ordinary curvature $F_A$ of a connection $D^A$ is defined to be the composition $d^A\circ d^A.$
However $\F_\A(\psi)$ is different from the composition 
$d^\A\circ d^\A$ which is not a tensor but a differential operator.
In fact, we have
\bgn{align}
d^\A\circ d^\A=&(d+\A)\circ (d+\A)=(d^A+ V)\circ (d^A+V)\\
=&F_A+ d^A\circ V + V\circ d^A+ V\circ V
\end{align} 
Then we have 
$$
(V\circ V)\cdot\psi =\frac 12[V\, \cdot V]\cdot\psi
$$
However, for $f\in C^\infty(M)$, $s\in \Gam(E)$, we have 
$$
(d^\A\circ d^\A)(fs\otimes\psi)=f(d^\A\circ d^\A)(s\otimes\psi)
-2\lan df, Vs\ran_{\tt}\cdot\psi
$$
Thus $(d^\A\circ d^\A)$ is not a tensor but a differential operator.
\end{remark}
$\End(E)\otimes(U^{-n}\oplus U^{-n+2})$ is decomposed as 
$$\End(E)\otimes(U^{-n}\oplus U^{-n+2})=(\End(E)\otimes U^{-n})\oplus (\End(E)\otimes U^{-n+2})$$
We denote by $\pi_{U^{-n}}$ the projection from $\End(E)\otimes(U^{-n}\oplus U^{-n+2})$ to the component $\End(E)\otimes U^{-n}$.
 The line bundle $U^{-n}$ becomes the trivial complex line bundle by using the basis $\psi\in U^{-n}$.
Then $\End(E)\otimes_\C U^{-n}$ is identified with $\End(E).$ 
Then it follows from (\ref{eq:End(E)}) that we have 
 \bgn{align}
\End(E)
\otimes U^{-n}
=&u(E)\oplus \Herm(E,h)
\end{align}
We define $\pi_{\Herm}$ to be the projection to the component $\Herm(E,h)$
and we denote by $\pi^{\Herm}_{U^{-n}}$ the composition $\pi^{\Herm}\circ\pi_{U^{-n}}$.
Then we define $\K_\A(\psi)$ by 
$$\K_\A(\psi):=\pi_{U^{-n}}^{\Herm}\F_\A(\psi)\in \Herm(E,h)$$
\bgn{definition}\label{Einstein-Hermitian condition}[Einstein-Hermitian condition]
A generalized Hermitian connection $\D^\A$ is {\it Einstein-Hermitian} if $\D^\A$
satisfies the following: 
$$
\K_\A(\psi)=\lam \,\id_E, \quad \text{\rm for a real constant }\lam
$$
\end{definition}
\bgn{remark}
If $\D^\A$ is an ordinary Hermitian connection $D^A$ over a K\"ahler manifold with
a K\"ahler form $\ome$, then 
the Einstein-Hermitian condition in Definition \ref{Einstein-Hermitian condition} coincides with the ordinary Einstein-Hermitian condition for $\psi=e^{\frac{-\sqrt{-1}}2\ome}.$
In fact, the ordinary connection $D^A$ is Einstein-Hermitian connection if its curvature $F_A$ satisfies
$$\sqrt{-1}\Lam_\ome F_A=+\lam \id_E$$ or equivalently 
$$
\frac{\sqrt{-1}n F_A\w\ome^{n-1}}{\ome^n}= +\lam \id_E,
$$
where $\lam$ is a real constant.
The projection $\pi_{U^{-n}}$ of $F_A$ is given by 
$$
\pi_{U^{-n}}F_A\cdot\psi =\frac{\lan F_A\cdot\psi, \, \ol\psi\ran_s}{\lan \psi, \,\ol\psi\ran_s}\psi
$$
since $\psi=e^{-\frac{\sqrt{-1}}2 \ome}$, we have 
$$
\frac{\lan F_A\cdot\psi, \, \ol\psi\ran_s}{\lan \psi, \,\ol\psi\ran_s}=\frac{F_A\w (-\sqrt{-1}\ome)^{n-1}}{(n-1)!}\frac{n!}{(-\sqrt{-1}\ome)^n}=\sqrt{-1} n\frac{F_A\w\ome^{n-1}}{\ome^n}
$$
Since $F_A$ is in $u(E)$, it follows that 
$$
\sqrt{-1} n\frac{F_A\w\ome^{n-1}}{\ome^n}\in \Herm (E, h).
$$
Thus under the identification $U^{-n}\cong \C\psi$, we have 
$$
\pi^{\Herm}_{U^{-n}}F_A\cdot\psi =\sqrt{-1} n\frac{F_A\w\ome^{n-1}}{\ome^n}\psi
$$
Hence $\pi^{\Herm}_{U^{-n}}F_A\cdot\psi=\lam\, \id _E$ is equivalent to 
$\sqrt{-1}\Lam_\ome F_A=+\lam \id_E$. 
\end{remark}
The unitary gauge group $U(E, h)$ acts on $\F_\A(\psi)$ by the adjoint action. Thus the Einstein-Hermitian condition is invariant under the action of the unitary gauge group.
Further our Einstein-Hermitian condition behaves nicely 
for the action of $b$-fields.
A real $d$-closed $2$-form $b$ acts on $\psi$ by 
$e^b\cdot\psi$ which is also a $d$-closed, nondegenerate, pure, spinor.
\bgn{definition}[$b$-field action of generalized connections]
Let $\D^\A=d+\A$ be a generalized connection of vector bundle $E$ over $(M, \psi).$  
A $d$-closed $2$-form $b$ acts on a generalized connection $\D^\A$ by 
$e^b\D^\A e^{-b}=d+e^b\A e^{-b}$, 
Then $e^b\A e^{-b}$ is given by 
$$
e^b\A e^{-b}=A+e^b Ve^{-b}=A+ad_b V +V,
$$
where $A+ ad_b V\in \End(E)\otimes T^*_M$ and 
$d+A+ \ad_b V$ is an ordinary connection of $E$ and 
$V\in \End(E)\otimes T_M$.
Note that $ad_b V$ is given by 
$$
\ad_bV =\sum_i V_i \otimes[b, v_i],
$$
for $V=\sum_i V_i \otimes v_i$, where $V_i\in u(E)$ and $v_i\in T_M$.
\end {definition}

\bgn{proposition}\label{prop: b-field action}
Let $\Ad_{e^b}\D^\A=d+\Ad_{e^b}\A$ be a generalized connection which is given by 
the $b$-field action on $\D^\A$. 
Then the curvature $\F_{\Ad_{e^b}\A}(\psi)$ is given by 
$$
\F_{\Ad_{e^b}\A}(\psi)=e^b\F_\A(e^{-b}\psi)
$$
 \end{proposition}
 \bgn{proof}
 The $b$-field action is given by 
  $\Ad_{e^b}\D^\A=d+A+\ad_bV + V,$ where $d+A+\ad_bV$ is a connection and its curvature is given by 
 $$
 F_A+d^A(\ad_b V)+\frac12[\ad_b V, \, \ad_b V]
 $$
 Then we have 
 \bgn{align}
\F_{\Ad_{e^b}\A}(\psi)= & (F_A+d^A(\ad_b V)+\frac12[\ad_b V, \, \ad_b V])\cdot\psi\\
+&d^A(V\cdot\psi)+[\ad_bV \cdot V]\cdot\psi+\frac12 [V\cdot V]\cdot\psi\\
\end{align}
Applying $\Ad_{e^b}V=V+\ad_b V$ and $e^bF_A e^{-b}=F_A$ and 
$[e^b Ve^{-b}\,\cdot\, e^b Ve^{-b}]=e^b[V\, \cdot\, V]e^{-b}$,  we have 

\bgn{align}
\F_{\Ad_{e^b}\A}(\psi)=&F_A\cdot\psi+ d^A(e^b Ve^{-b}\cdot\psi)+
\frac12[e^b Ve^{-b}\,\cdot\, e^b Ve^{-b}]\cdot\psi\\
=&e^bF_A\cdot e^{-b}\cdot\psi+ e^b d^A(V\cdot e^{-b}\psi)+\frac12e^b[V\cdot V]\cdot e^{-b}\cdot\psi\\
=&e^b\F_\A(e^{-b}\psi)
\end{align}
 \end{proof}
 Thus our Einstein-Hermitian condition is equivalent under  the action of $b$-field.
 \bgn{proposition}\label{prop: invariance of b-field action}
 Let $\D^\A$ be a generalized Hermitian connection of $E$ over $(M, \psi)$. Then  $\D^\A$ is Einstein-Hermitian over $(M,\psi)$ if and only if
 $\Ad_{e^b}\D^\A$ is an Einstein-Hermitian generalized connection 
 of $E$ over $(M, e^{b}\psi)$.
 \end{proposition}
 \bgn{proof}
 We denote by $\pi_{U^{-n}_{e^b\psi}}$ the projection to the component 
 $U^{-n}_{e^b\psi}:=e^b\cdot U^{-n}_{\J_\psi}.$ 
 From Proposition \ref{prop: b-field action}, we have 
 \bgn{align}
 \pi^{\Herm}_{U^{-n}_{e^b\psi}}\F_{\Ad_{e^b}\A}(e^b \psi)=&\pi^{\Herm}
 \frac{\lan \F_{\Ad_{e^b}\A}(e^b \psi), e^b\cdot\ol\psi\ran_s}{\lan e^b\psi, \,e^b\ol\psi\ran_s}\\
=&\pi^{\Herm}
 \frac{\lan e^b\F_{\A}(\psi), e^b\cdot\ol\psi\ran_s}{\lan e^b\psi, \,e^b\ol\psi\ran_s}
 \end{align}
 Since $\lan e^b\psi, \,e^b\ol\psi\ran_s=\lan \psi, \, \ol\psi\ran_s$, we have 
 $$
 \pi^{\Herm}_{U^{-n}_{e^b\psi}}\F_{\Ad_{e^b\A}}(e^b \psi)=
 \pi^{\Herm}
 \frac{\lan \F_{\A}(\psi), \,\cdot\ol\psi\ran_s}{\lan \psi, \,\ol\psi\ran_s}
=\pi^{\Herm}_{U^{-n}_{\J_\psi}}\F_\A(\psi).
 $$
Hence $\K_\A(\psi)$ is invariant under the action of $b$-fields.
Thus $\pi^{\Herm}_{U^{-n}_{e^b\psi}}\F_{\Ad_{e^b\A}}(e^b \psi)=\lam\id_E$
 if and only if $\pi^{\Herm}_{U^{-n}_{\J_\psi}}\F_\A(\psi)=\lam\id_E$. 
  Thus the result follows.
  \end{proof}
\subsection{The first Chern class of $E$ and $\tr\F_\A(\psi)$}
\bgn{theorem}\label{th: Chern class}
$\frac{-1}{2\pi\sqrt{-1}}\tr\F_\A(\psi)$ is a $d$-closed differential form on $M$ which is a representative of the class $[c_1(\det E)]\cup [\psi]\in H^\bullet(M)$.
\end{theorem}
\bgn{proof}As in Definition \ref{def:Fa(psi)}, $\tr\F_\A(\psi)$ is given by 
$$
\tr\F_\A(\psi):=\tr F_A\cdot\psi+ \tr d^A(V\cdot\psi)+\tr\frac12[V\, \cdot V]\cdot\psi
$$
We have  $\tr[V\cdot V]\cdot\psi=0$. 
By using local trivializations of $E,$ we have 
$$
\tr d^A(V\cdot\psi)=\tr \(d(V_\a\cdot\psi)+[A_\a\cdot V_\a]\)\cdot\psi 
=d (\tr V)\cdot\psi.
$$
Since $\tr\, V\cdot \psi$ is a globally defined form on $M,$ we have 
$$
\frac{-1}{2\pi\sqrt{-1}}[\tr\F_\A(\psi)]=\frac{-1}{2\pi\sqrt{-1}}[\tr F_A\cdot\psi]=[c_1(\det (E)]\cup[\psi].
$$
\end{proof}
\bgn{proposition}
Let $\D^\A$ be a generalized Einstein-Hermitian connection which satisfies $\pi^{\Herm}_{U^{-n}}\F_\A(\psi)=\lam \id_E$. 
Then $\lam$ is given in terms of the first Chern class $c_1(E)$ of $E$ 
and the class $[\psi]$ by 
$$
\frac{-\lam r}{2\pi\sqrt{-1}}\int_M\lan \psi, \, \ol\psi\ran_s =
\int_M \lan c_1(E)\w\psi ,\,\ol\psi\ran_s
$$
\end{proposition}
\bgn{proof}
It follows from Theorem \ref{th: Chern class} that we have 
$$
\frac{-1}{2\pi\sqrt{-1}}\int_M\tr\lan \F_\A(\psi), \, \ol\psi\ran_s =
\int_M\lan c_1(E)\w\psi, \, \ol\psi\ran_s.
$$
The $U^{-n}$-component of $\F_\A(\psi)$ is written as 
$$
\F_\A(\psi)=+\lam \id_E\psi+ \xi\psi,
$$
where $\lam\in \R$ and $\xi\in u(E)$. 
Then we have 
\bgn{align*}
\frac{-1}{2\pi\sqrt{-1}}\tr\lan \F_\A(\psi), \, \ol\psi\ran_s
=&\frac{-1}{2\pi\sqrt{-1}}\tr \lan +\lam \psi, \, \ol\psi\ran_s+ 
\frac{-1}{2\pi\sqrt{-1}}\tr\lan \xi \psi, \,\ol\psi\ran_s\\
=&\frac{+\lam r}{2\pi\sqrt{-1}}\lan \psi, \, \ol\psi\ran_s
+\frac{-1}{2\pi\sqrt{-1}}(\tr\xi)\lan \psi, \, \ol\psi\ran_s
\end{align*}
We have 
$$
\frac{\lan c_1(E)\w\psi, \, \ol\psi\ran_s}{\lan \psi, \, \ol\psi\ran_s}
=\frac{c_1(E)\w(\sqrt{-1}\ome)^{n-1}}{(n-1)!}\frac{n!}{(\sqrt{-1}\ome)^n}
=\frac{c_1(E)\w \ome^{n-1}}{\ome^n}\frac{n}{\sqrt{-1}}\in \sqrt{-1}\R
$$
Since $\xi \in u(E)$, $\tr\xi$ is pure imaginary. 
Thus we have the result.
\end{proof}
\section{Generalized holomorphic vector bundles}
\label{Generalized holomorphic vector bundles}
Let $E$ be a complex vector bundle over a \complex manifold $(M, \J)$. 
{\it A generalized holomorphic structure} of $E$ is a differential operator 
$$
\ol\pa^E_\J : E\to E\otimes \ol \L_\J,
$$
which satisfies 
$$
\ol\pa^E_\J (f s) = s\otimes(\ol\pa_\J f)   + f(\ol\pa^E_\J s), \qquad 
\ol\pa^E_\J\circ\ol\pa^E_\J=0, \qquad \text{\rm for   }f\in C^\infty(M), s\in E.
$$
A complex vector bundle equipped with a generalized holomorphic structure is called as
{\it a generalized holomorphic vector bundle}
\cite{Gua_2010}.
Then a generalized holomorphic structure $\ol\pa^E_\J$ is extended to be  
an operator $\ol\pa^E_\J : E \otimes \w^i\ol \L_\J \to E\otimes \w^{i+1}\ol\L_\J$ by 
$$
\ol\pa^E_\J (s \a) =(\ol\pa^E_\J s )\w\a+ s (\ol\pa_\J \a),\qquad \text{\rm for   }
s\in E, \a\in \w^i\ol\L_\J.
$$
Then we obtain an elliptic complex which is the Lie algebroid complex.
$$
0\to E\to E\otimes \ol\L_\J \to  E\otimes \w^2\ol\L_\J\to \cdots \to  E\otimes \w^i\ol\L_\J \to 
E\otimes \w^{i+1}\ol\L_\J \to \cdots\to 0
$$
We denote by 
$C^{0,\bullet}=\{C^{0,i}\}$ the Lie algebroid complex  
$(E\otimes \w^i\ol{\L_\J}, \ol\pa^E_\J)$. 
\subsection{The canonical generalized connection}
Let $(E, h)$ be an Hermitian vector bundle over a generalized complex manifold $(M,\J)$. 
We assume that $E$ admits a generalized holomorphic structure $\ol\pa^E_\J$.
Then there is a unique generalized Hermitian connection $\D_h$ such that 
$$
\D_h=\D^{1,0}+\ol\pa_\J^E,
$$
 where $\D^{1,0}: E\to E\otimes\L_\J$ denotes the $(1,0)$-component of $\D_h$ with respect to $\J$. 
 We call $\D_h$ as the canonical generalized connection of Hermitian vector bundle over a
 \complex manifold. 
 This is an analogue of the canonical connection of Hermitian vector bundle over a
 complex manifold. 
 In fact, $\D_h$ is determined by 
 $$
 \pa_\J h(s_1, s_2) =h(\D^{1,0}s_1, s_2)+h(s_1, \ol\pa^E_\J s_2), \qquad s_1, s_2\in E.
 $$
 We denote by $\K_h(\psi)$ the mean curvature of the generalized connection $\D_h.$
\bgn{definition}\label{Einstein-Hermitian metric}[Einstein-Hermitian metric]
If the canonical generalized connection $\D_h$ satisfies the Einstein-Hermitian condition,
$$
\K_h(\psi)=\lam \id_E\psi,
$$
then $h$ is called {\it an Einstein-Hermitian metric.}
\end{definition}


\subsection{Weak generalized holomorphic subbundles}
\label{Weak generalized holomorphic subbundles}
Let $(E, \ol\pa^E_\J)$ be a generalized holomorphic vector bundle over a generalized complex manifold $(M, \J).$
Let $F$ is a subbundle of $E$. 
If $F$ admits a generalized holomorphic structure $\ol\pa_\J^{F}$  such that 
$j\circ \ol\pa_\J^{F}=\ol\pa_\J^E\circ j$, where 
$j: F\to E$ denotes the inclusion of $F$ into $E$, 
then $F$ is called a generalized holomorphic subbundle, that is, 
$$
\xymatrix{&{F}\ar@{->}^{\ol\pa_\J^{F}\,\,\,\,}[r]\ar@{->}[d]&{F}\otimes\ol\L_\J\ar@{->}[d]\\
&{E}\ar@{->}^{\ol\pa^{E}_\J\,\,\,\,}[r]&{E}\otimes\ol\L_\J
}
$$
In order to define the notion of stability of $E,$ we need "a subbundle with singularities" which is already introduced in Uhlenbeck-Yau's paper as "a weak holomorphic subbundle". 
We shall define a weak generalized holomorphic subbundle: 
\bgn{definition}
Let $E$ be a complex vector bundle with an Hermitian metric $h.$
An element $\pi$ of $L_1^2(\End(E))$ is called {\it a weak generalized holomorphic subbundle} of $E$ if the followings hold in $L^1(\End(E))$:
\bgn{equation}\label{weakly hol.sub}
\pi^*=\pi=\pi^2, \qquad (\id_E-\pi)\circ \ol\pa^E_\J \pi =0
\end{equation}
\end{definition}
We assume that there is another \complex structure $\J_2$ such that 
$(M, \J, \J_2)$ is a generalized K\"ahler manifold. 
From a view point of bihermitian structure $I_\pm$, a generalized holomorphic vector bundle $(E, \ol\pa_\J^E)$
gives a locally free sheaf with respect to both complex structure $I_+$ and $I_-.$
We also have 
\bgn{proposition}
A weak generalized holomorphic subbundle $\pi$ is a weak holomorphic subbundle with respect to 
both $I_+$ and $I_-.$
\end{proposition}
\bgn{proof}
Let $\pi_T: \TT \to T_M$ be the projection from $\TT$ to $T_M.$
Then $\pi_T$ gives the identification $\ol\L_\J^+\cong T_{I_+}^{1,0}$ and 
$\ol\L_\J^-\cong T_{I_-}^{1,0},$ where $T^{1,0}_{I_\pm}$ denotes the $(1,0)$-component of $T_M^\C$
with respect to $I_\pm,$ respectively. Thus we have  
the lift $\til v_\pm$ of $v\in T_M^\C$ defined by $\pi_T(\til v_\pm )=v$, which is explicitly given by 
\bgn{align}
\til v_\pm =&v^{1,0}_\pm \pm g(v^{1,0}_\pm, \, )+b(v^{1,0}_\pm,\,)\\
 +&v^{0,1}_\pm\pm g(v^{0,1}_\pm, \,)+b(v^{0,1}_\pm, \, ),
\end{align}
where $(I_\pm, g, b)$ denotes the bihermitian structure corresponding to the generalized K\"ahler structure and 
$v_\pm^{1,0}\in T_{I_\pm}^{1,0}$ and $v^{0,1}_{I_\pm}\in T^{0,1}_{I_\pm}$ and 
\bgn{align}
&v^{1,0}_\pm \pm g(v^{1,0}_\pm, \, )+b(v^{1,0}_\pm,\,)\in \ol{\L}_\J^\pm\\
&v^{0,1}_\pm\pm g(v^{0,1}_\pm, \,)+b(v^{0,1}_\pm, \, )\in \L_\J^\pm
\end{align}

We denote by $\ol\pa^\pm$ the ordinary $\ol\pa$-operators with respect to $I_\pm.$
Since $df=\pa_\J f +\ol\pa_\J f$ for a function $f,$
we have 
\bgn{align}
\lan \ol\pa_\J f,\,\, \til v_\pm\ran_{\tt} =&\lan \ol\pa_\J f, \,\(v^{0,1}_\pm\pm g(v^{0,1}_\pm, \,)+b(v^{0,1}_\pm, \, )
\)\ran_{\tt}\\
=&\lan df, \,\(v^{0,1}_\pm\pm g(v^{0,1}_\pm, \,)+b(v^{0,1}_\pm, \, )
\)\ran_{\tt}\\
=&\lan df, \,\, v^{0,1}_\pm\ran_{\tt}\\
=&\ol\pa^{\pm}f(v)
\end{align}

Then we have
$$\ol\pa^{\pm}_v s= \lan\ol\pa_\J^E s,\,\,\til v_\pm\ran_{\tt}$$
where $\til v_\pm \in C^\pm$ is the lift of $v.$ 
Thus $\pi$ satisfies 
$$
(\id_E-\pi)\circ \ol\pa^\pm \circ\pi =0
$$
Hence $\pi$ is a weak holomorphic subbundle with respect to both $I_\pm.$
\end{proof}
The following result ensures that a weak holomorphic subbundle gives rise to a subsheaf.

\bgn{theorem}[Uhlenbeck-Yau]\label{Uhlenbeck-Yau}
A weak holomorphic subbundle $\pi$ of $\E$ represents a coherent subsheaf $\F$ of $\E$.
More precisely,
there is a coherent subsheaf $\F$ of $\E$ and an analytic subset $S\subset M$ such that 
\bgn{enumerate}
\item codim $S\geq 2$
\item $\pi|_{M\bsh S}$ is $C^\infty$ and satisfies (\ref{weakly hol.sub}).
\item $\F|_{M\bsh S}=\pi(\E|_{M\bsh S})\subset \E|_{M\bsh S}$ is locally free, i.e., a holomorphic subbundle.
\end{enumerate}
\end{theorem}
Then from Theorem \ref{Uhlenbeck-Yau}, we have
\bgn{proposition}
We denote by $\E_\pm$  locally free sheaves which are given by a generalized holomorphic vector bundle $(E, \ol\pa_\J^E)$ with respect to $I_\pm, respectively.$
Let $\pi$ be a weak generalized holomorphic subbundle of $(E, \ol\pa_\J^E).$
Then there are coherent subsheaves $\F_{\pm}$ of $\E_\pm$ which satisfy   
\bgn{enumerate}
\item codim $S\geq 2$
\item $\pi|_{M\bsh S}$ is $C^\infty$ and satisfies (\ref{weakly hol.sub}).
\item $\F_{\pm}|_{M\bsh S}=\pi(\E_\pm|_{M\bsh S})\subset \E_\pm$ is locally free, i.e., a holomorphic subbundle.
\end{enumerate}

\end{proposition}
\section{Definition of $\psi$-stability and $\psi$-polystability}
\label{Definition of stability and polystability}
Let $(E, \ol\pa^E_\J)$ be a  generalized holomorphic vector bundle. 
Then the degree of $E$ is given by 
$$
\deg(E):= \int_M i^n\lan c_1(E)\cdot\psi, \,\,\ol\psi\ran_s,
$$
where $c_1(E)$ denotes a $d$-closed $2$-form representing the first Chern class of $E$ which acts on $\psi.$
Then the slope $\mu(E)$ of $E$ is the ratio
$$
\mu(E):=\frac{\deg E}{\rank E}
$$
Let $h$ be an Hermitian metric on $E$. Then we shall define the stability of $E$ by using weak generalized holomorphic subbundles.
A weak generalized holomorphic subbundle $\pi$ gives a generalized holomorphic subbundle $F$ on a complement $M\bsh S$ of codim $2$ 
subset $S$ from Proposition \ref{Uhlenbeck-Yau}. 
Then we define $\rank(\pi)$ to be the rank of the generalized subbundle $F$ on the complement. 
Since $F$ is a generalized holomorphic subbundle with the induced Hermitian metric $h|_F,$ the canonical generalized connection 
$\D^F$ of $(F, \ol\pa^F_\J, h|_F)$ gives the first Chern form $c_1(\pi, h|_F)$ on the complement $M\bs S$. 
Then we define
$$
\deg(\pi):= \int_{M\bs S} i^n\lan c_1(\pi, h|_F)\cdot\psi, \,\,\ol\psi\ran,
$$
and then $$\mu(\pi):=\frac{\deg (\pi)}{\rank (\pi)}.$$
\bgn{remark} The degree of $F$ is well-defined which is given in terms of the second fundamental form (c.f. Section \ref{The second fundamental forms}). 
It turns out that $\deg(\pi)$ is finite which coincides with the one given by the first Chern class of the coherent sheaves 
$\F_\pm.$ 
\end{remark}

\bgn{definition}[Stability]
A generalized holomorphic vector bundle $(E,\ol\pa^E_\J)$ is {\it stable} if and only if
for every weak generalized holomorphic subbundle $\pi $ with $0< \rank\,\pi < \rank E$, we have 
$$
\mu(\pi) <\mu(E).
$$
If we have 
$$
\mu(\pi)\leq\mu(E),
$$
then $E$ is {\it semistable}.
If $E$ is decomposed into the direct sum $\oplus_i E_i$ of generalized holomorphic subbundles with the same slope $\mu(E)$, then 
$E$ is {\it polystable}.
\end{definition}
\section{The second fundamental forms}\label{The second fundamental forms}
Let $E$ be a generalized holomorphic vector bundle with an Hermitian metric $h$ over a generalized complex manifold $(M,\J)$.
We denote by $S$ a generalized holomorphic subbundle of $E$. 
Then we have the short exact sequence
$$
0\to S\to E\to Q\to 0.
$$
The quotient $Q$ is identified with the orthogonal complement $S^\perp$ by using $h$ and we have the decomposition $E=S\oplus S^{\perp}.$
Let $\D^E$ be the canonical generalized Hermitian connection of $E$ with respect to $h.$ 
Then we have the decomposition 
$\D^Es =\D^S s+ H^S(s)$ for $s\in \Gam(S)$
where $\D^S(s)\in S$ and $H^S(s)\in S^\perp.$
We also  denote by $\pa_0$ the operator $\D^{1,0}.$
\bgn{proposition}\label{leths}
Let $h_S$ be the Hermitian metric of $S$ which is the restriction of $h$ to $S.$ 
We denote by $\pi_s$ the orthogonal projection from $E$ to $S$. Then we have\\ 
(1) $\D^S $ is the canonical connection of the generalized holomorphic vector bundle $S$ with respect to $h_S$.\\
(2) $H^S$ is a section $\pa_0 \pi_s$of $\L_\J\otimes \Hom(S, S^\perp)$, where $\pa_0=\D^{1,0}$ acts on $\End(E).$
\end{proposition}
\bgn{proof}
Since $S$ is a generalized holomorphic subbundle, $\D^{0,1}s=\ol\pa^E s\in \Gam(\ol\L_\J\otimes S)$ for $s\in \Gam(S)$. Then (1) follows.
Since $\D^Es=\D^Ss+ H^S s$, it follows that $(H^S s)^{0,1}=0.$
Hence we obtain $H^S\in \L_\J\otimes\Hom (S, S^\perp).$ Then
$H^S$ is given by 
\bgn{align}
H^S=&(1-\pi_s)\circ \D^E\circ \pi_s=(1-\pi_s)\circ \pa_0\circ\pi_s
\end{align}
Since $(\pa_0\pi_s) s=\pa_0 s-\pi_s (\pa_0 s)$ for $s\in \Gam(S),$
 we see that 
$$
H^S s- \pa_0\pi_s s =(1-\pi_s)\circ \pa_0\circ\pi_s s -(\pa_0\pi_s) s
=\pa_0 s-\pi_s (\pa_0 s)-\pa_0 s+\pi_s(\pa_0 s)
=0
$$
Thus $H^S=\pa_0\pi_s.$
\end{proof}
We also define $\D^{S^\perp}$ and $H^{S^\perp}$ by 
$$
\D^E s_Q =H^{S^\perp}s_Q+\D^{S^\perp}s_Q , 
$$
for $s_Q\in \Gam(S^\perp)$,
where  $H^{S^\perp}s_Q\in S,$ $\D^{S^\perp}s_Q\in S^\perp.$
Then we also have 
\bgn{proposition}\label{HSperp}
(1) $H^{S^\perp}$ is a section $-\ol\pa_\J \pi$$\in \ol \L_\J\otimes\Hom(S^\perp, S)$. 
\\
(2) $H^S$ and $H^{S^\perp}$ satisfy the following
$$
h(H^S s, \, s')+ h(s, \,\,H^{S^\perp}s')=0
$$
\end{proposition}
\bgn{proof}
The results follows from the similar way as before.
\end{proof}

Let $\{e_i\}$ be a local basis of $\L_\J$. Then $H^S$ and $H^{S^\perp}$ are written as follows 
\bgn{equation}\label{eq:HS}
H^S=\sum_i H^S_i  e_i, \qquad H^{S^\perp}=-\sum_i (H^S)^* \ol e_i,
\end{equation}
where $(H^S)^*$ denotes the adjoint of $H^S$ with respect to $h$.
\section{From Einstein-Hermitian metrics to $\psi$-stability}\label{From Einstein-Hermitian metrics to stability}
At first, we assume that $\psi=e^{-\sqrt{-1}\ome}$
and then we shall show that Einstein-Hermitian generalized holomorphic vector bundle is $\psi$-polystable. 
Then we shall show the result in general cases of $\psi=e^{b-\sqrt{-1}\ome}$ by using the invariance under the action of $b$-field.
Let $(E, \ol\pa_\J^E)$ be a generalized holomorphic bundle with an Hermitian metric $h$
and $S$ a generalized holomorphic subbundle of $(E, \ol\pa_\J^E).$
 Let $\D^E$ be 
 the canonical generalized connection of $(E,\ol\pa_\J^E, h)$ over a generalized K\"ahler manifold $(M, \J, \J_\psi)$.
A generalized connection $\D^E$ is written as 
$$
\D^E=d+\A=d+A+V,
$$
where $A\in T^*_M\otimes\End(E)$ is an ordinary connection form and $V\in T_M\otimes\End(E).$
We identify the quotient bundle $Q$ with the orthogonal complement $S^\perp$.
Then the decomposition $E=S\oplus S^\perp$ gives the following decomposition of $\A$: 
$$
\A=
\bgn{pmatrix}\A_{ss}&\A_{sq}\\ \A_{qs}&\A_{qq}
\end{pmatrix}
$$
where $\A_{ss}\in\Hom(S, S),\, \A_{sq}\in\Hom(Q,S), \A_{qs}\in\Hom(S,Q)$ and $ \A_{qq}\in \Hom(Q,Q).$
The connection form $A$ and the section $V$ are also decomposed 
\bgn{equation}\label{matrix decomposition A,V}
A=\bgn{pmatrix}A_{ss}&A_{sq}\\ A_{qs}&A_{qq}
\end{pmatrix}\qquad\qquad
V=\bgn{pmatrix}V_{ss}&V_{sq}\\ V_{qs}&V_{qq},
\end{pmatrix}
\end{equation}
Then the second fundamental forms are given by
\bgn{equation}\label{eq:Hs}
H^S=\A_{qs}=A_{qs}+V_{qs}, \qquad H^{S^\perp}=\A_{sq}=A_{sq}+V_{sq}
\end{equation}
 Let $\F_\A(\psi)=F_A\cdot\psi+ d^A(V\cdot\psi)+ V\cdot V\cdot\psi$ be the generalized curvature of $\D^E$ of $(E, h).$
If we take $\psi=e^{\frac{\ome}{\sqrt{-1}}},$
then we have 

$$
\K_\A(\psi) =\pi_{U^{-n}}^{\Herm}\F_\A(\psi)=\pi_{U^{-n}}^{\Herm}\(F_A\cdot\psi + \frac12[V\cdot V]\cdot\psi\).
$$
(Note that $\pi^{\Herm}_{U^{-n}}d^\A(V\cdot\psi)=0$.)

Let $\{e_i\}$ be a local basis of $\L_\J=\L_\J^+\oplus\L_\J^-$. When we take the decomposition $e_i=e_i^++e_i^-$, then we assume that $\{e_i^\pm\}$ satisfies 
\bgn{equation}\label{laneipm}
\lan e_i^\pm, \,\, \ol{e_i^\pm}\ran_s =\pm\del_{i,j}
\end{equation}
Then we have
\bgn{lemma}\label{lem:key lemma 1}
\bgn{align}
\lan \ol e_i \cdot  e_j \cdot \psi, \,\,\ol\psi\ran_s =&-2\del_{i,j}\lan \psi, \,\,\ol\psi\ran_s\\
\lan e_i \cdot \ol e_j \cdot \psi, \,\,\ol\psi\ran_s =&+2\del_{i,j}\lan \psi, \,\,\ol\psi\ran_s
\end{align}
\end{lemma}
\bgn{proof}
Since $e_i\cdot\psi=e_i^-\cdot\psi$ and $\ol {e_i^-}\cdot\psi=0$ and 
$(e_i^- \cdot \ol {e_j^-}+\ol {e_i^-} \cdot  e_j^-)=2\lan e_i^-, \,\, \ol {e_j^-}\ran_{\tt},$\\
 we have
\bgn{align}
\lan \ol e_i \cdot  e_j \cdot \psi, \,\,\ol\psi\ran_s =
&\lan \ol {e_i^-} \cdot  e_j^- \cdot \psi, \,\,\ol\psi\ran_s \\
=&\lan (e_i^- \cdot \ol {e_j^-}+\ol {e_i^-} \cdot  e_j^-) \cdot \psi, \,\,\ol\psi\ran_s\\
=&+ 2\lan e_i^-, \,\, \ol {e_j^-}\ran_{\tt}\lan \psi,\,\,\ol\psi\ran_s\\
=&- 2\lan \psi, \,\,\ol\psi\ran_s
\end{align}
\bgn{align}
\lan  e_i \cdot \ol e_j \cdot \psi, \,\,\ol\psi\ran_s =
&\lan  {e_i^+} \cdot \ol{e_j^+} \cdot \psi, \,\,\ol\psi\ran_s \\
=&\lan (e_i^+ \cdot \ol {e_j^+}+\ol {e_i^+} \cdot  e_j^+) \cdot \psi, \,\,\ol\psi\ran_s\\
=&+ 2\lan e_i^+, \,\, \ol {e_j^+}\ran_{\tt}\lan \psi,\,\,\ol\psi\ran_s\\
=&+ 2\lan \psi, \,\,\ol\psi\ran_s
\end{align}
\end{proof}
\bgn{lemma}\label{kihon lemma}
Let $\a$ be a differential form which acts on $\psi$. 
Then we have
$$\pi_{U^{-n}}(\a\cdot\psi) =\frac{\lan \a\cdot\psi, \,\,\ol\psi\ran_s}{\lan \psi, \,\,\ol\psi\ran_s}\psi
$$
\end{lemma}
\bgn{proof}
Let $\pi_{U^p}$ denotes the projection to the component
$U^p.$
Then we have 
$$
\lan \pi_{U^p}(\a\cdot\psi), \,\, \ol\psi\ran_s=0
$$
for $p\neq -n.$
We denote by $f\psi$ the component $\pi_{U^{-n}}(\a\cdot\psi),$ where $f$ is a function.
Then we have 
$$\frac{\lan \a\cdot\psi, \,\,\ol\psi\ran_s}{\lan \psi, \,\,\ol\psi\ran_s}\psi
=\frac{\lan f\psi, \,\,\ol\psi\ran_s}{\lan \psi, \,\,\ol\psi\ran_s}\psi
=f\psi=\pi_{U^{-n}}(\a\cdot\psi).
$$
\end{proof}
The composition $H^S\cdot H^{S^\perp}$ between $H^S\in\L_\J\otimes \Hom(S, S^\perp)$ and $H^{S^\perp}\in \ol{L_\J}\otimes\Hom(S^\perp, S)$ is a section of $\L_\J\w\ol{\L_\J}\otimes\End(S)$.
Then $(H^S\cdot H^{S^\perp})\cdot\psi$ is a section of $(U^{-n}\oplus U^{-n+1})\otimes\End(S).$
Then the trace of the projection $\pi^{\Herm}_{U^{-n}}(H^{S^\perp}\cdot H^S)\cdot\psi$ is given by the following: 
\bgn{lemma}\label{lem: positive}
 \bgn{align}\tr\,\pi_{U^{-n}}^{\Herm }(H^{S^\perp}\cdot H^S)\cdot\psi=&+ \|H^S\|^2\psi\\
 \tr\,\pi_{U^{-n}}^{\Herm }(H^{S}\cdot H^{S^\perp})\cdot\psi=&-  \|H^S\|^2\psi
 \end{align}
\end{lemma}
\bgn{proof}
By using (\ref{eq:HS}) and (\ref{eq:Hs}), we have
\bgn{align}
H^{S^\perp}\cdot H^{S}=\A_{sq}\cdot\A_{qs}=-\sum_{i,j}(H^S_i)^*(H^S_j)\ol e_i  e_j\\
H^{S}\cdot H^{S^\perp}=\A_{qs}\cdot\A_{sq}=-\sum_{i,j}(H^S_i)(H^S_j)^* e_i  \ol e_j
\end{align}
From Lemma \ref{lem:key lemma 1}, we have 
\bgn{align}
\tr\lan(H^{S^\perp}\cdot H^S)\cdot\psi, \, \ol\psi\ran_s
=&-\sum_{i,j}\tr((H^S_i)^*(H^S_j))
\lan \ol e_i \cdot e_j \cdot \psi, \,\,\ol\psi\ran_s \\
=&+2\sum_i\|H^S_i\|^2\lan\psi, \,\,\ol\psi\ran_s\\
=& \|H^S\|^2\lan \psi, \,\,\ol\psi\ran_s,
\end{align}
Note that it follows from (\ref{laneipm}) that $e_i=e_i^++e_i^-$ and $\|e_i\|^2=2$ and $\|H^S\|^2=+2\sum_i\|H^S_i\|^2.$
From Lemma \ref{kihon lemma}, we have 

\bgn{align}
\tr\,\,\pi_{U^{-n}}^{\Herm }(H^{S^\perp}\cdot H^S)\cdot\psi=& \|H^S\|^2\psi.
\end{align}
We also have
\bgn{align}
\tr\lan\pi_{U^{-n}}^{\Herm }(H^{S}\cdot H^{S^\perp})\cdot\psi, \,\,\ol\psi\ran_s=&-\sum_{i,j}\tr((H^S_i)(H^S_j)^*)
\lan e_i \cdot \ol e_j \cdot \psi, \,\,\ol\psi\ran_s \\
=&-  \|H^S\|^2\psi.
\end{align}
Thus we have 
\bgn{align}
\tr\,\,\pi_{U^{-n}}^{\Herm }(H^{S}\cdot H^{S^\perp})\cdot\psi=&- \|H^S\|^2 \psi.
\end{align}
\end{proof}
\bgn{lemma}\label{lem:hontouni key}
Let $P$ be a skew Hermitian endmorphism. Then we have 
$$
\pi^{\Herm}_{U^{-n}}P(\t \w v)\cdot \psi =0,
$$
where $v\in T_M$ and $\t\in T_M^*$ are real. 
\end{lemma}
\bgn{proof}
This is a point-wise calculation. 
Since $\psi=e^{-\sqrt{-1}\ome}$, we have 
$v-\sqrt{-1}i_v\ome\in \ol{\L_{\J_\psi}}$ and $\t-\sqrt{-1}\ome^{-1}(\t)\in \L_{\J_{\psi}}.$
Then we have 
\bgn{align}
4\pi_{U^{-n}}(\t\w v)\cdot\psi =&
\pi_{U^{-n}}(\t-\sqrt{-1}\ome^{-1}(\t))\cdot(v-\sqrt{-1}i_v\ome)\cdot\psi\\
=&2\lan(\t-\sqrt{-1}\ome^{-1}(\t)),\,\, (v-\sqrt{-1}i_v\ome)\ran_{\tt}\psi\\
=&2\t(v)\psi,
\end{align}
where $\t(v)$ denotes the coupling $\t\in T^*_M$ and $v\in T_M.$
Since $\t(v)$ is real, $ P \t(v) $ is skew-Hermitian. 
Thus we have 
$\pi^{\Herm}_{U^{-n}}P(\t \w v)\cdot\psi=0.$
\end{proof}
Let $\F^E_\A(\psi)$ be the curvature of $E$ and $\K^E_\A(\psi)$ the mean curvature of $E$.
Since $\frac12[V\cdot V]\cdot\psi =V\cdot V\cdot\psi$ and $\frac12[A\cdot A]\cdot\psi=A\cdot A\cdot\psi$ for $V\in T_M\otimes\End(E)$ and $A\in T_M^*\otimes\End(E),$
then the component of
$\F^E_\A(\psi)$ of End $(S)$ is given by 
\bgn{align}
\pi_s\circ \F_\A(\psi)\circ \pi_s
=&\pi_s\circ (F_A+ V\cdot V)\cdot\psi\circ \pi_s  \\
=&F^S\cdot\psi+V_{ss}\cdot V_{ss}\cdot\psi + A_{sq}\cdot A_{qs}\cdot\psi + V_{sq}\cdot V_{qs}\cdot\psi,\\
\end{align}
where $\pi_s$ is the orthogonal projection to $S$ and $i_S: S\to E$ denotes the inclusion and
$F^S$ denotes the ordinary curvature of Hermitian connection 
 $A_{ss}$ of the subbundle $S$.
It must be noted that the curvature of $\D^S$ of the subbundle $S$ is given by 
$$\F^S_{\A_S}(\psi)=F^S\cdot\psi+V_{ss}\cdot V_{ss}\cdot\psi.$$
Hence we obtain
$$
\pi_S\circ\K^E_\A(\psi)\circ \pi_s =\K^S_{\A_S}(\psi)+\pi^{\Herm}_{U^{-n}}(A_{sq}\cdot A_{qs}\cdot\psi + V_{sq}\cdot V_{qs}\cdot\psi)
$$

\bgn{lemma}\label{lem:Av+ VA}
$$2\tr\,\pi_{U^{-n}}^{\Herm}\,(A_{sq}\cdot A_{qs}\cdot\psi + V_{sq}\cdot V_{qs}\cdot\psi)=\tr\,\pi_{U^{-n}}^{\Herm }\,(H^{S^\perp}\cdot H^S-H^S\cdot H^{S^\perp})\cdot\psi
$$
\end{lemma}
\bgn{proof}
As in (\ref{eq:Hs}), we have 
$H^{S^\perp }=\A_{sq}, \,\,H^S=\A_{qs}$, 
where $\A_{sq}=A_{sq}+V_{sq},\,\,\, \A_{qs}=A_{qs}+V_{qs}$.
Taking a real basis $\{v_i\}$ of $T_M$ and the real dual basis $\{\t^i\}$ of $T^*_M$, we obtain 
the following 
\bgn{align}
A_{sq}=&\sum_i A_{sq, i}\t^i, \qquad A_{qs}=\sum_i A_{qs, i}\t^i \\
V_{sq}=&\sum_i V_{sq, i}v_i,\qquad V_{qs}=\sum_i V_{qs, i}v_i,
\end{align}
where $A_{sq, i}=-A_{qs, i}^*, \quad V_{sq, i}=-V_{qs, i}^*$.
Thus we have 
\bgn{align}
H^{S^\perp}\cdot H^S=\A_{sq}\cdot\A_{qs}=&(A_{sq}+V_{sq})\cdot(A_{qs}+V_{qs})\\
=&(A_{sq}\cdot A_{qs}+V_{sq}\cdot V_{qs})\\
+&(A_{sq}\cdot V_{qs}+V_{sq}\cdot A_{qs})
\end{align}
Then we have 
\bgn{align}\label{(1)}
(A_{sq}\cdot A_{qs}+V_{sq}\cdot V_{qs})=-&\sum_{i,j}(A_{qs, i}^*A_{qs, j}\t^i\w\t^j+V_{qs, i}^*V_{qs, j}v_i\w v_j)
\end{align}
From $\t^i\cdot v_j+v_j\cdot\t^i=\del_{i,j}$, we obtain 
\bgn{align}
(A_{sq}\cdot V_{qs}+V_{sq}\cdot A_{qs})=&-\sum_{i,j}(A^*_{qs, i}V_{qs, j}\t^iv_j
+V_{qs, i}^*A_{qs, j}v_i\t^j)\\
=-&\sum_{i,j}(A^*_{qs, i}V_{qs, j}\t^iv_j
+V_{qs, j}^*A_{qs, i}v_j\t^i)\\
=-&\sum_{i,j}(A^*_{qs, i}V_{qs, j}
-V_{qs, j}^*A_{qs, i})\t^iv_j-\sum_i V^*_{qs, i}A_{qs, i}\\
=&-\sum_{i, j}P_{i,j}\t^iv_j-\sum_i V^*_{qs, i}A_{qs, i},
\end{align}
where $P_{i,j}$ denotes the skew-Hermitian endmorphism $A^*_{qs, i}V_{qs, j}-V_{qs, j}^*A_{qs, i}.$
We also have 
\bgn{align}
H^S\cdot H^{S^\perp}=
\A_{qs}\cdot\A_{sq}=&(A_{qs}+V_{qs})\cdot(A_{sq}+V_{sq})\\
=&(A_{qs}\cdot A_{sq}+V_{qs}\cdot V_{sq})\\
+&(A_{qs}\cdot V_{sq}+V_{qs}\cdot A_{sq})
\end{align}
Then we have 
\bgn{align}\label{(1)'}
(A_{qs}\cdot A_{sq}+V_{qs}\cdot V_{sq})=
-&\sum_{i,j}(A_{qs, i}A_{qs, j}^*\t^i\w\t^j+V_{qs, i}V_{qs, j}^*v_i\w v_j)\\
=&\sum_{i,j}(A_{qs, j}A_{qs, i}^*\t^i\w\t^j+V_{qs, j}V_{qs, i}^*v_i\w v_j)
\end{align}
We also calculate $(A_{qs}\cdot V_{sq}+V_{qs}\cdot A_{sq})$.
\bgn{align}
(A_{qs}\cdot V_{sq}+V_{qs}\cdot A_{sq})=&-\sum_{i,j}A_{qs, i}V_{qs, j}^*\t^iv_j
+V_{qs, i}A_{qs, j}^*v_i\t^j\\
=-&\sum_{i,j}A_{qs, i}V_{qs, j}^*\t^iv_j
+V_{qs, j}A_{qs, i}^*v_j\t^i\\
=-&\sum_{i,j}(A_{qs, i}V_{qs, j}^*
-V_{qs, j}A_{qs, i}^*)\t^iv_j-\sum_i V_{qs, i}A_{qs, i}^*\\
=&-\sum_{i,j}P'_{i,j}\t^iv_j-\sum_i V_{qs, i}A_{qs, i}^*,
\end{align}
where $P'_{i,j}$ is a skew-Hermitian endomorphism $(A_{qs, i}V_{qs, j}^*-V_{qs, j}A_{qs, i}^*)$.
From (\ref{(1)}) and (\ref{(1)'}),  we have 
\bgn{align}
\tr(A_{sq}\cdot A_{qs}+V_{sq}\cdot V_{qs})=-\tr(A_{qs}\cdot A_{sq}+V_{qs}\cdot V_{sq}))
\end{align}
From Lemma \ref{lem:hontouni key}, we have
$\pi_{U^{-n}}^{\Herm}P_{i,j}\t^iv_i\cdot\psi=0$ and $\pi_{U^{-n}}^{\Herm}P'_{i,j}\t^iv_j\psi=0$. Thus we obtain 
\bgn{align}
&\tr\pi_{U^{-n}}^{\Herm}(A_{sq}\cdot V_{qs}+V_{sq}\cdot A_{qs})-\tr\pi_{U^{-n}}^{\Herm}(A_{qs}\cdot V_{sq}+V_{qs}\cdot A_{sq})\\
&=-\sum_i \tr\pi_{U^{-n}}^{\Herm}\(V^*_{qs, i}A_{qs, i}- V_{qs, i}A_{qs, i}^*\)\psi\\
&=-\sum_i \tr\pi_{U^{-n}}^{\Herm}\(A_{qs, i}V^*_{qs, i}- V_{qs, i}A_{qs, i}^*\)\psi=0
\end{align}
Note that $\(A_{qs, i}V^*_{qs, i}- V_{qs, i}A_{qs, i}^*\)$ is a skew-Hermitian endmorphism.
Hence we have 
\bgn{align}
\tr\pi_{U^{-n}}^{\Herm }(H^{S^\perp}\cdot H^S-H^S\cdot H^{S^\perp})\cdot\psi
=&2\tr\pi_{U^{-n}}^{\Herm}(A_{sq}\cdot A_{qs}\cdot\psi + V_{sq}\cdot V_{qs}\cdot\psi).
\end{align}
\end{proof}
\bgn{proposition}\label{trpiScirc}
$$
\tr\,\pi_S\circ\K^E_\A(\psi)\circ \pi_S =\tr\,\K^S_{\A_S}(\psi)+\frac12 \tr\,\pi_{U^{-n}}^{\Herm }(H^{S^\perp}\cdot H^S-H^S\cdot H^{S^\perp})\cdot\psi
$$
\end{proposition}
\bgn{proof}
The result follows from Lemma \ref{lem:Av+ VA}.
\end{proof}
\bgn{proposition}\label{prop:trpiScircKE}
$$\tr\(\pi_s\circ\K^E_\A(\psi)\circ \pi_s\) =\tr\,\K^S_{\A_S}(\psi)+\|H^S\|^2\psi
$$
\end{proposition}
\bgn{proof}
From Lemma \ref{lem: positive}, we already have 
\bgn{align}\tr\,\pi_{U^{-n}}^{\Herm }(H^{S^\perp}\cdot H^S)\cdot\psi=&+  \|H^S\|^2\psi\\
 \tr\,\pi_{U^{-n}}^{\Herm }(H^{S}\cdot H^{S^\perp})\cdot\psi=&-  \|H^S\|^2\psi
 \end{align}
 Then the result follows from Proposition \ref{trpiScirc}.
\end{proof}
We assume that a generalized holomorphic vector bundle $E$ admits an Einstein-Hermitian metric. 
Then we have 
$\K^E_\A(\psi) =\lam \id_E.$
The component of
$\K^E_\A(\psi)$ of End $(S)$ is given by
$$
\lam \id_S =\K^S_{\A_S}(\psi)+ \pi_{U^{-n}}^{\Herm}\(A_{sq}\cdot A_{qs}\cdot\psi + V_{sq}\cdot V_{qs}\cdot\psi\)
$$
Let $S$ be a generalized holomorphic subbundle of $E$ with rank $S=p.$
From Proposition \ref{prop:trpiScircKE}, we have
$$
p\lam =\tr \,\K^S_{\A_S}(\psi)+ \|H^S\|^2\psi
$$
The degree of $E$ is given by
$$\deg (E):=\int_M\frac{i^n}{2\pi}\tr\,\lan\K^E_\A(\psi), \,\,\ol\psi\ran_s
$$
The slope is given by $\mu(E):=\frac{\deg(E)}{r}$.　Since vol$(M)=i^n\lan \psi, \,\,\ol\psi\ran_s$, we have
$\mu(E)=\frac{\lam}{2\pi}\vol(M)$.
Then we obtain the key inequality:
\bgn{align}\label{key inequality}
\mu(S) =&\frac{\deg(S)}{p}=\frac1p\int_M \frac{i^n}{2\pi}\tr\lan\K_{\A_S}^S(\psi), \,\,\ol\psi\ran_s\notag\\
=&\mu(E)-\int_M  \frac{1}{2\pi} \|H^S\|^2\vol_M\leq \mu(E)
\end{align}

\bgn{proposition}\label{slope inequality EH}
Let $\psi=e^{-\sqrt{-1}\ome}.$
If a generalized holomorphic vector bundle $E$ over a generalized K\"ahler manifold $(M, \J,\J_\psi)$ admits an Einstein-Hermitian metric, then for a generalized holomorphic subbundle $S$ of $E$, we have 
$$
\mu(S)\leq \mu(E)
$$
\end{proposition}
\bgn{proof}
The inequality follows from (\ref{key inequality}).
\end{proof}
We shall extend our result to the general cases of $\psi=e^{b-\sqrt{-1}\ome}.$
\bgn{proposition}
Let $\psi=e^{b-\sqrt{-1}\ome}.$
If a generalized holomorphic vector bundle $E$ over a generalized K\"ahler manifold $(M, \J,\J_\psi)$ admits an Einstein-Hermitian metric, then 
for a generalized holomorphic subbundle $S$ of $E$, we have 
$$
\mu(S)\leq \mu(E)
$$
\end{proposition}
\bgn{proof}
Let $\D^E$ be the canonical connection of $(E, \ol\pa_\J^E, h)$  over a $(M, \J, \J_\psi).$ Since $b$ is a $d$-closed form, the action of $b$-field gives 
a generalized K\"ahler manifold $(M, \J_b, \J_{\psi_b}),$ where 
$\J_b:=\Ad_{e^{-b}}\J$ and $\psi_b=e^{-b}\psi=e^{-\sqrt{-1}\ome}.$
There is also an action of $b$-field on $\ol\pa_\J^E$ by 
$$
\ol\pa_{\J_b}:=e^{-b}\circ \ol\pa^E_\J\circ e^b.
$$
Then $\ol\pa_{\J_b}^E$ is a generalized holomorphic structure with respect to $\J_b.$
As in Proposition \ref{prop: b-field action},
the action of $b$-filed gives a generalized connection $\D^E_b$ by  
$\D^E_b =d+\Ad_{e^{-b}}\A= d+ A- b(V)+ V,$
where $b(V)\in T^*_M\otimes\End(E)$ denotes the contraction of $b$ and $V.$
Then $\D^E_b$ is the canonical connection of a generalized holomorphic vector bundle 
$(E, \ol\pa_{\J_b}^E, h).$
From Proposition \ref{prop: b-field action}, 
We have 
$\F_\A^E(\psi)=\Ad_{e^{-b}}\F_{\A_b}^E$.
Thus $\deg(E)$ is invariant under the action of $b$-field. 
Let $S$ be a generalized holomorphic subbundle. 
Then $b$-filed also acts on $(S, \ol\pa_\J^S)$ to give a generalized holomorphic subbundle 
$(S, \ol\pa_{\J_b}^S)$ of $(E, \ol\pa_{\J_b}^E).$
Since $\deg (S)$ is also invariant under the action of $b$-field, 
we obtain 
$\deg(S)\leq \deg(E)$
from Proposition \ref{slope inequality EH}. 
\end{proof}
\bgn{proposition}
Let $\psi=e^{b-\sqrt{-1}\ome}.$
If a generalized holomorphic vector bundle $E$ over a generalized K\"ahler manifold $(M, \J,\J_\psi)$ admits an Einstein-Hermitian metric, then 
for a weak generalized holomorphic subbundle $\pi$ of $E$, we have 
$$
\mu(\pi)\leq \mu(E)
$$
\end{proposition}
\bgn{proof}
Let $F$ be a generalized holomorphic subbundle which is given by $\pi$ on the complement $M\bsh S,$ where $S$ is a subset of $M$ of codim $2$.
From Proposition \ref{leths} and Proposition \ref{prop:trpiScircKE}, on $M\bsh S$, we have 
$$\tr\(\pi\circ\K^E_\A(\psi)\circ \pi\) =\tr\,\K^F_{\A_F}(\psi)+\|\pa_0\pi\|^2\psi.$$
Since $\pi\in L^2_1\End(E),$ $\pa_0\pi$ is square-integrable.
Thus if $(E, h)$ has an Einstein-Hermitian metric, then we also have
\bgn{align}\label{key inequality for pi}
\mu(\pi) =&\frac{\deg(S)}{p}=\frac1p\int_M \frac{i^n}{2\pi}\tr\lan\K_{\A_S}^S(\psi), \,\,\ol\psi\ran_s\notag\\
=&\mu(E)-\int_M  \frac{1}{2\pi} \|\pa_0\pi\|^2\vol_M\leq \mu(E)
\end{align}
 Thus we have $\mu(\pi)\leq \mu(E).$
Since $\mu(\pi)$ is also invariant under the action of $b$-field, we have the result.
\end{proof}
\bgn{proposition}
If a generalized holomorphic vector bundle $E$ admits an Einstein-Hermitian metric, then 
$E$ is polystable.
\end{proposition}
\bgn{proof}
Let $E$ be an  irreducible holomorphic vector bundle with an Einstein-Hermitian metric. 
Then it follows from (\ref{key inequality}) that $\mu(S)\leq \mu(E)$ for a generalized holomorphic subbundle $S$. 
If the second fundamental form $H^S$ vanishes, $E$ is decomposed into $S\oplus S'$. 
Since $E$ is irreducible, $H^S$ does not vanish. Thus we have the strict inequality 
$\mu(S)<\mu(E)$. Hence $E$ is stable. 
We assume that $E$ is decomposed into $E_1\oplus\cdots \oplus E_m$, where each $E_i$ is irreducible. 
Then each $E_i$ is stable with the same Einstein factor $\lam.$
\end{proof}
\bgn{remark}
From Proposition \ref{prop:trpiScircKE} and Proposition \ref{HSperp}, we have 
$$\tr\(\pi_S\circ\K^E_\A(\psi)\circ \pi_s\) =\tr\,\K^S_{\A_S}(\psi)+\|\ol\pa\pi\|^2\psi
$$
Since the first Chern form of the generalized subbundle is given by $\tr\,\K^S_{\A_S}(\psi)$ and 
$\pi \in L^2_1(\End(E)),$
it turns out that 
the first Chern form is integrable which coincides with the one given by the first Chern class of the corresponding sheaves. 
\end{remark}
\section{Variation formula of mean curvature}\label{Variation formula of mean curvature}
Let $(E, h_0)$ be a complex Hermitian vector bundle over a generalized K\"ahler manifold $(M, \J,\J_\psi)$.
Let $\ol\pa_{\J}^E$ be a generalized holomorphic structure on $E$ and 
$\D_0$ the canonical connection of $(E, h_0, \ol\pa^E_\J).$ 
We denote by
$\End(E, h_0)$ Hermitian endmorphisms of $E$
 with respect to $h_0$ and we also denoted by 
 $\Herm^+(E, h_0)$ positive-definite Hermitian endmorphisms with respect to $h_0.$
For $f\in \Herm^+(E, h_0)$, we define an Hermitian metric $h_{f}$ by 
$$
h_{f}(e_1, e_2):=h_0(fe_1, e_2), 
$$
for $e_1, e_2\in E.$
We denote by $\D_{f}$ the canonical connection of $(E, h_{f}, \ol\pa^E_\J)$
with the connection form $\A_{f}.$
\bgn{lemma}\label{D1,0f=fD1,0f{-1}}
$\D^{1,0}_{f}={f}^{-1}\circ \D^{1,0}_0\circ {f}$
\end{lemma}
\bgn{proof}
From the definition of $\D_{f}$, we have 
\bgn{align}\label{paJhte1e2=htDft}
\pa_\J h_{f}(e_1, e_2) =h_{f}(\D_{f}^{1,0}e_1, \, e_2)+h_{f}(e_1, \,\ol\pa^E_\J e_2)
\end{align}
We also have 
$$
\pa_\J h_0({f} e_1, \, e_2)=h_0(\D^{1,0}_0({f}e_1), \, e_2)+ h_0({f} e_1, \, \ol\pa^E_\J e_2)
$$
Since $ h_0(e_1, \, e_2) =h_{f}({f}^{-1}e_1, \, e_2)$, we have
\bgn{align}\label{ht-1D10}
\pa_\J h_{f}( e_1, \, e_2)=h_{f}(f^{-1}\circ\D^{1,0}_0({f}e_1),\,\, \, e_2)+ h_{f}(e_1, \, \ol\pa^E_\J e_2)
\end{align}
Thus from (\ref{paJhte1e2=htDft}) and (\ref{ht-1D10}), we obtain 
$\D_{f}^{1,0}=f^{-1}\circ\D_0^{1,0}\circ {f}.$
Hence the result follows.
\end{proof}
Then the connection form $\A_{f}^{1,0}$ is given by 
\bgn{equation}\label{Af=f-1paJf}
\A_{f}^{1,0}={f}^{-1}\pa_\J {f}+ {f}^{-1}\A_0^{1,0}{f},
\end{equation}
where $\A_0^{1,0}$ denotes the $(1,0)$-component of the connection form of the connection $\D_0.$
Let  $\{f_t\}$ be a smooth family of $\Herm^+(E, h_0)$ with $f_0=\id_E.$
We denote by $\K_{f_t}(\psi)$ the mean curvature of the generalized connection 
$\D_{f_t}.$
Then the derivative of variation of the mean curvature $\K_{f_t}(\psi)$ 
is given by 
\bgn{proposition}\label{ddtKftpsi}
For all $t_0,$ we have
$$
\pi^{\Herm_{f_t}}\(\frac{d}{dt}\K_{f_t}(\psi)\)\Big|_{t={t_0}}=\pi^{\Herm_{f_t}}_{U^{-n}}\l(d^{\D_{f_t}}\overset{\,\,\,\scriptscriptstyle\bullet}{\A}_{f_t}\cdot\psi\r)\Big|_{t={t_0}},
$$
where ${\overset{\,\,\,\scriptscriptstyle\bullet}{\A}}_{f_t}=\frac{d}{dt}\A_{f_t}$ and $\pi^{\Herm_{f_t}}$ denotes the 
projection from $\End(E)$ to the Hermitian endmorphisms with respect to $h_{f_t}.$
\end{proposition}
We need the following Proposition to prove Proposition \ref{ddtKftpsi},
which implies that the mean curvature is regarded as the moment map in \cite{Goto_2017}.
\bgn{proposition}\label{variation of fixed h0}
Let $h_0$ be a fixed Hermitian metric on $E$ and 
$\{\A_t\}$ a family of Hermitian generalized connection of $(E, h)$ with the respect to $h_0.$
Then the derivative of variation of the mean curvature is given by 
\bgn{equation}\label{fracddtKAt}
\frac{d}{dt}\K_{\A_t}(\psi)=\pi^{\Herm_{h_0}}_{U^{-n}}\l(d^{\A_{t}}{\overset{\,\,\,\scriptscriptstyle\bullet}{\A}}_{t}\cdot\psi\r),
\end{equation}
\end{proposition}
\bgn{proof}
By using an action of $b$-field, our Proposition is reduced to the cases of $\psi=e^{-\sqrt{-1}\ome}$.
Thus it suffices to show our Proposition in the cases of $\psi=e^{-\sqrt{-1}\ome}.$
Then the mean curvature $\K_{\A_t}(\psi)$ is given by 
\bgn{align}
\K_{\A_t}(\psi)=&\pi_{U^{-n}}^{\Herm_{h_0}}\(F_{A_t}\cdot\psi +\frac12[V_t\cdot V_t]\cdot\psi\),
\end{align}
where $\A_t=A_t+V_t$ and $d^A:=d+A_t$ denotes the ordinary connections and $V_t\in T_M\otimes u(E).$
In fact, from Lemma \ref{lem:hontouni key}, we see that 
$\pi^{\Herm_{h_0}}_{U^{-n}}d^{A_t}(V_t\cdot\psi)=0$ in the cases of $\psi=e^{-\sqrt{-1}\ome}.$
Then we have 
\bgn{align}
\frac{d}{dt}\K_{\A_t}(\psi) =\pi^{\Herm_{h_0}}_{U^{-n}}\(d^A(\overset{\bullet}{A_t}\cdot\psi)+ [\overset{\bullet}{V_t}\cdot V_t]\cdot\psi\)
\end{align}
Applying Lemma \ref{lem:hontouni key} again, we obtain
$\pi^{\Herm_{h_0}}_{U^{-n}}\(d^{A_t}(\overset{\bullet}{V_t}\cdot\psi) +
[V_t\cdot\overset{\bullet}{A_t}]\cdot\psi\)=0.$ 
Then the right hand side of (\ref{fracddtKAt}) is given by 
\bgn{align}
\pi^{\Herm_{h_0}}_{U^{-n}}d^{\A_t}\overset{\bullet}{\A_t}\cdot\psi=&
\pi^{\Herm_{h_0}}_{U^{-n}} \(d^{A_t}\overset{\bullet}{A_t}\cdot\psi+ [\overset{\bullet}{V_t}\cdot V_t]\cdot\psi\).
\end{align}
Thus we have the result.
\end{proof}

\bgn{proof}[{\sc Proof of Proposition \ref{ddtKftpsi}}]
There is a square $f^{\frac12}\in \Herm^+(E, h_0)$ such that
$h_f(e_1, e_2)=h_0(f^{\frac12}e_1, f^{\frac12}e_2).$
Since the relation between $\pi^{\Herm_{f}}$ and $\pi^{\Herm_{h_0}}$ is given by 
\bgn{equation}\label{piHermhf}
\Ad_{f^{\frac12}}\circ\pi^{\Herm_{f}}=\pi^{\Herm_{h_0}}\circ \Ad_{f^{\frac12}},
\end{equation}
then $f^{\frac12}\circ\D_f\circ f^{\frac{-1}2}$ is an Hermitian generalized connection with respect to $h_0.$
In fact, 
$dh_f(e_1, e_2)= h_f(\D_f e_1, e_2)+ h_f(e_1, \D_fe_2)$ yields the following: 
\bgn{align}
dh_0(e_1, e_2)=&dh_f(f^{\frac{-1}2}e_1, f^{\frac{-1}2}e_2)\\
=&h_f(\D_f\circ f^{\frac{-1}2}e_1,\,\, f^{\frac{-1}2}e_2)+
h_f(f^{\frac{-1}2} e_1,\,\,\D_f\circ f^{\frac{-1}2}e_2)\\
=&h_0(f^{\frac12}\circ\D_f\circ f^{\frac{-1}2}e_1,\,\,e_2)+
h_0(e_1, \,\,f^{\frac12}\circ\D_f\circ f^{\frac{-1}2}e_2)
\end{align}
We denote by $\h\D_f$ the generalized connection $f^{\frac12}\circ\D_f\circ f^{\frac{-1}2}.$
Since $f^{\frac12}$ acts on the curvature by the Adjoint action,
we also have 
\bgn{equation}\label{Fffrac12D}
\F_{\h\D_f}(\psi)=
f^{\frac12}\circ \F_{\D_f}(\psi)\circ f^{\frac{-1}2}=\Ad_{f^{\frac12}}\(\F_{\D_f}(\psi)\)
\end{equation}
We shall go back to our proof of Proposition \ref{ddtKftpsi}.
We shall reduce the cases of Proposition \ref{ddtKftpsi} to the one of Proposition \ref{variation of fixed h0} by using the gauge transformation $\Ad_{f^{\frac12}}$. 

When we replace $h_0$ by $h_{f_{t_0}}$
and use $\til{f_t} =f_{t_0}^{-1}f_t$, Proposition \ref{variation of fixed h0} reduces to the case $t_0=0.$
Thus it suffices to show Proposition \ref{variation of fixed h0} in the case $t_0=0$.
From (\ref{piHermhf}) and (\ref{Fffrac12D}), we obtain 
\bgn{align}
\Ad_{ f_t^{\frac12}}\circ\pi^{\Herm_{{f_t}}}_{U^{-n}}\(\F_{\D_{f_t}}(\psi)\)=&
\pi^{\Herm_{{f_{0}}}}_{U^{-n}}\circ\Ad_{ f_t}^{\frac12}\(\F_{\D_{ f_t}}(\psi)\)\\
=&\pi^{\Herm_{{f_{0}}}}_{U^{-n}}\(\F_{\h\D_{f_t}}(\psi)\)
\end{align}
Since $\pi^{\Herm_{{f_t}}}_{U^{-n}}\F_{\D_{f_t}}(\psi)=\K_{f_t}(\psi)$, applying $\Ad_{f_t^{\frac{-1}2}}$ to the both sides and taking the differential at $t=0$, we have 
 \bgn{align}\label{fracddtpiHermftU-nFDft}
 \frac{d}{dt}\K_{\D_{f_t}}(\psi)\big|_{t=0}=
 \Ad_{ f^{\frac{-1}2}_{0}}\(\frac{d}{dt}
\pi^{\Herm_{f_0}}_{U^{-n}}\F_{\h\D_{f_t}}(\psi)|_{t=0}\)
+\frac{d}{dt}\Ad_{f_t^\frac{-1}2}(\K_{\h\D_{f_{0}}}(\psi))\Big|_{t=0}
  \end{align}
Since $f_{0}=\id_E,$ we have 
$\F_{\h\D_{f_{0}}}(\psi)=\F_{\D_{ f_{0}}}(\psi)$. 
Thus we have 
\bgn{align}\label{fracddtAdtillftdeax-12K}
\frac{d}{dt}\Ad_{ f_t^\frac{-1}2}(\K_{\h\D_{f_{0}}}(\psi))\Big|_{t=0}=\Big[\frac{d}{dt} f_t^{\frac{-1}2}\Big|_{t=0}, \,\,\,\,\,(\K_{\D_{f_{0}}}(\psi))\Big]
\end{align}

Since $f_t$ and $(\K_{\h\D_{f_{0}}}(\psi))$ are Hermitian with respect to $h_{f_{0}}=h_0$, 
the bracket of (\ref{fracddtAdtillftdeax-12K}) is Skew-Hermitian. 
Thus 
$$
\pi^{\Herm_{ f_{0}}}\(\frac{d}{dt}\Ad_{ f_t^\frac{-1}2}(\K_{\h\D_{f_{0}}}(\psi))\Big|_{t=t_0}\)=0
$$
From (\ref{fracddtpiHermftU-nFDft}) and $f_0=\id_E$, we have 
\bgn{equation}
\pi^{\Herm_{ f_{0}}} \(\frac{d}{dt}\K_{\D_{f_t}}(\psi)\)\Big|_{t=0}=
\( \frac{d}{dt}
\pi^{\Herm_{{f_{0}}}}_{U^{-n}}\F_{\h\D_{f_t}}(\psi)|_{t=0}\)
\end{equation}
Since $\h\D_{f_t}$ is an Hermitian generalized connection with respect to  the metric $h_{0},$ we can apply our formula (\ref{fracddtKAt}) in Proposition \ref{variation of fixed h0}
to obtain 
\bgn{equation}\label{fracddtKDftpsit=t0}
\pi^{\Herm_{ f_{0}}}\(\frac{d}{dt}
\K_{\D_{f_t}}(\psi)\Big|_{t=0}\)
=\pi^{\Herm_{{f_{0}}}}_{U^{-n}}\(d^{\h\D_{f_{0}}}\overset{\,\,\bullet}{\widehat{\A}}_{f_{0}}\cdot\psi\)
\end{equation}
From $d^{\h\D_{f_t}}\overset{\,\,\bullet}{\widehat{\A}}_{f_t}
=f^{\frac12}_t\circ (d^{\D_{f_t}}\overset{\,\,\,\bullet}{{\A}}_{f_t})\circ f^{-\frac12}_t,
$
it follows from $f_0=\id_E$ that we have 
$
d^{\h\D_{ f_{0}}}\overset{\,\,\bullet}{\widehat{\A}}_{ f_{0}}
=
d^{\D_{{ f_{0}}}}{\overset{\,\,\,\scriptscriptstyle\bullet}{\A}}_{{ f_{0}}}.
$
From (\ref{fracddtKDftpsit=t0}) we obtain 
$$
\pi^{\Herm_{{f_{0}}}}\(\frac{d}{dt}\K_{f_t}(\psi)\)
\Big|_{t=0}=\pi^{\Herm_{f_{0}}}_{U^{-n}}\l(d^{\D_{f_t}}{\overset{\,\,\,\scriptscriptstyle\bullet}{\A}}_{f_t}\cdot\psi\r)\Big|_{t=0}.
$$
Hence we obtain the result.
 \end{proof}
 \bgn{proposition}For all $t=t_0,$ we have 
 $$
\pi^{\Herm_{f_t}}\(\frac{d}{dt}\K_{f_t}(\psi)\)\Big|_{t={t_0}}=\pi^{\Herm_{f_t}}_{U^{-n}}\l(\ol\pa^E_{\J}\overset{\,\,\,\scriptscriptstyle\bullet}{\A}_{f_t}\cdot\psi\r)\Big|_{t={t_0}},
$$
 \end{proposition}
 \bgn{proof}
 Since $\D_{f_t}$ is a family of canonical connections, i.e., $\D^{0,1}_{f_t}=\ol\pa^E_\J.$
 Thus $\overset{\,\,\,\scriptscriptstyle\bullet}{\A}_{f_t}=\overset{\,\,\,\scriptscriptstyle\bullet}{\A}^{1,0}_{f_t}\in\L_\J\otimes\End(E)$.
 Since the $\pi_{U^{-n}}$-component is given by $(1,1)$-component, i.e., $(\L_\J\cdot\ol\L_\J)\cdot\psi$, we have
 $$
 \pi^{\Herm_{f_t}}_{U^{-n}}\l(d^{\D_{f_t}}\overset{\,\,\,\scriptscriptstyle\bullet}{\A}_{f_t}\r)\Big|_{t={t_0}}
 =\pi^{\Herm_{f_t}}_{U^{-n}}\l(\ol\pa^E_{\J}\overset{\,\,\,\scriptscriptstyle\bullet}{\A}_{f_t}\r)\Big|_{t={t_0}}$$ 
 \end{proof}
 \bgn{proposition}\label{fracddtKftpsi=piHermftU-n}
$$
\pi^{\Herm_{f_t}}\(\frac{d}{dt}\K_{f_t}(\psi)\)=\pi^{\Herm_{f_t}}_{U^{-n}}
\(\ol\pa^E_\J\circ\D^{1,0}_{f_t}(f_t^{-1}\overset{\scriptscriptstyle\bullet}{f_t})\cdot\psi\),
$$
 \end{proposition}
 \bgn{proof}
 The generalized connection $\D^{1,0}_{f_t}$ is given by the gauge transformation from Lemma \ref{D1,0f=fD1,0f{-1}},
 $$
 \D^{1,0}_{f_t}=f_t^{-1}\circ\D^{1,0}_0\circ f_t
 $$
 Thus we have 
 $$
 {\overset{\,\,\,\scriptscriptstyle\bullet}{\A}}_{f_t}=\D^{1,0}_{f_t}(f_t^{-1}{\overset{\scriptscriptstyle\bullet}{f_t}})
 $$
 Then the result follows from Proposition \ref{ddtKftpsi}.
 \end{proof}
 \bgn{remark}We also have the following:
 $$(f_t^{-1}{\overset{\scriptscriptstyle\bullet}{f_t}})
=h_t^{-1}\frac{\pa}{\pa t}h_t$$
\end{remark}
Taking the trace of the both side of Proposition \ref{fracddtKftpsi=piHermftU-n}, we obtain
\bgn{proposition}\label{trace}
$$\(\frac{d}{dt}\tr\,\K_{f_t}(\psi)\)=\pi^{\Herm_{f_t}}_{U^{-n}}\(\ol\pa^E_\J\circ\D^{1,0}_{f_t}\tr\,(f_t^{-1}{\overset{\scriptscriptstyle\bullet}{f_t}})\cdot\psi\)
$$
\end{proposition}
\bgn{proof}
The result follows from Proposition \ref{ddtKftpsi} and
\ref{fracddtKftpsi=piHermftU-n}, since we have 
$$
 \tr \frac{d}{dt}\Ad_{\til f_t^\frac{-1}2}(\K_{\h\D_{f_{t_0}}}(\psi))\Big|_{t=t_0}=\tr[\frac{d}{dt}\til f_t^{\frac{-1}2}, \,\,\,(\F_{\h\D_{f_{t_0}}}(\psi))]=0
$$
\end{proof}
\bgn{lemma}\label{Herm ST ST}
Let $h_0$ be an Hermitian metric of $E$ and $h_f$
 be an Hermitian metric given by $f\in\Herm^+(E, h_0)$ as before. 
 Let $T_f$ be a section of $\Herm(E, h_0)$ satisfying $[f, T_f]=0.$
 Then for any section $S\in \End(E),$ we have 
 $$
 h_0\(\pi^{{\ss\Herm_{h_f}}}(S), \,\,T_f\) =h_0(S, \,\,T_f).
 $$ 
\end{lemma}
\bgn{proof}
From (\ref{piHermhf}), we have  
$$
\Ad_{f^{\frac12}}\circ\pi^{\Herm_{f}}=\pi^{\Herm_{h_0}}\circ \Ad_{f^{\frac12}}
$$
 Since $h_0(\Ad_{f^{\frac12}}A, \,\,B) =h_0(A, \,\, \Ad_{f^{\frac{1}2}}B)$ for $A, B\in \End(E)$ and $[f, T_f]=0,$
 $T_f\in \Herm(E, h_0),$
 we have
\bgn{align}
h_0\(\pi^{{\ss\Herm_{h_f}}}(S), \,\,T_f\)=&h_0(\pi^{\Herm_{h_0}}\Ad_{f^{\frac12}}S, \,\,\Ad_{f^{\frac{-1}2}} T_f)\\
=&h_0(\Ad_{f^{\frac12}}S,\,\,T_f)\\
=&h_0(S, T_f)
\end{align}
\end{proof}

\bgn{proposition}\label{Kftpsiint}
Let $h_0$ be an Hermitian metric and $h_f:=h_0f$ an Hermitian metric given by 
$f\in \Herm^+(E, h_0)$. Let $T_{f}\in \Herm(E, h_0)$ be an Hermitian endmorphism satisfying $[T_{f}, f]=0.$
Then we have  
\bgn{align}
h_0(\K_{f}(\psi)-\K_{0}(\psi), \,\, T_{f}) =&h_0(\ol\pa^E_\J(\A_{f_1}-\A_{0})\cdot\psi, \,\,\,T_{f}\ol\psi)\\
=&h_0(\ol\pa^E_\J( f^{-1}\pa_0 f\cdot\psi), \,\,T_f\cdot\ol\psi)
\end{align}\end{proposition}
\bgn{proof}
Let $\xi=\log (h_0^{-1}h_f)$. Then we have a family
$\{f_t= e^{\xi t}\}$ of $\Herm^+(E, h_0)$ satisfying $f_{1}=f$ and $f_0=\id_E.$
Since $f_t=e^{\xi t}$ commutes with $T_{f},$
we have 
$$
h_0(\K_{f}(\psi), \,\, T_{f}\cdot\ol\psi)=
h_0(\Ad_{f^{\frac12}}(\K_{f}(\psi)), \,\, T_{f}\cdot\ol\psi)
$$
Since $\Ad_{f^{\frac12}}(\K_{f}(\psi))=\K_{\h\D_{f}},$
We have
$$
\Ad_{f^{\frac12}}(\K_{f}(\psi))-\K_{f_{0}}(\psi)=\K_{\h\D_{f}}-\K_{\h\D_{{0}}}.
$$
Then we have 
\bgn{align}
\K_{\h\D_{f}}(\psi)-\K_{\h\D_{0}}(\psi)=\int^{1}_{0}\frac{d}{dt}\K_{\h\D_{f_t}}(\psi)dt
\end{align}
Since $\h\D_{f_t}$ is a family of generalized Hermitian connections with respect to the fixed $h_0$, 
we apply (\ref{fracddtKAt}) to obtain 
\bgn{align}
\int^{1}_{0}\pi^{\Herm_{h_0}}_{U^{-n}}\frac{d}{dt}\K_{\h\D_{f_t}}(\psi)ds
=&\int^1_0\pi^{\Herm_{h_0}}_{U^{-n}}d^{\h\A_{f_t}}\overset{\bullet}{\h\A_{f_t}}\cdot\psi\\
=&\int^1_0\pi^{\Herm_{h_{0}}}_{U^{-n}}\Ad_{f^{\frac12}_t}\(d^{\A_{f_t}}\overset{\bullet}{\A_{f_t}}\cdot\psi\)\\
\end{align}
Since $T_f\in \Herm(E, h_0),$ we have
\bgn{align}
h_0(\K_{f}(\psi)-\K_{0}(\psi), \,\, T_{f}\cdot\ol\psi)=&h_0^{\top}(\int^{1}_{0}\frac{d}{dt}\K_{\h\D_{f_t}}(\psi)dt, \,\, T_f\cdot\ol\psi)\\
=&h_0^{\top}(\int^{1}_{0}\pi^{\Herm_{h_0}}\frac{d}{dt}\K_{\h\D_{f_t}}(\psi)dt, \,\, T_f\cdot\ol\psi)\\
=&\int^1_0h_0^{\top}(\Ad_{f^{\frac12}_t}\(d^{\A_{f_t}}\overset{\bullet}{\A_{f_t}}\cdot\psi\), \,\,T_f\cdot\ol\psi)dt\\
\end{align}
Since $[f_t, T_f]=0, f\in \Herm(E, h_f)$ and 
$\overset{\bullet}{\A_{f_t}}=\overset{\bullet}{\A_{f_t}^{1,0}},$ we have
\bgn{align}
h_0(\K_{f}(\psi)-\K_{0}(\psi), \,\, T_{f}\cdot\ol\psi)=&\int^1_0h_0(\(\ol\pa_\J^E\overset{\bullet}{\A_{f_t}^{1,0}}\cdot\psi\), \,\,\Ad_{f^{\frac12}_t}(T_f)\cdot\ol\psi)dt\\
=&\int^1_0h_0(\(\ol\pa_\J^E\overset{\bullet}{\A_{f_t}^{1,0}}\cdot\psi\), \,\,(T_f)\cdot\ol\psi)dt\\
=&h_0(\ol\pa^E_\J(\A_{f_1}-\A_{0})\cdot\psi, \,\, T_f\cdot\ol\psi)\\
=&h_0(\ol\pa^E_\J (f^{-1}\pa_0 f)\cdot\psi, \,\, T_f\cdot\ol\psi)
\end{align}
Thus we obtain the result.
\end{proof}
\bgn{remark}
For instance, we will apply Proposition \ref{Kftpsiint} to $T_f=\log f$ or $f^\sig,$ for $0\leq\sig\leq 1.$
\end{remark}
\section{Construction of Einstein-Hermitian metrics on stable generalized holomorphic bundles}
\label{Construction of Einstein-Hermitian metrics on stable generalized holomorphic bundles}
\subsection{The continuity method}
We shall use the continuity method to obtain Einstein-Hermitian metrics on polystable generalized holomorphic bundles.
We will use the same notations as before.
Given an Hermitian metric $h_0$ on a generalized holomorphic vector bundle $E$, we denote the canonical generalized connection 
by $\D_0=\pa_0+\ol\pa^E$ $\D_0$ and the mean curvature by $\K_0(\psi)$.
An Hermitian metric $h_f$ is given by 
$h_f(s_1, s_2)=h_0(fs_1, s_2)$, where $f\in \Herm^+(E, h_0)$ and $s_1, s_2\in \Gam(E).$
Then we have the canonical generalized connection $\D_{f}$ which is associated with
$h_f$ and $\K_f(\psi)$ is the mean curvature of $\D_f$.
Then
the Einstein-Hermitian condition is given by
$$
\K_f(\psi)=\lam \id_E\psi,
$$
where $\lam$ is the Einstein constant.
By using $\psi,$ $\Herm(E)\otimes U^{-n}$ is identified with $\Herm(E).$
For abuse of notation, we often consider $\K_f(\psi)$ as the Hermitian endmorphism
under the identification.
For $\e\in [0,1]$, we introduce the following equation on which the 
continuity method is applied
\bgn{equation}\label{eq:continuity method1}
L_\e(f):=\K_f(\psi)-\lam \id_E +\e(\log f )=0
\end{equation}
Note that the solution of the equation gives an Einstein-Hermitian metric if $\e=0.$ 
We define the subset ${\mathcal S}\subset[0,1]$ by
$$
{\mathcal S}:=\{ \, \e\in [0,1]\, |\, \text{\rm the equation (\ref{eq:continuity method1}) has an solution}\,\}
$$
Let $\h\K_f(\psi):=\K_f(\psi)-\lam \id_E$. Our construction of the solution is divided into four steps \par\medskip 
\noindent 
Step 0. The subset ${\mathcal S}$ contains $1$. \par\noindent
Step 1. ${\mathcal S}\subset [0,1]$ is a nonempty, open set.
\par\noindent
Step 2. $(0,1]\subset{\mathcal S}$\par\noindent
Step 3. ${\mathcal S}=[0,1]$\\
We shall show these steps one by one in the following subsections.
\subsection{Preliminary results for Step 0 and Step 1}
According to the decomposition $(\TT)^\C=\L_\J\oplus\ol\L_\J,$ the exterior 
derivative $d$ is decomposed into $d=\pa+\ol\pa.$
In a generalized K\"ahler manifold, we also use the decomposition 
$\L_\J=\L_\J^+\oplus \L^-_\J$ to obtain  
$\pa=\del_++\del_-$ and $\olpa=\ol\del_++\ol\del_-.$
These four differential operators act on differential forms and we have the adjoint operators $\del_\pm^*$ and 
$\ol\del^*_\pm$.
We denote by $\trian_d$ the Laplacian of $d.$
We also have the Laplacians $\trian_{\pa_\J}$ and 
$\trian_{\ol\pa_J}$, respectively.
Let $\trian_{\del_\pm}$ be the Laplacians $\del_\pm\del_\pm^*+\del_\pm^*\del_\pm$ which acts on differential forms on $M.$
Then the following generalized K\"ahler identity holds 
$\del_+^*=-\del_+$ and $\del_-^*=\del_-.$
By the generalized K\"ahler identity, 
we have 
$$
\trian_d=2\trian_{\pa_\J}=2\trian_{\ol\pa_\J}=4\trian_{\del_-}=4\trian_{\del_+}=4\trian_{\ol\del_-}=4\trian_{\ol\del_-}
$$
Let $h$ be an Hermitian metric on $E$ and $\xi\in \Herm(E, h).$
We denote by $h_s$ the $1$-parameter family of Hermitian metrics of $E$ which are given by 
$h_s(e_1, e_2):= h(e^{\xi s}e_1, e_2)$, where $s$ is a parameter and
$h_0=h.$
Let $(E, \ol\pa_\J^E)$ be a generalized holomorphic vector bundle and  $\D_s=d+\A_s$ the $1$-parameter family of canonical generalized Hermitian connections with respect to $h_s.$
Then we have 
\bgn{lemma}\label{dotA=fracddsAs}
$${\overset{\,\,\,\scriptscriptstyle\bullet}{\A}}=\frac{d}{ds}\A_s|_{t=0}=\D^{1,0}_h\xi
$$
\end{lemma}
\bgn{proof}
From Lemma \ref{D1,0f=fD1,0f{-1}}, the connection form $\A_s$ is given by 
$$
\A_s=(e^{-\xi s}\pa_\J e^{\xi s}) s+e^{-\xi s}\A_h^{1,0}e^{\xi s},
$$
where $\A_h=\A_0$ and $\D_h=d+\A_h.$
Thus we have 
$$
\frac{d}{ds}\A_s|_{t=0}=\pa_\J \xi +[\A_h^{1,0},\,\,\xi]=\D^{1,0}_h\xi.
$$
\end{proof}
\bgn{proposition}\label{fracddsKhs|t=0trianh}
$$\pi^{\Herm_h}\(\frac{d}{ds}\K_{h_s}(\psi)|_{t=0}\)=\trian_h (\xi\cdot\psi),$$
where $\trian_h$ is the Laplacian $\D^{1,0}_-(\D^{1,0}_-)^*+(\D^{1,0}_-)^*\D^{1,0}_-$ actin on $\End(E).$
\end{proposition}
\bgn{proof}
The key point is Proposition \ref{fracddtKftpsi=piHermftU-n}:
$$\pi^{\Herm_h}\(\frac{d}{ds}\K_{h_s}(\psi)|_{s=0}\)=\pi_{U^{-n}}^{\Herm_h}d^{\D_h}({\overset{\,\,\,\scriptscriptstyle\bullet}{\A}}\cdot\psi),
$$
where ${\overset{\,\,\,\scriptscriptstyle\bullet}{\A}}=\frac{d}{ds}\A_{h_s}|_{s=0}\in \End(E)\otimes(\TT).$
From Lemma \ref{dotA=fracddsAs}, we have 
$$
\pi^{\Herm_h}\(\frac{d}{ds}\K_{h_s}(\psi)|_{s=0}\)=\pi^{\Herm}d^{\D_h}(\D^{1,0}_h\xi\cdot\psi)
$$
Then we also have
 $$\pi_{U^{-n}}^{\Herm}d^{\D_h}(\D^{1,0}_h\xi\cdot\psi)=\pi_{U^{-n}}^{\Herm}\D^{0,1}_h(\D^{1,0}_h\xi\cdot\psi)$$
 Since $\xi\cdot\psi\in U^{-n},$ it follows $(\D^{1,0}_{h\,+})(\xi\cdot\psi)=0.$
 Thus we have 
 \bgn{align}
 \pi^{\Herm_h}\(\frac{d}{ds}\K_{h_s}(\psi)|_{t=0}\)=&\pi_{U^{-n}}^{\Herm}\D^{0,1}_h(\D^{1,0}_{h\, -})(\xi\cdot\psi) \\
 =&\pi_{U^{-n}}^{\Herm}(\D^{0,1}_{h\,-})(\D^{1,0}_{h\, -})(\xi\cdot\psi)
 \end{align}
 Applying the following generalized K\"ahler identity
 $$
( \D^{0,1}_h)_-=(\D^{0,1}_{h\,-})^*, \qquad ( \D^{0,1}_{h\,+})=-(\D^{0,1}_{h\,+})^*,
 $$
 We obtain 
 \bgn{align}
 \pi^{\Herm_h}\(\frac{d}{ds}\K_{h_s}(\psi)|_{t=0}\)=&(\D^{1,0}_{h\,-})^*(\D^{1,0}_{h\, -})(\xi\cdot\psi)\\
 =&\((\D^{1,0}_{h\,-})^*(\D^{1,0}_{h\, -})+(\D^{1,0}_{h\,-})(\D^{1,0}_{h\, -})^*\)(\xi\cdot\psi)\\
 =&\trian_h(\xi\cdot\psi),
 \end{align}
 since $(\D^{1,0}_{h\,-})^*(\xi\psi)=0$ and $\trian_h$ is a real self adjoint operator.
\end{proof}
\subsection{Step 0}
We use the same notations as before.
Let
$Q:=\h\K_{f}(\psi)-\h\K_{f_0}(\psi)$.
Then we have
\bgn{proposition}\label{tr Qtrian}
$$
\tr \,Q =\trian\tr\log f
$$
\end{proposition}
\bgn{proof}
From Proposition \ref{trace},
the trace of $Q$ is given by
$$
\tr \,Q =\int^1_0\pi_{U^{-n}}^{\Herm_{f_t}}\l(\ol\pa^E\tr\,\dot{\A_{f_t}}\cdot\psi\r)dt
$$
Then we see 
$$
\tr \Dot{\A_{f_t}}=\tr\,\frac{d}{dt} f_t^{-1}(\pa_0f_t)\psi =\pa(\tr\log {\overset{\scriptscriptstyle\bullet}{f_t}})
$$
Then we have 
$$
\tr\, Q=\int^1_0\pi_{U^{-n}}^{\Herm_{f_t}}\l(\ol\pa\pa(\tr\log {\overset{\scriptscriptstyle\bullet}{f_t}})\cdot\psi\r)dt
$$
By using the generalized K\"ahler identity, 
we have 
$$
\pi_{U^{-n}}\ol\pa\pa(\tr\log {\overset{\scriptscriptstyle\bullet}{f_t}})\psi =\trian \tr\log {\overset{\scriptscriptstyle\bullet}{f_t}}\psi,
$$
where $4\trian$ is the ordinary Laplacian, which does not depend on $h.$
Since $\tr\log {\overset{\scriptscriptstyle\bullet}{f_t}}$ is a real function, $\pi^{\Herm_{f_t}}\trian \tr\log {\overset{\scriptscriptstyle\bullet}{f_t}}\psi
=\trian \tr\log {\overset{\scriptscriptstyle\bullet}{f_t}}\psi.$
Thus we have 
$$
\tr\, Q=\int^1_0 (\trian \tr\log {\overset{\scriptscriptstyle\bullet}{f_t}}\psi) dt=\trian\tr\log f
$$
Hence we obtain 
the result.
\end{proof}
\bgn{lemma}\label{for everyHermitian}
For every Hermitian metric $h_0$, there exists an Hermitian metric $h_f$ such that 
$$\tr\,\h\K_f(\psi)=0.$$
\end{lemma}
\bgn{proof}
Since 
$$
\int_M\lan \tr\h\K_{f_0}(\psi), \ol\psi\ran_s=0, 
$$
we have a solution $f$ such that 
$$
 \trian\tr\log f+\tr\,\K_{f_0}(\psi)=0
$$
Thus we have $\tr\,\h\K_f(\psi)=0.$
\end{proof}
Hence it follows from Lemma \ref{for everyHermitian}
that we can choose $h_0$ which satisfies $\tr\h\K_0(\psi)=0$. 

Then the trace of the equation (\ref{eq:continuity method1}) is written as 
$$
\tr\, (\K_{f}(\psi)-\lam \id_E)+\trian\tr\log f+\e\tr\,\log f =0
$$
Since $\h\K_f(\psi):=(\K_{f}(\psi)-\lam \id_E),$
then we have the following equation:
\bgn{equation}\label{trhK0(psi)}
\tr\, \h\K_f(\psi)+\trian\tr\log f+\e\tr\,\log f =0
\end{equation}
\bgn{proposition}
If $f$ satisfies {\rm(\ref{eq:continuity method1})} for $\e>0$, then $\tr\log f=0,$ that is, $\det f =1$.
\end{proposition}
\bgn{proof}
Since $f$ satisfies (\ref{trhK0(psi)}) and $\tr\,\h\K_0(\psi)=0,$ we have 
$$
\trian\tr\log f+\e\tr\,\log f =0
$$
Since $\trian +\e$ is invertible for $\e>0$, it follows $\tr\log f=0$. 
Then $\det f=1.$
\end{proof}
\bgn{proposition}
There exists an Hermitian metric $h_0$ and $f_1\in \Herm^+(E, h_0)$ such that 
$f_1$ satisfies (\ref{eq:continuity method1}) for $\e=1$ and $\tr\,\h\K_0(\psi)=0$
\end{proposition}
\bgn{proof}
From Lemma \ref{for everyHermitian}, we have an Hermitian metric $h$ such that 
$\tr\h\K_h(\psi)=0$. Let $\h\K_h(\psi)=\h\K_h\psi,$ where $\h\K_h\in \Herm(E, h_0).$
Then we define an Hermitian metric $h_0$ by 
$$
h_0 =h ( e^{\h\K_h}\,\,,\,\,)
$$
Let $f_1:=e^{-\h\K_h}.$
Then we have 
$$
\h\K_h+\log f_1=0
$$
Since $h=h_0(e^{-\h\K_h},\,)$ and $\h\K_h(\psi)=\h\K_{f_1}(\psi),$
we have a $h_0$ and $f$ satisfying (\ref{eq:continuity method}) for $\e=1.$
Since $\tr\log f_1=-\tr\,\h\K_h=0,$ we have 
$$
\tr \, \h\K_{f_1}(\psi)-\tr\,\h\K_{h_0}(\psi) =\trian\tr\log f_1=0
$$
Thus we see $\tr\,\h\K_{h_0}(\psi)=0.$
\end{proof}
\subsection{Step 1}
For a smooth section $f\in \Herm(E, h_0)$ and an integer $s\geq 0$, we define 
the norm $\|s\|_{L^2_s}^2$ by 
$$
\|f\|^2_{L^2_s}:=\sum_{i=0}^s\int_M |\nab^i f|^2\vol_M,
$$
where $|\nab^sf|$ denotes the point-wise norm of the $i$-th covariant derivative.
We denote by $H_s$ the Sobolev space $L_s^2\Herm(E, h_0)$ which is the completion of 
$\Herm(E, h_0)$ with respect to the Sobolev norm $\|\,\|_{L^2_s}$.
Let $L(f, \e)$ be the equation (\ref{eq:continuity method1})
$$
L(f, \e):=\K_f(\psi)-\lam (\id_E)+(\e\log f)=0
$$
Since $L(f, \e)\in \Herm(E, h_f),$ it follows from (\ref{piHermhf}) that 
$$
\Ad_{f^{\frac12}}L(f,\e)\in \Herm(E, h_0).
$$
Then
we define an operator $F: H_{s+2}\times [0,1]\to H_s$ by 
$$F(f, \e):=\Ad_{f^{\frac12}}L(f,\e),$$
where we identify $H_s\cdot\psi$ with the Sobolev space $H_s.$
Let $f$ be a solution of $L(f,\e)=0.$
Then the differential $DF$ of $F$ along $f$-direction at $(f,\e)$ is given by 
$$
D F(f, \e)(\xi)=\frac{d}{ds}F(fe^{\xi s}, \,\e)|_{s=0},
$$
where $\xi\in H_{s+2}.$
Since $L(f, \e)=0$ and $L(f, \e)\in \Herm(E, h_0),$ 
applying (\ref{piHermhf})and Proposition \ref{fracddsKhs|t=0trianh}, we have
\bgn{align}
D F(f, \e)(\xi) =&\pi^{\Herm_{h_0}}\Ad_{f^{\frac12}}\(\frac{d}{ds}\K_{f_s}(\psi)+\e D(log f)(\xi)\)+[\Dot{f}^{\frac12}, L(f, \e)]\\
=&\Ad_{f^{\frac12}}\pi^{\Herm_{h_{f_s}}}\(\frac{d}{ds}\K_{f_s}(\psi)+\e D(log f)(\xi)\)\\
=&\Ad_{f^{\frac12}}\(\trian_f(\xi\cdot\psi)+\e D(log f)(\xi)\),
\end{align}
where $\trian_f:=\trian_{h_f}$ is the Laplacian as in Proposition \ref{fracddsKhs|t=0trianh}.

In order to show that ${\mathcal S}$ is a open set, we shall show that the operator $\trian_f+\e D(\log f): H_{s+2}\to H_s$ is an isomorphism. At first we need the following Proposition in order to show that
the derivative $D(\log f)$ gives an injective map.
\bgn{proposition}\label{ifDlogfxi=0}
If $h_f(D(\log f)(\xi), \xi )=0,$ then $\xi =0.$
\end{proposition}
\bgn{proof}
Let $\{f_t\}$ be a smooth family of $\Herm(E, h_0)$.
We shall diagonalize $f_t$ by using a unitary basis of $E$ on an open dense subset $W.$
It is already know that 
there exits an open dense set $W$ in $M$ such that there is a unitary basis of $\{s_i(t)\}_{i=1}^r$ of $E|_W$ with respect to 
$h_t:=h_{f_t}$ satisfy the followings:   $s_i(t)$ are smoothly depending on $t,$
and $f_t$ is written as 
$$
f_t=\sum_i e^{\lam_i(t)}s_i(t)\otimes s_i(t)^*,
$$
where $\{s_i(t)^*\}_{i=1}^n$ denotes the dual basis of $\{s_i(t)\}_{i=1}^r.$
(This local diagonalization method is well-explained in Appendix 7.4 \cite{Lu-Tel_1995}). (see Proposition \ref{local diagonalization}).

Then $\log f_t$ is given by 
$$
\log f_t=\sum_i {\lam_i(t)}s_i(t)\otimes s_i(t)^*.
$$
Then $D(\log f)(\xi)$ is 
$$
D(\log f)(\xi)=\frac{d}{dt}\log f_t|_{t=0},
$$
where $\xi =f_t^{-1}\dot {f_t}.$
The derivative of $\{s_i(t)\}$ is given by 
$\frac{d}{dt}s_i(t) =\sum_{i,j }P_{i,j}s_j(t)$ for $P=(P_{i,j}).$
Then we have 
\bgn{align}
D(\log f)(\xi)=&\sum_i\lam_i'(0)s_i(0)\otimes s_i(0)^* +[P, \log f]\\
=&\sum_i\lam_i'(0)s_i(0)\otimes s_i(0)^*+\sum_{i\neq j}P_{i,j}(\lam_j(0)-\lam_i(0))s_i(0)\otimes s_j(0)^*
\end{align}
Since $\xi=f_t^{-1}{\overset{\scriptscriptstyle\bullet}{f_t}}$, $\xi$ is given by 
$$
\xi=\sum_i\lam_i'(0)s_i(0)\otimes s_i(0)^*+\sum_{i\neq j}P_{i,j}(e^{\lam_j(0)-\lam_i(0)}-1)s_i(0)\otimes s_j(0)^*
$$
We denote by $x_{i,j}$ the real number $(\lam_j(0)-\lam_i(0))$.
Then we have 
\bgn{align}
h_f (D(\log f)(\xi), \,\,\xi)=&\sum_i|\lam'_i(0)|^2+\sum_{i\neq j}|P_{i,j}|^2x_{i,j}(e^{x_{i,j}}-1)
\end{align}
Since $x_{i,j}(e^{x_{i,j}}-1)\geq0$, $h_f (D(\log f)(\xi), \,\,\xi)=0$ implies that $\xi=0.$

\end{proof}

Then we have 
\bgn{proposition}
The operator $\trian_f+\e D(\log f) : H_{s+2}\to H_s$ is an isomorphism.
\end{proposition}
\bgn{proof}
The operator $\trian_f+\e D(\log f)$ is an elliptic operator of index $0.$
Thus it suffices to show that the operator $\trian_f+\e D(\log f)$ is injective.
We have 
\bgn{align}
\int_M h_f((\trian_f+\e D(\log f))\xi\psi, \,\,\xi\ol\psi)=&\int_Mh_f(\trian_f\xi\psi , \,\,\xi\ol\psi)+\e h_f(D(log f)(\xi)\psi, \,\,\xi\ol\psi)\\
=&\int_Mh_f((\D^{1,0}_{h_f\, -})\xi \psi, \, (\D^{1,0}_{h_f\, -})\xi\ol\psi)+\e h_f(D(\log f)\xi\psi, \,\,\xi\ol\psi)\\
=&\|(\D^{1,0}_{h_f\, -})\xi\|_{L^2}^2+\int_M\e h_f(D(\log f)\xi\psi, \,\,\xi\ol\psi).
\end{align}
Thus $(\trian_f+\e D(\log f))\xi=0$ implies that $h_f(D(\log f)\xi, \,\,\xi)=0.$
From Proposition \ref{ifDlogfxi=0}, we have $\xi =0$
Hence $\trian_f+\e D(\log f) : H_{s+2}\to H_s$ is injective.
Since $\trian_f+\e D(\log f) : H_{s+2}\to H_s$, is an elliptic operator whose index is $0.$
Then the result follows.
\end{proof}
\bgn{proposition}
${\mathcal S}\subset (0,1]$ is an open set.
\end{proposition}
\bgn{proof}
From the implicit function theorem, we obtain that 
the region ${\mathcal S}$ is open set in $(0,1].$
\end{proof}
\subsection{Preliminary results for the Step 2 and Step 3}

An Hermitian metric $h$ on $E$ yields a bilinear form $ h$ of 
$\End(E)\otimes \w^\bullet T^*_M$ by
$$
 h(A\otimes \a, \,\, B\otimes \b) := i^n\Re\{\,\tr (A B^{*_h})(\a\w\sig\ol\b)\}\in \w^\bullet T^*_M,
$$
where $A, B\in \End(E)$, $B^{*_h}$ is the adjoint of $B$ with respect to $h$
and $\a,\b\in \w^\bullet T^*_M$ and $\sig$ is the Clifford involution.
We denote by
$h^{top}$ the $2n$-component of $h$, i.e., the form of top degree. 
Since $(\log f)\psi\in\Herm_f\cap U^{-n}$, we have
\bgn{lemma}\label{htoppiHermftU-n}
$$
h^{\top}( \pi^{\Herm_{f_t}}_{U^{-n}}\ol\pa^E
(\dot{A}_{f}\psi), \,\, (\log f)\psi)=h^{\top}( \ol\pa^E
(\dot{A}_{f}\psi), \,\, (\log f)\psi)
$$
\bgn{proof}
The result follows from Lemma \ref{Herm ST ST}, applying 
$T_f=\log f.$

\end{proof}

\end{lemma}
The following three Propositions are important for Step.2 and Step.3.
Let $\a=f\,\vol_M$ and $\b=g\,\vol_M$ be two differential forms of top degree.
Then if $f>g$, we denote it by $\a>\b$.
\bgn{proposition}\label{prop:htop(Q}
Let $h_t=h_{f_t}$ be a family of Hermitian metrics with $f:=f_1$. 

We denote by $Q$ the difference $\K_{f}(\psi)-\K_{f_0}(\psi).$
Then we have
$$h^{top}_0(Q, (log f)\psi)\geq \frac12 \trian |\log f|^2\vol_M,
$$
\end{proposition}
 \bgn{proof}
 From Proposition \ref{Kftpsiint} and $f_1=f$, 
 we have 
\bgn{align} 
h^{\top}(Q, (log f)\psi)=&\int^1_0h^{\top}(\ol\pa^E(\overset{\,\,\,\bullet}{\A}_{f_t}\psi) dt, \,\, (\log f)\psi)dt\\
=&h^{\top}(\ol\pa^EA_{f_1}\psi, \,\,(\log f)\psi)-
h^{\top}(\ol\pa^E A_{f_0}\psi, \,\,(\log f)\psi)\\
=&h^{\top}(\ol\pa^E(f^{-1}(\pa_0 (f\psi))), \,\,(\log f)\psi)
\end{align}
Then we have
\bgn{equation}\label{htopQlogfpsipitopolpaf}
h^{\top}(Q, (log f)\psi)=\pi^{\top}\ol\pa_\J h((f^{-1}(\pa_0 (f\psi))), \,\,(\log f)\psi)+h^{\top}((f^{-1}(\pa_0 (f\psi))), \,\,\ol\pa^E(\log f)\psi)
\end{equation}
We shall use the method of "local diagonalization"\cite{Lu-Tel_1995}
again.
\bgn{proposition}\label{local diagonalization}
Let $(E,h)$ be an Hermitian vector bundle and $f$ an Hermitian endmorphism with respect to $h.$
Then there exists an open dense subset $W\subset M$ such that for every $x\in W$  the following holds: 
There exists an open neighborhood $U$ of $x$ a unitary basis $\{s_i\}$ for $E$ defined over $U$ and functions 
$\lam_i\in C^\infty(U)$ such that 

$$
f(y) =\sum_{i=1}^r e^{\lam_i(y)}s_i(y)\otimes s^i(y),
$$
where $\{s^i(y)\}$ denotes the dual basis and $y\in U,$
that is, $f$ is given in the following form of the diagonal matrix :

\bgn{equation}\label{pmatrix f}
f=
\bgn{pmatrix}
e^{\lam}_1&\cdots&0\\
0&e^{\lam_2}&\cdots\\
&\vdots& \\
0& \cdots & e^{\lam_r}
\end{pmatrix}
\end{equation}
\end{proposition}
Then we have
\bgn{align}
\pa_0 (f\psi))=&\pa(f\psi)+[\A_0 , f\psi]\\
=&\sum_i \pa(\lam_i\psi)e^{\lam_i}s_i\otimes s_i^* +\sum_{i,j}\A_{i,j}\psi(e^{\lam_j}-e^{\lam_i})s_i\otimes s_j^*\\
\\
f^{-1}\pa_0 (f\psi))=&
\sum_i \pa(\lam_i\psi)s_i\otimes s_i^* +\sum_{i,j}\A_{i,j}\psi(e^{\lam_j-\lam_i}-1)s_i\otimes s_j^*\\
\\
\pa_0(\log f\psi)=&\pa \lam_i s_i\otimes s_i^*+[\A, \log f]\psi
\end{align}
Thus we have 
\bgn{align}
h(f^{-1}\pa_0 (f\psi), \,\, \log f \psi) =
\sum_i \lam_i(\pa \lam_i)\psi\w\sig\ol\psi+ h( f^{-1}[\A^{1,0},f]\psi, \,\, \log f \psi)
\end{align}
Since $\tr( f^{-1}[\A^{1,0},f]\log f)=0,$ the first term of the right hand side of (\ref{htopQlogfpsipitopolpaf})
is given by 
$$
\pi^{\top}\ol\pa h(f^{-1}\pa_0 (f\psi), \,\, \log f \psi) =
\sum_i \ol\pa(\lam_i(\pa \lam_i))\vol_M
=\frac12\trian |\log f|^2 \vol_M
$$
We shall show that 
the second term of the R.H.S of (\ref{htopQlogfpsipitopolpaf})
$h((f^{-1}(\pa_0 (f\psi))), \,\,\ol\pa^E(\log f)\psi)
$ is greater than or equal to $0.$
Applying the local diagonalization, we have  
\bgn{align}
h((f^{-1}(\pa_0 (f\psi))), \,\,\ol\pa^E(\log f)\psi)
=\sum_i h(\pa \lam_i \psi , \,\, \pa\lam_i \psi)
+h (f^{-1}[\A^{1,0}, \, f]\psi, \,\, [\A^{1,0},\, \log f]\psi)
\end{align}
Since 
$h(\pa \lam_i \psi, \,\,\pa\lam_i \psi) =\lan \pa \lam_i \psi, \,\,\ol\pa\ol\psi\ran_s\geq 0$, thus the first therm is non-negative.
The last term is also nonnegative from the following lemma:
\bgn{lemma}\label{h(f-1[A1,0 f]psi]}
$$
h (f^{-1}[\A^{1,0}, \, f]\psi, \,\, [\A^{1,0},\, \log f]\psi)
\geq 0
$$
\end{lemma}
\bgn{proof}
The result follows from the direct calculation by using the method of local diagonalization.
\end{proof}
Then the second term of the R.H.S of (\ref{htopQlogfpsipitopolpaf})$\geq 0.$
Thus Proposition \ref{prop:htop(Q}
follows. 
\end{proof}

\bgn{proposition}\label{forsig>0wehave} For $\sig >0$, we have 
$$\pi_{\top}d h( f^{-1}(\pa_0 f\psi)), \, f^\sig \psi) =\frac1\sig \lan\trian(\tr\,f^\sig \psi), \,\,
\ol\psi\ran_s,
$$
Where $\pi_{\top}$ denotes the projection to $\w^{2n}T^*_M.$
Note that $\trian:=\trian_{\del_-}$ of $(M, \J, \J_\psi)$ does not depend on $h_f$.
\end{proposition}
\bgn{proof}
By using the generalized K\"ahler identity, we obtain
$$
\frac1{\sig}\lan \trian (\tr f^\sig \psi), \,\, \ol\psi\ran_s=\frac1{\sig}\lan \ol\pa_\J\pa_\J (\tr f^\sig \psi), \,\,\ol\psi\ran_s
$$
We see that 
 $$\tr\,f^{\sig-1}[\A^{1,0}, f]=\tr\,[\A^{1,0}, f^\sig]=0$$
we have 
$$
\tr(f^{\sig-1}\pa_0 f) =\tr(f^{\sig-1}\pa f)+\tr f^{\sig-1}[\A^{1,0}, f]=\tr(f^{\sig-1}\pa f)
=\frac1\sig \tr (\pa_0 f^\sig)
$$
The left hand side is given by 
\bgn{align}
\pi_{\top}d\, h( f^{-1}(\pa_0 f\psi)), \, f^\sig \psi)=&
\pi_{\top}d\,h( f^{\sig-1}(\pa_0 f\psi)), \,  \psi)\\
=&\pi_{\top}  h( \ol\pa(f^{\sig-1}(\pa_0 f\psi))), \,  \psi)\\
=& h( \pi_{U^{-n}}(\ol\pa f^{\sig-1}(\pa_0 f\psi))), \, \,\, \psi)\\
=&\tr\lan \ \ol\pa (f^{\sig-1}(\pa_0 f\psi))), \,\, \ol\psi\ran_s\\
=&\frac1\sig \tr\lan \ol\pa\pa_0f^\sig\psi , \,\,\ol\psi\ran_s\\
=&\frac1\sig\lan \ol\pa_\J\pa_\J \tr f^\sig\psi, \,\,\ol\psi\ran_s
\end{align}
Thus we obtain 
$$\pi_{\top}d\, \h h( f^{-1}(\pa_0 f\psi)), \, f^\sig \psi) =\frac1\sig \lan\trian(\tr\,f^\sig \psi), \,\,
\ol\psi\ran_s,
$$
\end{proof}

\bgn{proposition}\label{prop:Let f in Herm}
Let $f\in \Herm(E, h_0)$ and $0<\sig \leq 1$.
Then we have the followings: 
\bgn{enumerate}
\item[(1)]
$$h_0^{\top}(f^{-1}(\pa_0 f\psi), \,\,\pa_0(f^\sig\psi)) \geq 
|f^{\frac{-\sig}2}(\pa_0 f^{\sig} \psi)|^2_{h_0}\vol_M
$$
\item[(2)]
If $f$ satisfies the equation (\ref{eq:continuity method1}) for some $\e>0$, then 
we have 
$$
\frac1\sig \trian (\tr\, f^\sig)\vol_M+ \e h^{\top}_0((\log f)\psi, \,\, f^\sig \psi)+|f^{\frac{-\sig}2}(\pa_0 f^{\sig} \psi)|^2_{h_0}\vol_M
\leq -h^{\top}_0(\h \K_0(\psi), \,\, f^\sig \psi)
$$
\end{enumerate}
\end{proposition}
\bgn{proof} 
(1) We use the same notations as before. 
$$
\pa_0(f^\sig\psi) =\pa(f^\sig\psi)+[\A^{1,0}, f^\sig]\psi
$$
Then the left hand side of (1) is given by 
\bgn{align}
h_0^{\top}(f^{-1}(\pa_0 f\psi), \,\,\pa_0(f^\sig\psi))=&
h_0^\top(f^{-1}(\pa f\psi)+f^{-1}[\A^{1,0}, f]\psi, \,\,\pa (f^\sig\psi)+[\A^{1,0}, f^\sig]\psi)\\
=&h_0^\top(f^{-1}(\pa f\psi), \,\,\pa (f^\sig\psi))+h_0^\top(f^{-1}(\pa f\psi),[\A^{1,0},\,\, f^\sig]\psi)\label{h0topf-1}\\
+&h_0^\top(f^{-1}[\A^{1,0}, f]\psi, \,\,\pa (f^\sig\psi))+h_0^\top(f^{-1}[\A^{1,0}, f]\psi, \,\,[\A^{1,0}, f^\sig]\psi)
\label{A10f}
\end{align}
Then we see 
\bgn{align}
&h_0^\top(f^{-1}(\pa f\psi),[\A^{1,0},\,\, f^\sig]\psi)=0\\
&h_0^\top(f^{-1}[\A^{1,0}, f]\psi, \,\,\pa (f^\sig\psi))=0
\end{align}
The right hand side of (1) is given by 
\bgn{align}
|f^{\frac{-\sig}2}(\pa_0 f^{\sig} \psi)|^2_{h_0}\vol_M
=&h_0^\top(f^{-\sig}\pa_0 (f^\sig\psi), \,\,\pa_0(f^\sig\psi))\\
=&h_0^\top(f^{-\sig}\pa (f^\sig\psi), \,\,\pa(f^\sig\psi)\label{h0topf-sig}\\
+&h_0^\top(f^{-\sig}[\A^{1,0}, f^\sig]\psi, \,\,[\A^{1,0}, f^\sig]\psi)\label{A10fsig}
\end{align}
The first term of (\ref{h0topf-1}) and (\ref{h0topf-sig}) are written as 
\bgn{align}
&h_0^\top(f^{-1}(\pa f\psi), \,\,\pa (f^\sig\psi))=\sig h_0^\top( f^{\sig-2}\pa(f\psi), \,\, \pa (f\psi))\\
&h_0^\top(f^{-\sig}\pa_0 (f^\sig\psi), \,\,\pa_0(f^\sig\psi))=\sig^2h_0(f^{\sig-2}\pa(f\psi), \,\,\pa(f\psi))
\end{align}
since $f\in \Herm(E).$
Since $0\leq \sig \leq1, $ we have 
$$
h_0^\top(f^{-1}(\pa f\psi), \,\,\pa (f^\sig\psi))\geq h_0^\top(f^{-\sig}\pa_0 (f^\sig\psi), \,\,\pa_0(f^\sig\psi)).
$$
In order to estimate the second term of (\ref{A10f}) and (\ref{A10fsig}), we diagonalize $f$ at each point of $M$
as in (\ref{pmatrix f}).
Let $\A^{1,0}=\sum_k \A^{1,0}_{(k)}e_k$, where $\{e_k\}\in \L_\J$ denotes the orthonormal basis.
Then $(i,j)$-components of the matrix $[\A^{1,0}_{(k)}, f]$, $f^{-1}[\A^{1,0}_{(k)}, f]$ and $[\A^{1,0}_{(k)}, f^\sig]$ are given by 
\bgn{align}
[\A^{1,0}_{(k)}, f]_{i,j}=&a_{i,j}(e^{\lam_j}-e^{\lam_i)}\\
f^{-1}[\A^{1,0}_{(k)}, f]_{i,j}=&a_{i,j}(e^{\lam_j-\lam_i}-1)\\
[\A^{1,0}_{(k)}, f^\sig]_{i,j}=&a_{i,j}(e^{\sig\lam_j}-e^{\sig\lam_i})
\end{align}
Since $[\A^{1,0}_{(l)}, f^\sig]$ is Hermitian, 
\bgn{align}
h_0^\top(f^{-1}[\A^{1,0}, f]\psi, \,\,[\A^{1,0}, f^\sig]\psi)&=\sum_{k,l}\tr (f^{-1}[\A^{1,0}_{(k)}, f]\,[\A^{1,0}_{(l)}, f^\sig])\lan e_k\cdot\psi, \,\,\ol e_l\cdot\ol\psi\ran_s\\
&=\sum_k 2\tr (f^{-1}[\A^{1,0}_{(k)}, f]\,[\A^{1,0}_{(k)}, f^\sig])\vol_M, 
\end{align}
since $i^n\lan e_k\cdot\psi, \,\,\ol e_l\cdot\ol\psi\ran_s=2\del_{k,l}\vol_M.$
Then we have 
\bgn{align}\label{h0topf-1A10k}
h_0^\top(f^{-1}[\A^{1,0}_{(k)}, f]\psi, \,\,[\A^{1,0}_{(k)}, f^\sig]\psi)=&
\sum_{i,j}(e^{\lam_j-\lam_i}-1)(e^{\sig\lam_j}-e^{\sig\lam_i})| a_{i,j}|^22\vol_M
\end{align}
since $a_{i,j}=-a_{j,i}^*.$
(\ref{A10fsig}) is also given by 
\bgn{align}
h_0^\top(f^{-\sig}[\A^{1,0}, f^\sig]\psi, \,\,[\A^{1,0}, f^\sig]\psi)=
&\sum_{k}\tr(f^{-\sig}[\A^{1,0}_{(k)}, f^\sig] \,[\A^{1,0}_{(k)}, f^\sig])2\vol_M
\end{align}
We also have 
\bgn{align}\label{sumktrf-sig}
\sum_{k}\tr(f^{-\sig}[\A^{1,0}_{(k)}, f^\sig] \,[\A^{1,0}_{(k)}, f^\sig])2\vol_M=\sum_{i,j}(e^{\sig\lam_j -\sig\lam_i} -1)(e^{\sig\lam_j}-e^{\sig\lam_i})|a_{i,j}|^2
2\vol_M
\end{align}
Let $\lam=\lam_j-\lam_i. $
Since $(e^\lam -1)(e^{\sig\lam_j}-e^{\sig\lam_i})\geq (e^{\sig\lam}-1)(e^{\sig\lam_j}-e^{\sig\lam_i})$ for $0\leq \sig \leq 1,$ it follows from (\ref{h0topf-1A10k}
) and (\ref{sumktrf-sig})that, 
$$
h_0^\top(f^{-1}[\A^{1,0}_{(k)}, f]\psi, \,\,[\A^{1,0}_{(k)}, f^\sig]\psi)\geq h_0^\top(f^{-\sig}[\A^{1,0}, f^\sig]\psi, \,\,[\A^{1,0}, f^\sig]\psi)
$$
Hence we have 
$$h_0^{\top}(f^{-1}(\pa_0 f\psi), \,\,\pa_0(f^\sig\psi)) \geq 
|f^{\frac{-\sig}2}(\pa_0 f^{\sig} \psi)|^2_{h_0}\vol_M
$$

Proof of (2) :
Since $f$ satisfies the equation (\ref{eq:continuity method}) for some $\e>0$, we have 
$$
\h\K_f(\psi)-\h\K_{0}(\psi) +\h\K_{0}(\psi)+\e(\log f)\psi =0
$$
Then we have 

\bgn{align}
h_0^\top(\h\K_f(\psi)-\h\K_{0}(\psi), \,\,f^\sig\psi)+\e h_0^\top((\log f)\psi, \,\, f^\sig\psi) =-h_0^\top(\h\K_{0}(\psi), \,\, f^\sig\psi)
\end{align}

From Proposition \ref{Kftpsiint}, we obtain
$$
h_0^\top(\h\K_f(\psi)-\h\K_{0}(\psi), \,\, f^\sig \psi ) =h_0^\top(\ol\pa^E(f^{-1}(\pa_0 (f\psi))), f^\sig\psi)
$$
From Proposition \ref{forsig>0wehave} and Proposition \ref{prop:Let f in Herm} (1), we obtain 
\bgn{align}
h_0^\top(\ol\pa^E(f^{-1}(\pa_0 (f\psi))), f^\sig\psi)
=&\pi_{\top}d {h}_0 (f^{-1}\pa_0(f\psi), \,\,f^\sig\psi)+h_0^\top (f^{-1}\pa_0(f\psi), \,\,\pa_0(f^\sig\psi))\\
\geq& \frac1\sig \lan\trian(\tr\,f^\sig \psi), \,\,
\ol\psi\ran_s,
+|f^{\frac{-\sig}2}(\pa_0 f^{\sig} \psi)|^2_{h_0}\vol_M
\end{align}
Hence we obtain 
$$
\frac1\sig \lan\trian(\tr\,f^\sig \psi), \,\,
\ol\psi\ran_s,
+|f^{\frac{-\sig}2}(\pa_0 f^{\sig} \psi)|^2_{h_0}\vol_M+
\e h_0^\top((\log f)\psi, \,\, f^\sig\psi) \leq -h_0^\top(\h\K_{f_0}(\psi), \,\, f^\sig\psi)
$$
\end{proof}

\subsection{Step 2}
Step 2 is divided into the following three steps.
In Step 2-1, we shall show a $C^0$ estimate of $\log f_\e$, i.e., $\e\|\log f_\e\|_{C^0}<C$.
In Step 2-2, we shall show a $C^1$-estimate of $f_\e$ in terms of $m_\e:=\|\log f_\e\|_{C^0}$.
 In Step 2-3, we shall obtain an $L_2^p$-estimate of $f_\e$  for all $p>0$.
\subsubsection{Step 2-1}
We need the following well-known result through this section (c.f.\cite{Lu-Tel_1995}). 
\bgn{lemma}\label{trian flea lam f+mu}
Let $f\in C^2(M,\R)$ be a function satisfying $f\geq 0$ and 
$$
\trian f\leq \lam f+\mu,
$$
for $0\leq \lam \in \R$ and $\mu \in \R.$
Then there is a positive constant $C$ depending only on $(M, \J, \J_\psi)$ and $\lam$ such that 
$$
\|f\|_{C^0}\leq C(\|f\|_{L^1}+\mu)
$$
\end{lemma}
\bgn{proposition}\label{LetfeinHerm+(E, h0)}
Let $f_\e\in \Herm^+(E, h_0)$ be a solution of $L_\e(f)\psi=0$ for some $\e>0.$
Then we have the followings:
\bgn{enumerate}
\item[(1)]
$$
\frac12\trian (|\log f_\e|^2)+\e|\log f_\e|^2\leq |\h \K_0||\log f_\e|
$$
\item[(2)]Let $m_\e$ be the $C^0$-norm of $\log f_\e$. Then
$$m_\e\leq \frac1\e \|\h\K_0\|_{C^0}$$
\item[(3)] $${\|f_\e\|_{C^0}\leq C( \|\log f_\e\|_{L^2}+\|\h \K_0\|_{C^0})}$$
\end{enumerate}
\end{proposition}
\bgn{proof}
Since $f_\e$ is a solution of $L_\e(f_\e)=0,$ we have
$$
\K_{f_\e}(\psi)-\K_0(\psi)+ \h\K_0(\psi) +\e(\log f_\e)=0
$$
Then by using  $\log f_\e$ and $h_0$, we have 
$$
h_0 (\K_{f_\e}(\psi)-\K_0(\psi)+\e(\log f_\e), \,\,\log f_\e)=-h_0(\h\K_0(\psi), \log f_\e) 
$$
Then using the Schwartz inequality, we obtain 
$$ h_0(\K_{f_\e}(\psi)-\K_0(\psi), \,\log f_\e)+\e |\log f_\e|^2
\leq |h_0(-\h\K_0, \,\log f_\e)|
\leq
|\h\K_0|\,|\log f_\e| 
 $$
From Proposition \ref{prop:htop(Q}, we have 
$$
\frac12\trian (|\log f_\e|^2)
\leq
h_0(\K_f(\psi)-\K_0(\psi), \,\log f_\e)$$
Then we obtain (1)
$$
\frac12\trian (|\log f_\e|^2)+\e |\log f_\e|^2\leq
|\h K_0|\,|\log f_\e|
$$
(2) follows from (1) by using the Maximum principle.
In fact, Let $x\in M$ be a point such that $|\log f_\e|^2$ attaints its maximum. 
Then $\trian (|\log f_\e(x)|^2)\geq 0$. 
It follows from (1) that 
$$\e |\log f_\e(x)|^2\leq|\h \K_0|\,|\log f_\e(x)|$$
Thus we have $m_\e=\|\log f_\e\|_{C^0}\leq \e^{-1}\|\h\K^0\|_{C^0}$
From (1), we have 
$$
\trian (|\log f_\e|^2)\leq 2|\h \K_0|\,|\log f_\e|\leq |\log f_\e|^2+|\h\K_0|^2\leq |\log f_\e|^2+\|\h\K_0\|^2_{C^0}
$$
where $\|\h K_0\|_{C^0}$ is a constant. 
Applying Lemma \ref{trian flea lam f+mu}, we obtain 
$$
\|(\log f_\e)^2\|_{C^0}\leq C(\|\log f_\e\|^2_{L^2}+\|\h \K_0\|^2_{C^0})
$$
Thus we have (3)
$$
\|\log f_\e\|_{C^0}\leq C(\|\log f_\e\|_{L^2}+\|\h \K_0\|_{C^0})
$$
\end{proof}
\subsubsection{Step 2-2}
We use the following notations :
\bgn{align}
m_\e:=&\|\log f_\e\|_{C^0}\\
\zeta_\e:=&\frac{d}{d\e}f_\e\\
\eta_\e:=&f_\e^{-\frac12}\circ \zeta_\e\circ f^{-\frac12}\\
\h\eta:=&f_\e^{-1}\frac{d}{d\e}f_\e
\end{align}
Since $\det f_\e=1,$ we have $C \leq|\log f_\e|$, where 
$C$ are constant depending only on $m_\e$.
In this section, $C(m_\e)$ means a positive constant which depends only on $m_\e$.
Note that $C(m_\e)$ may be changed as we proceed.
In order to obtain results of Step 2-2, we need the following inequality :

\bgn{proposition}\label{trianetae2etaedf}
$$\trian(|\eta_\e|^2)+2\e|\eta_\e|^2+|d^f\eta_\e|^2\leq -2h_0(\log f_\e, \, \eta_\e)
$$
\end{proposition}
\bgn{proof}
Since $L(f_\e, \e)=0,$ we have 
$$
\frac{d}{d\e}L(f_\e, \e)=0.
$$
Then as in Section 1 and subsection 2-1, the result follows from calculation.
\end{proof}
\bgn{lemma}\label{tretae=0}
$$\tr\, \eta_\e=0.$$
\end{lemma}
\bgn{proof}
The result follows from $\det f_\e=1$ for all $\e$ since 
$\tr\,\eta_\e =\tr\, f_\e^{-1}\frac{d}{d\e}f_\e.$\\
\end{proof}
\bgn{proposition}
For all $\e$, it holds 
$$
\|\zeta_\e\|_{C^0}:=\|\frac{d}{d\e}f_\e\|_{C^0}\leq C(m_\e)\e^{-1}
$$
\end{proposition}
\bgn{proof} From Proposition \ref{trianetae2etaedf}, we have 
$$\trian(|\eta_\e|^2)+2\e|\eta_\e|^2+|d^f\eta_\e|^2\leq -2h_0(\log f_\e, \, \eta_\e)
\leq 2|\log f_\e||\eta_\e|
$$
Thus we have 
$$\trian(|\eta_\e|^2)+2\e|\eta_\e|^2\leq -2h_0(\log f_\e, \, \eta_\e)
\leq 2|\log f_\e||\eta_\e|
$$
Integration over $M$ yields 
$$
\e\int_M|\eta_\e|^2\leq m_\e\int_M |\eta_\e|
$$
Thus we have 
$$
\e\|\eta_\e\|^2_{L^2}\leq m_\e\|\eta_\e\|_{L^1}\leq m_\e \Vol_M^{\frac12}\|\eta_\e\|_{L^2}
$$
Hence we obtain 
$$
\|\eta_\e\|_{L^2}\leq C(m_\e)\e^{-1}
$$
From Proposition \ref{trianetae2etaedf}, we also have 
$$\trian(|\eta_\e|^2)\leq -2h_0(\log f_\e, \, \eta_\e)
\leq 2|\log f_\e||\eta_\e|\leq 2m_\e|\eta_\e|\leq m_\e|\eta_\e|^2+m_\e
$$
Applying Lemma \ref{trian flea lam f+mu}, we have 
$$
\|\eta_\e\|_{C^0}\leq C(m_\e)(\|\,|\eta_\e|^2\|_{L^1}+m_\e)\leq 
C(m_\e)(\|\eta\|^2_{L^2}+m_\e)
$$
Since $\|\eta_\e\|_{L^2}<C(m_\e)\e^{-1}, $ 
we have $\|\eta_\e\|_{C^0}<C(m_\e)\e^{-1}$. 
Since $\eta_\e=f_\e^{\frac12}\circ \zeta_\e\circ f_\e^{\frac12}$, 
we have the result.
\end{proof}
\subsubsection{Step 2-3}
Let $\e_0$ be a constant in $(0,1).$
We set $m:=\dstyle{\frac{\|\K_0\|_{C^0}}{\e_0}}.$
Then $m_\e:=\|\log f_\e\|_{C^0}\leq m$ for all $\e\in(\e_0, 1]$.
We also have $C(m_\e)\e^{-1}\leq C(m),$ where $C(m)$ is a constant depending only on $m.$
We shall show that the following Proposition in the section. 
\bgn{proposition}\label{Letfeandzeta}
Let $f_\e$ and $\zeta_\e:=\frac{d}{d\e}f_\e$ be as before. 

Then for all $p>0$ and $\e\in (\e_0, 1]$, we have 
\bgn{enumerate}
\item[(1)]
$$
\|\zeta_\e\|_{L_2^p}\leq C(m)\(1+\|f_\e\|_{L_2^p}\)
$$
 
\item[(2)] 
$$
\|f_\e\|_{L_2^p}\leq e^{C(m)(1-\e)}\(1+\|f_1\|_{L_2^p}\)
$$
\end{enumerate}
\end{proposition}
\bgn{proof}
(1) From Proposition \ref{fracddsKhs|t=0trianh}, the differential of $\h\K_{f_\e}(\psi)$ is given by 
$$
\pi^{\Herm_{f_\e}}\frac{d}{d\e}\h\K_{f_\e}(\psi) =\trian_{f_\e}(\xi\psi), 
$$
where $\xi =f_\e^{-1}\zeta_\e$.
Since $\pi_{U^{-n}}\ol\pa^E_\J \pa_{f_\e}=\trian_{f_\e}$ for $U^{-n}$ and 
$\pa_{f_\e}=f_\e^{-1}\circ \pa_0 \circ f_\e$, from $ (\ref{Af=f-1paJf}),$
we have 
\bgn{align}
\trian_{f_\e}(\xi_\e\psi) =&\pi_{U^{-n}}\ol\pa^E \pa_{f_\e}(\xi\psi)=\pi_{U^{-n}}\ol\pa^E \circ f_\e^{-1}\circ \pa_0
\circ f_\e(\xi_\e\psi)\\
=&\pi_{U^{-n}}\ol\pa^E f_\e^{-1}\circ\pa_0 (\zeta_\e\psi)\notag\\
=&\pi_{U^{-n}}f_\e^{-1}\trian_0(\zeta_\e\psi)+\pi_{U^{-n}}(\ol\pa^E f_\e^{-1})(\pa_0 \zeta_\e\psi),\notag
\end{align}
where $\trian_0=\pi_{U^{-n}}\ol\pa^E\pa_0$ denotes the Laplacian with respect to $h_0.$
Since $L(f_\e, \e)=0$, we have 
\bgn{align}\label{fracddeLfe,e}
0=\frac{d}{d\e}L(f_\e, \e)=&\pi^{\Herm_{f_\e}}\frac{d}{d\e}\h\K_{f_\e}(\psi)+\e\frac{d}{d\e}(\log f_\e)(\psi)+(\log f_\e)(\psi)\\
=&f_\e^{-1}\trian_0(\zeta_\e\psi)+\pi_{U^{-n}}(\ol\pa^E f_\e^{-1})(\pa_0 \zeta_\e\psi)\notag\\
+&\e\frac{d}{d\e}(\log f_\e)(\psi)+(\log f_\e)(\psi)\notag
\end{align}
Since the operator $f^{-1}\trian_0+ \e D(\log f)$ is invertible, we have the estimate
\bgn{align}
\|\zeta_\e\|_{L_2^p}\leq& C\|f_\e^{-1}\trian_0(\zeta_\e \psi)+\e D(\log f_\e)(\zeta_\e\psi)\|_{L^p}\\
\end{align}
Since $\e D(\log f_\e)(\zeta_\e\psi)=\e\frac{d}{d\e}(\log f_\e)(\psi),$
it follows from (\ref{fracddeLfe,e}) that we have 
\bgn{align}
\|\zeta_\e\|_{L^p}\leq& C \|\pi_{U^{-n}}(\ol\pa^E f_\e^{-1})(\pa_0 \zeta_\e\psi)+(\log f_\e)(\psi)\|_{L^p}\\
\leq &C\|\pi_{U^{-n}}(\ol\pa^E f_\e^{-1})(\pa_0 \zeta_\e\psi)\|_{L_2^p}+\|(\log f_\e)(\psi)\|_{L^p}\\
\end{align}
Since $\|f_\e^{-1}\|_{C^0}$ is bounded, using the H\"older inequality, we have 
\bgn{align}
\|\zeta_\e\|_{L_2^p}\leq&C \|f_\e\|_{L_1^{2p}}\, \|\zeta_\e\|_{L_1^{2p}}+C
\end{align}
Then using the general formula as in the book by Aubin, we  have 
$$
\|\zeta_\e\|_{L_2^p}\leq C \|f_\e\|_{L_1^{2p}}\, \|\zeta_\e\|_{L_1^{2p}}+C
\leq C \|f_\e\|^{\frac12}_{L_2^{p}}\, \|\zeta_\e\|^{\frac12}_{L_2^{p}}+C
$$
Since $C x^{\frac12}y^{\frac12}\leq \frac12(C^2x+y),$ we have 

 $$
 \|\zeta_\e\|_{L_2^p}\leq \frac12 C^2 \|f_\e\|_{L_2^p}+\frac12\|\zeta_\e\|_{L_2^p}+C
 $$
Thus we obtain 
$$
\|\zeta_\e\|_{L_2^p}\leq C(m)\(1+\|f_\e\|_{L_2^p}\).
$$
(2) If $x(\e)$ is a smooth function with values in a normed vector space, then we have 
$$
\frac{d}{d\e}\|x(\e)\|\geq -\|\frac{d}{d\e}x(\e)\|.
$$
This follows from the triangle inequality. 
Applying this to $x(\e):=f_\e$ with the norm $\|\,\|_{L_2^p},$ we have 
$$
\frac{d}{d\e}\|f_\e\|_{L_2^p}\geq -\|\zeta_\e\|_{L_2^p}\geq -C(m)\(1+\|f_\e\|_{L_2^p}\)
$$
from (1).
Then we have 
\bgn{align}
\log(1+x(1))-\log(1+x(\e))=&\int^\e_1\frac{1}{1+x(s)} \frac{d}{ds}x(s) ds\\
\geq &-C(m)(1-\e)
\end{align}
Hence 
$$
\frac{1+x(1)}{1+x(\e)}\geq e^{-C(m)(1+\e)}
$$
Thus 
$\|f_\e\|_{L_2^p}\leq 1+\|f_\e\|_{L_2^p}\leq e^{C(m)(1+\e)}(1+\|f_1\|_{L_2^p}).$
\end{proof}
\bgn{proposition} We have the followings
\bgn{enumerate}
\item[(1)]
The set of solutions ${\mathcal S}$ contains $(0, 1].$
\item[(2)] If we have a uniform $C^0$-upper bound $\|f_\e\|_{C^0}<C$ for all $\e\in (0,1],$
then there is a solution $f_0$ of the equation $L(f_0, 0)=0,$
that is, $h_{f_0}$ is an Einstein-Hermitian metric.
\end{enumerate}
\end{proposition}
\bgn{proof}(1) From Proposition \ref{Letfeandzeta}, there is a sequence $\{\e_j\}$ such that 
$f_{\e_j}$ converges strongly $f_{\e_0}$ in $L_1^p$ for $p>2n$, when $\e_j\to \e_0.$
Then by using the Stokes theorem, we obtain that 
$$
(L(f_{\e_0},\e_0), \,\,\a)_{L^2}=0
$$
for any smooth section $\a\in \Herm(E).$
Then by using the elliptic regularity, it turns out that $f_{\e_0}$ is smooth which gives rise to an Einstein-Hermitian metric
$h_{f_{\e_0}}$. 
(2) is also shown by the same way as in (1).
\end{proof}
\subsection{Step 3 }
Let $f_\e$  be a solution of the equation for $\e>0$.
\bgn{equation}\label{eq:continuity method}
\K_{f_\e}(\psi)-\lam \id_E+\e\log f_\e=0
\end{equation}
We assume that $\|f_\e\|_{C^0}\to\infty$ $(\e\to 0).$
Let $\lam_\e(x)$ be the largest eigenvalue of $\log f_\e(x)$ at $x\in M.$
We denote by $M_\e$ maximum of $\lam_\e$, i.e.,
$M_\e:=\max\{\, \lam_\e(x)\, |\, x\in M\,\}$.
Set $\rho_\e=e^{-M_\e}.$
Since $\tr \log f_\e=0,$ we see that $M_\e=O(\e^{-1})$ from Lemma \ref{LetfeinHerm+(E, h0)}. We can assume $\rho_\e\leq 1.$
\bgn{proposition}
There is a sequence $\{\e_i\}$ which goes to $0$ $(i\to \infty)$ and 
$\rho(\e_i)$ also goes to $0$ such that $f_i:=\rho(\e_i)f_{\e_i}$ satisfies the following properties:
\bgn{enumerate}
\item[(1)] When $i$ goes to infinity, the $f_i$ converges weekly in $L_1^2$ to a $f_\infty\neq 0.$
\item[(2)] There is a sequence of numbers $\{\sig_j\}$ which goes to $0$ ($j\to \infty$)
such that the sequence $\{f_\infty^{\sig_j}\}$ converges weekly in $L_1^2$ to a $ f_\infty^0$, 
where $0<\sig_j \leq 1.$
\item[3)] $\pi:=\id_E- f_\infty^0$ is a weekly generalized holomorphic subbundle of $E.$
\end{enumerate}
\end{proposition}

\bgn{proof}
(1) 
Since every eigenvalue of $f_i$ is less than or equal to $1$, the $L^2$-norm of $f_i$ is bounded by 
the volume of $M$. Thus it suffices to show that the $L^2$-norm of $d_0f_i$ is bounded from above by a constant $C$
which does not depend on $\e.$
From Proposition \ref{prop:Let f in Herm}
 (2), we obtain 
 $$
\frac1\sig \trian (\tr\, f^\sig)\vol_M+ \e h^{\top}_0((\log f)\psi, \,\, f^\sig \psi)+|f^{\frac{-\sig}2}(\pa_0 f^{\sig} \psi)|^2_{h_0}\vol_M
\leq -h^{\top}_0(\h \K_0(\psi), \,\, f^\sig \psi)
$$

Let $\sig=1$. Then we have
$$
 \trian (\tr\, f_\e)\vol_M
\leq -h^{\top}_0(\h \K_0(\psi), \,\, f_\e \psi)- \e h^{\top}_0((\log f_\e)\psi, \,\, f_\e \psi)
$$
Since $\e\|\log f_\e\|<C$ , $h_0(\h \K_0(\psi), \,\, f_\e \psi)<C|f_\e|\,\,\vol_M$, we have 
$$\trian (\tr\, f_\e)\vol_M< C |f_\e|\,\,\vol_M<C\tr\, f_\e\,\,\vol_M
$$
This means 
$\max(f_\e)\leq \|f_\e\|_{L^1}.$
Since $\rho_\e f_\e$ has at least one eigenvalue $1,$ 
 we have 
$$
1\leq\max_M (\rho_\e|f_\e|)<C \rho_\e \|f_\e\|_{L^1}<C\Vol(M)^{\frac12}\rho_\e\|f_\e\|_{L^2}<C\|\rho_\e f_\e\|_{L^2},
$$
since $\rho_\e$ is a constant and we denotes by $C$ nonzero constants which does not depend on $\e$.
This implies that 
$$
0<\frac 1C\leq \|\rho_\e f_\e\|_{L^2}
$$
If $f_i:=\rho_{\e_i}f_{\e_i}$ converges strongly $f_\infty$ in $L^2$, then $f_\infty\neq 0.$
Further 
we see that 
$$
\|\rho_\e f_\e\|_{L^2}\leq \|\id_E\|_{L^2}=C
$$
Next we shall estimate $L^2_1$-norm of $\rho_\e f_\e.$
Since $f_\e=f_\e^*$ and $|\pa_0(\rho_\e f_\e)|=|\ol\pa(\rho_\e f_\e)|$, we have 
$$
|d_0(\rho_\e f_\e)|^2=2|\pa_0(\rho_\e f_\e)|^2\leq 2
|(\rho_\e f_\e)^{\frac{-1}2}\pa_0(\rho_\e f_\e)|^2,
$$
where $|(\rho_\e f_\e)^{\frac{-1}2}|\geq |\id_E|.$

By multiplying $\rho_{\e_i}^\sig$ on the both sides of Proposition \ref{prop:Let f in Herm}
 (2), we obtain 
$$
\frac1\sig \trian (\tr\, f_i^\sig)\vol_M+\e h^{\top}_0((\log f_{\e_i})\psi, \,\, f_i^\sig \psi)+|f_i^{\frac{-\sig}2}(\pa_0 f_i^{\sig} \psi)|^2_{h_0}\vol_M
\leq -h^{\top}_0(\h \K_0(\psi), \,\, f_i^\sig \psi),
$$
Note that $\rho_{\e_i}$ is a constant and $f_i:=\rho_{\e_i}f_{\e_i}$.
Let $\sig=1$. Applying $\int_M\trian (\tr\, f^\sig)\vol_M=0$, we obtain
\bgn{align}
\int_M |f_i^{\frac{-1}2}(\pa_0 f_i \psi)|^2_{h^{}_0}\vol_M&\leq -\int_M\e h^{\top}_0((\log f_{\e_i})\psi, \,\, f_i \psi)
-\int_Mh^{\top}_0(\h \K_0(\psi), \,\, f_i \psi)\\
&\leq C\|f_i\|_{L^1}<C,
\end{align}
because $\e|\log f_\e|<C$ and $\max |\h \K_0|<C.$
Thus we have 
$$
\int_M |d_0(\rho_{\e_i} f_{\e_i})|^2< C
$$
Proof of (2) 
We also apply Proposition \ref{prop:Let f in Herm}
 (2) to $f_i^\sig.$
 Then we obtain 
$$
\int_M |f_i^{\frac{-\sig}2}(\pa_0 f_i^{\sig} \psi)|^2_{h_0}\vol_M
<C
$$
where $C$ does not depend on $i$ and $\sig.$
Since $|( f_i^\sig)^{\frac{-1}2}|\geq |\id_E|,$
$$
|d_0( f_i^\sig)|^2=2|\pa_0( f_i^\sig)|^2\leq 2
|( f_i^\sig)^{\frac{-1}2}\pa_0( f_i^\sig)|^2,
$$
Thus we have 
$$
\|f_i^\sig\|_{L^2_1}<C.
$$
Then we have $\|f_\infty^\sig\|_{L^2_1}<C.$
$f_\infty^{\sig_j}$ converges to $f_\infty^0$ weak in $L^2_1$ when we take a subsequence $\{\sig_j\}$

Proof of (3) 
Since $f_i^{\sig_j}$ is bounded, $f_i^{\sig_j}$ converges strongly to $ f_\infty^0$ in $L^2$(we take a subsequence if necessary). 
Since $M$ is compact, $f_i^{\sig}$ converges strongly to $ f_\infty^0$ in $L^1$ also. 
Thus there is a subset $W$ of $M$ such that 
the measure of $M\bsh W$ is zero and $f_i^{\sig_j}$ converges $ f_\infty^{\sig_j}$ on $W$. 
Note that $f_\infty$ and $f_\infty^0$ are defined on $W$. 
Since $(f_i^{\sig_j})^*=f_i^{\sig_j}$, this means that $(f_\infty^0)^*=f_\infty^0$.
Since $\pi:=\id_E- f_\infty^0,$ we have $\pi^*=\pi.$
Since $f_\infty\in \Herm(E, h_0),$ we can diagonalize $f_\infty.$
Then we see 
$$
\lim_{j\to\infty}f_\infty^{\sig_j}=\lim_{j\to\infty}f_\infty^{2\sig_j}=f_\infty^0
$$
Thus we have in $L^1$
$$
\pi^2=\lim_{j\to\infty}(\id_E- f_\infty^{\sig_j})^2=
\lim_{j\to\infty}\(\id_E-2 f_\infty^{\sig_j}+ f_\infty^{2\sig_j}\)
=\id_E-2f_\infty^0+f_\infty^0 =\pi
$$
We shall show that $(\id_E-\pi)\circ\ol\pa_\J^E\pi=0.$
Since $(\id_E-\pi)\circ \pi=0,$ we have 
$$\ol\pa_\J^E((\id_E-\pi)\circ \pi)=\(\ol\pa_\J^E(\id_E-\pi)\)\circ \pi +(\id_E-\pi)\circ\ol\pa_\J^E\pi=0
$$
Thus we have 
$$
|(\ol\pa_\J^E(\id_E-\pi))\circ \pi |=|(\id_E-\pi)\circ\ol\pa_\J^E\pi|.
$$
Since $\pi^*=\pi,$ we have  
$$
|(\ol\pa_\J^E(\id_E-\pi))\circ \pi |=\Big|\((\ol\pa_\J^E(\id_E-\pi))\circ \pi\)^* \Big|
=|\pi\circ \pa_0(\id_E-\pi)|
$$
Thus it suffices to show that $\|\pi\circ \pa_0(\id_E-\pi)\|_{L^2}=0.$
Note that $L^1$-norm is bounded from above by $L^2$-norm. 

\bgn{lemma}\label{lem:futoushiki1}
For $0\leq \lam \leq 1$ and $0<s\leq \sig \leq 1,$
we have the following inequality:
$$
0\leq \frac{s+\sig}s(1-\lam^s)\leq \lam^{-\sig}
$$
\end{lemma}
\bgn{proof}
This follows from calculation (c.f. \cite{Lu-Tel_1995}).
\end{proof}
Since $f_i$ can be diagonalize such that every eigenvalue is positive and less than or equal to $1,$
applying Lemma \ref{lem:futoushiki1}
for $0<s\leq \frac\sig2\leq 1,$ we have 
\bgn{equation}\label{eq:fracs}
0\leq (\id_E-f_i^s)\leq \frac s{s+\frac{\sig}2}\,\,f_i^{\frac{-\sig}2}
\end{equation}
Applying (\ref{eq:fracs}) and Proposition \ref{prop:Let f in Herm}, we obtain 
\bgn{align}
\int_M |(\id_E-f_i^s)\circ \pa_0f_i^\sig|^2 &\leq \frac{s}{s+\frac\sig2}\int_M |f_i^{\frac{-\sig}2}\circ \pa_0f_i^\sig|^2\\
&\leq \frac{s}{s+\frac\sig2} C,
\end{align}
where $C$ does not depend on $s,\sig$ and $i$.
When $i$ goes to $\infty$, we have 
$$
\int_M |(\id_E-f_\infty^s)\circ \pa_0f_\infty^\sig|^2\leq \frac{s}{s+\frac\sig2} C.
$$
When $s$ goes to $0$ , we have 
$$
\int_M |\pi \circ\pa_0f_\infty^\sig|^2\leq \lim_{s\to 0}\frac{s}{s+\frac\sig2} C=0
$$
Then  when $\sig$ goes to $0$, we obtain 
$$
\int_M |\pi \circ\pa_0f_\infty^0|^2=\int_M|\pi \circ\pa_0(\id_E-\pi)|^2=0
$$
Thus we obtain the result.
\end{proof}
The $\pi$ is a weak generalized holomorphic subbundle of $E.$
We also denote it by $E_\pi.$
Then $\pi$ gives rise to a coherent subsheaf ${\mathcal E_\pi}$ with respect to both $I_\pm$ of the bihermitian structure 
\bgn{proposition}
The rank of $E_\pi$ satisfies 
$$
0<\rk E_\pi < r=\rk E.
$$
\end{proposition}
\bgn{proof}
Since $f_i$ converges to $f_\infty$ on almost everywhere, we see that $f_\infty^0\neq 0$ and $f_\infty^0\neq \id_E$ on almost everywhere.
Thus $\pi:= \id_E-f_\infty^0\neq 0$ and $\pi\neq\id_E$ on almost everywhere. 
Thus we have the result.
\end{proof}

\bgn{proposition}
Let $\mu(E)$ be the slope of $E$ and $\mu(E_\pi)$ the slope of $E_\pi.$
Then we have 
$$
\mu(E_\pi)\geq \mu(E).
$$
\end{proposition}
\bgn{proof}
As in before, 
$\pi$ gives a holomorphic subbundle $E_\pi|_W$ on an open dense subset $W$ of $M$.
Applying Proposition \ref{prop:trpiScircKE} to $\pi$ on $W$ and $S:=E_\pi|_W$, we have 
$$\tr\(\pi\circ\K^E_\A(\psi)\circ \pi\) =\tr\,\K^S_{\A_S}(\psi)+\|H^S\|^2\psi,
$$
where $\A$ is the canonical connection with respect to $h_0$ and 
$\K_\A^E$ is the mean curvature of $E$ and 
$\K^S_{\A_S}(\psi)$ denotes the mean curvature of $S:=E_\pi|_W$ and $H^S$ is the second fundamental form
of $S$ of $E|_W$ with respect to $h_0.$
From Proposition \ref{leths} (2), we obtain 
$$\tr\(\pi\circ\K^E_\A(\psi)\circ \pi\) =\tr\,\K^S_{\A_S}(\psi)+\|\pa_0\pi\|^2\psi,
$$
Recall that the degree of $E$ is defined by
$$\deg (E):=\frac1{\Vol_M}\int_Mi^n\frac{1}{2\pi}\tr\,\lan\K^E_\A(\psi), \,\, \ol\psi\ran_s
$$
Thus the degree of $E_\pi$ is given by
\bgn{align}
\Vol_M\deg(E_\pi) =\frac 1{2\pi}\int_M i^n\tr\,\lan\K^S_{\A_S}(\psi), \,\,\ol\psi\ran_s=
\frac 1{2\pi}\int_M i^n\tr\,\lan\pi\circ\K^E_\A(\psi)\circ\pi, \,\,\ol\psi\ran_s-\frac1{2\pi}\int_M\|\pa_0\pi\|^2\vol_M
\end{align}
where $\vol_M =i^n\lan \psi, \,\,\ol\psi\ran_s.$
Let $r=\rk E$ and $r':=\rk E_\pi.$
Then the slope is given by 
$$
\mu(E)=\frac1r\deg(E).
$$
Let $\h\K^E_\A(\psi)=\K^E_\A(\psi)-\mu(E)\psi\id_E.$
Then we have 
$$
\mu(E_\pi) =\frac1{r'}\deg(E_\pi)=\frac 1{2\pi r'\Vol_M}\int_M i^n\tr\,\lan\pi\circ\h\K^E_\A(\psi)\circ\pi, \,\,\ol\psi\ran_s-\frac1{2\pi r'\Vol_M}
\int_M\|\pa_0\pi\|^2
+\mu(E)
$$
Thus the inequality $\mu(E_\pi)\geq \mu(E)$ is equivalent to the following inequality:
$$
\int_M i^n\tr\,\lan\h\K^E_\A(\psi)\circ\pi, \,\,\ol\psi\ran_s\geq
\int_M\|\pa_0\pi\|^2\vol_M
$$
since $\tr\,(\pi\circ\h\K^E_\A(\psi)\circ\pi)=\tr(\h\K^E_\A(\psi)\circ\pi).$
Recall that $\id_E-f_i^\sig$ strongly converges in $L^2,$
$$
\pi =\lim_{\sig\to 0}\lim_{i\to\infty}(\id_E-f_i^\sig).
$$
Since $\tr\,\h\K^E_\A(\psi)=0,$ we have 
$$
\int_M i^n\tr\,\lan\h\K^E_\A(\psi)\circ\pi, \,\,\ol\psi\ran_s=-\lim_{\sig\to 0}\lim_{i\to\infty}
\int_M i^n\tr\,\lan\h\K^E_\A(\psi)\circ f_i^\sig, \,\,\ol\psi\ran_s
$$
Since $f_{\e_i}$ satisfies the equation (\ref{eq:continuity method}) for some $\e_i>0$, then from Proposition \ref{prop:Let f in Herm} (2),
we have 
$$\label{frac1sigtrian}
\frac1\sig \trian (\tr\, f_i^\sig)\vol_M+ \e h^{\top}_0((\log f_{\e_i})\psi, \,\, f_i^\sig \psi)+|f_i^{\frac{-\sig}2}(\pa_0 f_i^{\sig} \psi)|^2_{h_0}\vol_M
\leq -h^{\top}_0(\h \K_0(\psi), \,\, f_i^\sig \psi)
$$
since $\h \K_0(\psi)=\h\K_{\A}(\psi)$. \bgn{lemma}\label{real numer sig}
 For real numbers $\sig\leq 0$ and $\lam_i$ with $\sum_i \lam_i=0,$
we have 
$$
\sum_i \lam_i e^{\sig \lam_i}\geq 0
$$
\end{lemma}
\bgn{proof}
This follows directly from calculation.
\end{proof}
We can assume that $f_{\e_i}$ can be diagonalize with eigenvalues $\{e^{\lam_i}\}.$
Since $\tr(\log f_{\e_i})=0$, applying Lemma \ref{real numer sig}, we have 
$$
\e h^{\top}_0((\log f_{\e_i})\psi, \,\, f_i^\sig \psi)\geq 0
$$
Substituting $f=f_i$ into Proposition \ref{prop:Let f in Herm} (2), we have 
$$
\int_M|f_i^{\frac{-\sig}2}(\pa_0 f_i^{\sig} \psi)|^2_{h_0}\vol_M\leq 
-\int_M h^{\top}_0(\h \K_0(\psi), \,\, f_i^\sig \psi)
$$
since $\int_M \trian f_i^\sig =0.$
Thus we have 
$$\int_M|f_i^{\frac{-\sig}2}(\pa_0 f_i^{\sig} \psi)|^2_{h_0}\vol_M
\leq
|\int_M i^n\tr\,\lan\h\K^E_\A(\psi)\circ f_i^\sig, \,\,\ol\psi\ran_s|
$$
Thus we obtain 
$$
\int_M|\pa_0 (f_i^{\sig} \psi)|^2_{h_0}\vol_M=\int_M|\pa_0 (\id_E-f_i^{\sig})\psi|^2_{h_0}\vol_M
\leq
|\int_M i^n\tr\,\lan\h\K^E_\A(\psi)\circ f_i^\sig, \,\,\ol\psi\ran_s|,
 $$
since $f_i^\sig\geq \id_E$ for $0\leq \sig <1.$
Since $\pa_0 (\id_E-f_i^{\sig})$ converges to $\pa_0\pi$ weak in $L^2,$ we obtain 
\bgn{align}
\|\pa_0\pi\|^2_{L^2}\leq (\pa_0(\pi), \,\pa_0(\pi))_{L^2}&=\lim_{\sig\to0}\lim_{i\to0}
(\pa_0(\pi), \, \pa_0(\id_E-f_i^\sig))_{L^2}\\
&\leq \lim_{\sig\to0}\lim_{i\to0}
\|\pa_0(\pi)\|_{L^2}\,\|\pa_0(\id_E-f_i^\sig)\|_{L^2}\\
&\leq \lim_{\sig\to0}\lim_{i\to0}
|\int_M i^n\tr\,\lan\h\K^E_\A(\psi)\circ f_i^\sig, \,\,\ol\psi\ran_s|\\
&\leq |\int_Mi^n\tr\,\lan\h\K^E_\A(\psi)\circ\pi, \,\,\ol\psi\ran_s|,
\end{align}
since $\tr\, \h\K^E_\A(\psi)=0.$
Thus we obtain 
$\|\pa_0\pi\|^2_{L^2}\leq |\int_Mi^n\tr\,\lan\h\K^E_\A(\psi)\circ\pi, \,\,\ol\psi\ran_s|.$
Hence we have 
$$
\mu(E_\pi)\geq \mu(E).
$$
\end{proof}
\subsection{Proof of main theorem}\label{Proof of main theorem}
\indent{\sc Theorem} \ref{kobayashi-Hitchin correspondence}. \text{\rm [Kobayashi-Hitchin correspondence]}
{\it There exists an Einstein-Hermitian metric on a $\psi$-polystable generalized holomorphic vector bundle.
Conversely, a generalized holomorphic vector bundle admitting an Einstein-Hermitian metric is $\psi$-polystable.
}
\bgn{proof} Let $E$ be a $\psi$-polystable generalized holomorphic vector bundle. 
Then we start the continuity method by solving the equation (\ref{eq:continuity method1}).
From Step 2, we already have the set of solutions $\{f_\e\}$ for $\e\in (0,1].$
If we assume that $\|f_\e\|_{C^0}$ goes to infinity when $\e\to 0,$
then from Step 3, we have a weak generalized holomorphic subbundle $\pi$ such that 
$\mu(\pi)\geq\mu(E).$
This is a contradiction since $E$ is $\psi$-stable.
Thus $\|f_\e\|_{C^0}$ is bounded and then there exits $f_0:=\lim_{\e\to 0}f_\e$ such that 
$h_{f_0}$ gives an Einstein-Hermitian metric. 
Conversely, a generalized holomorphic vector bundle with an Einstein-Hermitian metric is $\psi$-stable 
from Section \ref{From Einstein-Hermitian metrics to stability}.
\end{proof}


\section{Einstein-Hermitian metrics and stable Poisson modules}
\label{Einstein-Hermitian metrics and stable Poisson modules}
\subsection{Stable Poisson modules}\label{Stable Poisson module}
Let $X:=(M, J)$ be a complex manifold with a holomorphic Poisson structure $\b$. 
The Poisson bracket $\{\, ,\,\}$ is given by the contraction by 
the $2$-vector $\b$ and $df\w dg$ for 
a holomorphic function $f, g,$
$$
\{ f, g\}:=\b(df\w dg)
$$
Let $E$ be a holomorphic vector bundle over $X$.
We denote by ${\mathcal E}$ the sheaf ${\mathcal O}(E)$ of germs of holomorphic sections of $E$.
Then a Poisson module structure of ${\E}$ is a map from ${\O}_X\times\E$ to $\E$ given by
$s\mapsto\{ f, s\}$ for $f\in{\O}_X$ and $s\in{\E}$ which satisfies 
\bgn{align}
\{ f, gs\}=&\{f, g\}s+g\{f, s\}\\
\{ fg, s\}=&f\{g, s\}+g\{f, s\}\\
\{\{f, g\}, s\} =&\{f, \{ g, s\}\}-\{g, \{f,s\}\}\label{Poisson module}
\end{align}
for $f, g\in {\O}_X$ and $s\in\E.$
An ${\O}_X$-module $\E$ with a Poisson module structure is called a Poisson module.

We denote by $\Theta$ the sheaf of germs of holomorphic vector  fields of $X$.
Then $\b$ gives a differential operator $\del_\b: {\O}_X \to \Theta_X$ by 
$$
\del_\b f= \b(df), \qquad f \in {\mathcal O}_X,
$$
where $\b(df)\in\Theta$ is the contraction of $\b\in\w^2\Theta$ and $df\in \Ome^1$.
Let $\del^\E_\b$ be a differential operator whose principal part is given by 
$\del_\b$, that is,
$$
\xymatrix {
&\del_\b^{\E}: {\E}\ar@{->}[r]&{\E}\otimes\Theta,
}
$$
which satisfies 
\bgn{equation}\label{pricipal part}
\del^{\E}_\b(fs)=f\del^{\E}_\b(s)+ (\del_\b f)s, \qquad s\in {\E}, \, f\in {\mathcal O}_X
\end{equation}
Then a Poisson module structure $\{\, ,\,\}$ gives a differential operator $\del_\b^\E$ by 
$$
\{ f, s\} =\lan\del_\b^\E(s), \, df\ran,
$$
where $\lan\,,\,\ran$ denotes the coupling between $\Theta$ and $\Ome^1.$
Then 
$\del_\b^\E$ gives the the following sequence:
$$
\xymatrix{
0\ar@{->}[r]&{\E}\ar@{->}^{\del_\b^{\E}\,\,\,\,\,\,}[r]&{\E}\otimes\Theta
\ar@{->}^{\del_\b^{\E}\,\,\,\,\,\,}[r]&{\E}\otimes\w^2\Theta\ar@{->}^{\,\,\,\del_\b^{\E}}[r]&\cdots
}
$$
Then the condition (\ref{Poisson module}) is equivalent to
$$\del^{\E}_\b\circ\del^{\E}_\b=0.$$
Thus it turns out that a Poisson module structure is equivalent to a differential operator 
$\del_\b^\E: \E\to \E\otimes \Theta$ which satisfies (\ref{pricipal part}) and
$\del^{\E}_\b\circ\del^{\E}_\b=0.$
Let $(M, \J_\b)$ be a generalized complex manifold given by Poisson deformation of $\b$ as in Example \ref{Poisson deformations}.
Then it is known that the followings (1) and (2) are equivalent 
\bgn{enumerate}
\item[(1)] a Poisson module $({\E},\del^{\E}_\b)$ on $X$
\item[(2)] a generalized holomorphic vector bundle $(E, \ol\pa^E_{\J_\b})$ over $(M, \J_\b)$
\end{enumerate}
In fact, the exterior derivative $d$ is decomposed into 
$d=\pa_{\J_\b}+\ol\pa_{\J_\b}$, where 
$\ol\pa_{\J_\b}=e^{\ol\b}\circ \ol\pa\circ e^{-\ol\b}-
\del_\b,$
where $\ol\pa$ is the ordinary $\ol\pa$-operator and
$e^{\ol\b}\circ \ol\pa\circ e^{-\ol\b}$ denotes the Adjoint action of $e^{\ol\b}$ on $\ol\pa$.
Thus a generalized holomorphic structure 
$\ol\pa^E_{\J_\b}$ is uniquely decomposed into 
$$
\ol\pa^E_{\J_\b}=e^{\ol\b}\circ \ol\pa^E\circ e^{-\ol\b}
-\del^\E_\b,
$$
where $\ol\pa^E$ denotes the ordinary holomorphic structure on $E$. 
Hence a Poisson structure $\del^\E_\b$ on a holomorphic vector bundle $E$ gives a generalized holomorphic structure, and vice versa.
\par
Then the cohomology group of the Lie algebroid complex complex 
$(E\otimes\w^\bullet \ol L_{\J_\b}, \ol\pa^E_{\J_\b})$ is given by the hyper-cohomology of the Poisson complex $( {\E}\otimes\w^\bullet\Theta, \del^{\E}_{\b}).$
The degree of a Poisson module $(\E, \del_\b^\E)$ is the degree of $\E$ and the slope is defined to be the slope of $E$ also.
A sub Poisson module $({\mathcal F}, \del_\b^{\F})$ is  a subsheaf ${\F}$ of ${\E}$ satisfying 
$$
\xymatrix{&{\F}\ar@{->}^{\del_\b^{\F}\,\,\,\,}[r]\ar@{->}_i[d]&{\F}\otimes\Theta\ar@{->}[d]\\
&{\E}\ar@{->}^{\del_\b^{\E}\,\,\,\,}[r]&{\E}\otimes\Theta
}
$$
where $i: \F\to \E$ denotes the inclusion.
Thus we have the stability of Poisson modules which coincides with the notion of  the stability of generalized holomorphic vector bundles
\bgn{definition}
A Poisson module $(\E, \del_\b^\E)$ is {\it stable} if for every sub Poisson module $(\F, \del^\F_\b)$ of $(\E, \del^\E_\b)$, 
the slope strict inequality holds
$$
\mu(\F)<\mu(\E)
$$
We also have the notion of semistability and polystability:
If $\E$ satisfies $\mu(F)\leq\mu(E)$, then $\E$ is {\it semistable}.
If $\E$ is a direst sum of Poisson modules $\oplus_i (\E_i, \del_\b^{\E_i})$ with the same slope $\mu(E)$, then 
$\E$ is {\it polystable}.
\end{definition}
Then we have the following:
\bgn{theorem}
Let $(E, \ol\pa^E_{\J_\b})$ be a generalized holomorphic vector bundle which is given by a Poisson module $(\E, \del_\b^\E)$.
Then the following (1) and (2) are equivalent \\\noindent
(1) $(E, \ol\pa^E_{\J_\b})$ admits an Einstein-Hermitian metric \\\noindent
(2) $(\E, \del_\b^\E)$ is polystable.
\end{theorem}
\bgn{proof}
The result follows from our main theorem.
\end{proof}

\subsection{Einstein-Hermitian metrics of Poisson modules}
Let $X=(M, J)$ be a Poisson manifold with a Poisson structure $\b$ and 
$\J_\b$ the induced generalized complex structure. 
Then as in (\ref{eigenspace decomposition}), $\J_\b$ gives the eigenspace decomposition
$(\TT)^\C=\L_{\J_\b}\oplus\ol{\L}_{\J_\b}$ and 
$\L_{\J_\b}=\w^{1,0}_\b\oplus T^{0,1}_J$ and 
$\ol{\L}_{\J_\b}=\w^{1,0}_\b\oplus T^{1,0}_J$, where 
and $\w^{1,0}_\b$ and $\w^{0,1}_\b$ are respectively given by 
\bgn{align}
\w^{1,0}_\b=&\{\t +[\b, \t]\,|\, \t\in \w^{1,0}_J\, \} \\
\w^{0,1}_\b=&\{\ol\t+[\ol\b, \ol\t]\,|\,\ol\t\in \w^{0,1}_J\, \},
\end{align}
where $\w^{1,0}_J$ denotes forms of type $(1,0)$ with respect to the ordinary complex structure $J$ and $\w^{0,1}_J$ is the complex conjugate of $\w^{0,1}_J$ and 
$T^{1,0}_J$ denotes vectors of type $(1,0)$ and $T^{0,1}_J$ is the complex conjugate of $T^{1,0}_J.$
The generalized complex structure $\J_\b$ gives the decomposition $d=\pa_\b+ \ol\pa_\b$.
For a function $f$, $\pa_\b f$ and $\ol\pa_\b f$ are explicitly given by
\bgn{align}
\pa_\b f =&e^{\b}\pa f e^{-\b}-[\ol\b, \ol\pa f ]\\
=&\overset{\w^{1,0}_\b}{\overbrace{\pa f+[\b, \pa f]}}-
\overset{T^{0,1}}{\overbrace{[\ol\b, \ol\pa f ]}}\in \L_{\J_\b}\label{pabf=}\\
\ol\pa_\b f=&\overset{\w^{0,1}_\b}{\overbrace{\ol\pa f+[\b, \ol\pa f]}}-
\overset{T^{1,0}}{\overbrace{[\b, \pa f ]}}\in \ol\L_{\J_\b},
\end{align}
where 
$df=\pa_\b f +\ol\pa_\b f.$

Let $(E, h)$ be an Hermitian vector bundle and 
$\ol\pa^E=\ol\pa+A^{0,1}$ a (ordinary) holomorphic structure of $E.$
Then we denote by $d^A=\pa^E+\ol\pa^E$ the canonical connection 
of $(E, h)$ with a connection form $A.$
Let $\ol\pa^E_\b$ be a generalized holomorphic structure of $E$ with respect to the generalized complex structure 
$\J_\b.$
Then as in Subsection \ref{Stable Poisson module}, there is a holomorphic structure $\ol\pa^E$ and 
a Poisson module structure $\del^E_\b$ such that 
$$
\ol\pa^E_\b =e^{\ol\b}\circ \ol\pa^E\circ e^{-\ol\b}-\del^E_\b.
$$
Then $\ol\pa^E$ and $\del_\b^E$ are respectively written as 
\bgn{align}
\ol\pa^E=&\ol\pa+ A^{0,1}\\
\del^E_\b=&[\b, \pa^E]+V^{1,0},
\end{align}
where $\pa^E$ is the $(1,0)$-component of the ordinary canonical connection $d^A$ and $A^{0,1}\in \w^{0,1}_J\otimes\End(E)$ denotes the connection form of $\ol\pa^E$ and $V^{1,0}\in T^{1,0}_J\otimes\End(E).$
Note that $\del^E_\b$ is a holomorphic operator but $\pa^E$ is not holomorphic. 
Since $\ol\pa^E_\b\circ\ol\pa^E_\b=0$ and $[e^{\ol\b},\b]=0$, we have 
\bgn{align}
&\ol\pa^E\circ\ol\pa^E=0\\
&\ol\pa^E\circ\del_\b^E+\del^E_\b\circ\ol\pa^E=0\\
&\del^E_\b\circ\del^E_\b=0
\end{align}
For given generalized holomorphic structure $\ol\pa^E_\b,$
we have the canonical generalized connection $\D^E_\b=\D^{1,0}_\b+\ol\pa^E_\b$  of $(E, h),$ where $\D^{1,0}_\b\in \L_{\J_\b}\otimes\End(E).$
Then we have
\bgn{proposition}\label{the canonical conn of poisson module}
$$
\D^E_\b =d^A+V
$$
where $V=V^{0,1}+V^{1,0}$ and 
$-V^{1,0}$ is the adjoint $(V^{1,0})^{*_h}$of $V^{1,0}$ with respect to $h$.
\end{proposition}
\bgn{proof}
The generalized holomorphic structure $\ol\pa^E_\b$ is explicitly given by
\bgn{align}
\ol\pa_\b s =&\overset{\w^{0,1}_\b}{\overbrace{\ol\pa^Es+[\ol\b, \ol\pa^Es]}}-\overset{T^{1,0}}{\overbrace{[\b, \pa^E s ]+V^{1,0}s}}\in \ol\L_{\J_\b}\otimes\End(E)
\end{align}
We define the following operator $\pa'_\b$ by
\bgn{align}
\pa'_\b s =&\overset{\w^{1,0}_\b}{\overbrace{\pa^Es+[\b, \pa^E s]}}-\overset{T^{0,1}}{\overbrace{[\ol\b,\,\, \ol\pa s ]+V^{0,1}s}}\in \L_{\J_\b}\otimes\End(E)
\end{align}
Since $d^A=\pa^E+\ol\pa^E$ is the ordinary canonical connection of $(E, h),$
we have 
\bgn{align}
\pa h(s_1, s_2) =&h(\pa^E s_1, s_2) + h(s_1, \ol\pa^Es_2),\\
\ol\pa h(s_1, s_2) =&h(\ol\pa^E s_1, s_2) + h(s_1, \pa^Es_2),
\end{align}
for sections $s_1, s_2\in E.$
Since $h$ is Hermitian,  we also have
\bgn{align}
[\b\,\,,\pa h(s_1, s_2)]=&h([\b,\,\,\pa^E s_1], \,\,s_2)+h(s_1, \,\,[\ol\b,\,\,\ol\pa^Es_2])\\
[\ol\b,\,\, \ol\pa h(s_1, s_2)]=&h([\ol\b,\,\,\ol\pa s_1], \,\,s_2)+h(s_1, \,\,[\b,\,\,\pa s_2])
\end{align}
Since $V^{1,0}=-(V^{0,1})^{*_h}$, we have $h(V^{1,0}s_1, s_2)+ h(s_1, V^{1,0}s_2)=0.$
From (\ref{pabf=}), we obtain 
$$
\pa_\b h (s_1,s_2) = \pa h(s_1, s_2) +[\b, \pa h(s_1, s_2)]-[\ol\b, \ol\pa h(s_1, s_2)]
$$
Thus we have 
$$
\pa_\b h(s_1, s_2) =h (\pa'_\b s_1, s_2) + h (s_1, \ol\pa^E s_2).
$$
We also have 
$$
\ol\pa_\b h(s_1, s_2) =h (\ol\pa^E_\b s_1, s_2) + h (s_1, \pa'_\b s_2).
$$
Thus $\D_\b^E=\pa'_\b.$
Hence the canonical connection $\D^E=\pa'_\b+\ol\pa_\b^E$ is given by 
$$
\D^E_\b =d^A+V
$$
\end{proof}
\bgn{proposition}
Let $(E, h)$ be an Hermitian vector bundle over a generalized K\"ahler manifold 
$(M, \J_\b, \J_\psi)$.
We denote by $\D_\b^E$ the canonical connection of $(E, h)$ as in 
Proposition \ref{the canonical conn of poisson module}. 
Then the curvature and the mean curvature of the canonical connection $\D_\b^E$
are respectively given by 
\bgn{align}
&\F_{D_\b}(\psi) =F_A\cdot\psi +
d^A(V\cdot\psi)
+ \frac12[V\cdot V]\cdot\psi\\
&\K_{D_\b}=\pi_{U^{-n}}^{\Herm}\F_{\D_\b}(\psi)
\end{align}
\end{proposition}
\bgn{proof}
The result follows from Proposition \ref{the canonical conn of poisson module}.
\end{proof}

\subsection{D-modules and Poisson modules over complex surfaces}
{\it A log Poisson structure} is 
a holomorphic Poisson structure which is the dual of a log symplectic structure. 
A log Poisson structure admits remarkable features, that is, the Poisson cohomology is calculated by using the singular cohomology of the complement of 
zero divisor of a log Poisson structure coupled with datas of the Jacobi rings of singularities of the zero divisor.
In particular,  log Poisson structures on complex surfaces provide important classes of Poisson structures. 
\par
Let $X=(M, J)$ be a complex surface. 
Then every holomorphic section $\b$ of the anti-canonical line bundle is a holomorphic Poisson structure.
If the zero divisor $C$ of $\b$ is smooth, $\b$ is a log Poisson structure.
Then $C$ is an elliptic curve which is of particular interest.

\bgn{proposition}\label{poisson and meromorphic connection}
Let $X$ be a Poisson surface with a simple normal crossing anticanonical divisor $C$.
A Poisson module on a complex surface $X$ is given by a meromorphic flat connection which admits single pole along the anti-canonical divisor
$C$. Thus a Poisson module $({\E, \del^{\E}_\b})$ gives a $\mathcal{D}$-module over $X$
\end{proposition}
\bgn{proof}
A Poisson structure $\b$ gives a map from the sheaf of germs of holomorphic $1$-forms $\Ome^1_X$ to 
$\Theta$ by the contraction of $\b$ and $\t\in \Ome^1.$
Then the map is extended to a map $i_\b$ from $\h\Ome^1$ to $\Theta,$
where $\h\Ome^1$ is the sheaf of meromorphic $1$-forms with single pole along $C$. 
Then it turns out that the map is isomorphism 
$$i_\b:\h\Ome^1\cong\Theta.$$
The isomorphism $i_\b$ is extended to an isomorphism $\w^p i_\b\: \w^p\h\Ome^1\cong \w^p\Theta$ by 
$$\t_1\w\cdots\w\t_p\mapsto i_\b\t_1\w\cdots\w i_\b\t_p.$$
Let $\h\Ome^p=\w^p\h\Ome^1.$
Then it turns out  that we have the isomorphism between the complex $(\h\Ome^\bullet, d)$ and the Poisson complex 
$(\w^\bullet, \del_\b)$ since $\w^{p+1}i_\b\circ d=\pm \del_\b\circ \w^p i_\b.$
Tensoring $\E$ with the both complexes, we obtain the isomorphism between the complex $(\E\otimes\w^\bullet\h\Ome^1, d^\E)$
the Poisson module complex 
$(\E\otimes \w^\bullet\Theta, \del^\E_\b)$, that is, 
$$
\xymatrix{
0\ar@{->}[r]&{\E}\ar@{->}[d]\ar@{->}^{d^{\E}\,\,\,\,\,\,}[r]&{\E}\otimes\h\Ome^1\ar@{->}[d]
\ar@{->}^{d^{\E}\,\,\,\,\,\,}[r]&{\E}\otimes\h\Ome^2\ar@{->}[d]\ar@{->}^{\,\,\,d^{\E}}[r]&\cdots\\
0\ar@{->}[r]&{\E}\ar@{->}^{\del_\b^{\E}\,\,\,\,\,\,}[r]&{\E}\otimes\Theta
\ar@{->}^{\del_\b^{\E}\,\,\,\,\,\,}[r]&{\E}\otimes\w^2\Theta\ar@{->}^{\,\,\,\del_\b^{\E}}[r]&\cdots
}
$$
Since $\del^\E_\b\circ \del^\E_\b=0,$ we have $d^\E\circ d^\E=0.$
Thus $d^\E$ give a flat connection whose connection form allows a single pole along the anti-canonical divisor $C$.
Then by using a meromorphic flat connection $\E$ becomes a $\D$-modules.
\end{proof}
\subsection{The Serre construction of stable and unstable Poisson modules}
The book \cite{OSS_1980} is a good reference of the Serre construction in this subsection. 
Let $X$ be a del Pezzo surface which is a projective complex surface whose anticanonical line bundle $K_X^{-1}\cong{\O(N)}$ for a positive integer $N.$
We shall construct examples of stable Poisson modules over $X$.
Let $J_Y$ be the ideal sheaf which is given by a set of points $Y=\{p_1,\cdots, p_m\}\subset X$.
Then we have a holomorphic vector bundle $E$ as an extension: 
\bgn{equation}\label{extension}
\xymatrix{
0\ar@{->}[r]&{\O_X}\ar@{->}[r]&E\ar@{->}[r]&{\O}_X(k)\otimes J_Y\ar@{->}[r]&0,
}
\end{equation}
where $c_1(E)=k, c_2(E)=m.$
Let ${\O}_X(k)\otimes J_Y=J_Y(k)$.
Then an extension is given by an element $e$ of Ext$^1( J_Y(k), {\O}_X)$.
The spectral sequence  whose $E_2$-term is given by
$\{H^p(X, {\mathcal Ext}^q(J_Y(k), {\O_X}))\}$ degenerates to Ext$^\bullet(J_Y(k), {\O}_X)$. 
Since ${\mathcal Hom}(J_Y(k), {\O}_X)= {\mathcal Hom}( {\O}_X(k), {\O}_X)={\O}_X(-k)$,
we have the local and global exact sequence:
$$
0\to H^1(X,  {\O}_X(-k)))\to  \text{\rm Ext}^1( J_Y(k),\,{\O}_X)
\to H^0(X, {\mathcal Ext}^1( J_Y(k),\,{\O}_X))\overset{d_2}\to H^2(X,{\O}(-k))\to \cdots
$$
Since ${\mathcal Ext}^1(J_Y(k),\,{\O}_X, )$ is the Skyscraper sheaf with support $Y=\{p_1, \cdots, p_m\}$, we see 
$$
H^0(X, {\mathcal Ext}^1(J_Y(k),\,{\O}_X ))=\oplus_{p_i\in Y} \C_{p_i}.
$$
Let $\{e_{p_i}\}\in \oplus_{p_i\in Y}\C_{p_i}$ be the image of $e\in \text{\rm Ext}^1( J_Y(k),\,{\O}_X).$
Then it is known that $e$ gives a holomorphic vector bundle if and only if each $e_{p_i}$ does not vanish \cite{GH_1978}.
From the Serre duality, we have  $H^2(X, {\O}_X(-k))=H^0(X, K_X(k))=H^0(X, {\O}(k-N)).$
If $k<N$, then $H^2(X,{\O}(-k))=0.$
Thus we have an extension $e$ such that $e_{p_i}\neq 0$ for all $p_i\in Y$. In the case of $k\geq N$,  $ \dim H^0(X, {\O}(k-N))\neq 0.$
If there is an element $\{ e_{p_i}\}\in \oplus_{p_i\in Y}\C_{p_i}$ such that $e_{p_i}\neq 0 $ for all $p_i\in Y$ and
$d_2(\{ e_{p_i}\})= 0,$ then
we obtain an extension as a vector bundle (c.f. Global duality in Chapter 5 of \cite{GH_1978}).

\bgn{proposition}\label{prop:Poisson module on surface}
Let $C$ be the anti-canonical divisor on $X$ which is the zero set of a Poisson structure $\b$. We assume that $C$ is a simple normal crossing divisor.
Let $E$ be a holomorphic vector bundle $E$ which is an extension in the cases of $k=N,$
\bgn{equation}\label{extension k=N}
\xymatrix{
0\ar@{->}[r]&{\O_X}\ar@{->}[r]&E\ar@{->}[r]&{\O}_X(N)\otimes J_Y\ar@{->}[r]&0,
}
\end{equation}
Then $E$ has a Poisson module structure if  
$Y$ is included in the anti-canonical divisor $C$.
\end{proposition}

\bgn{remark}
Hitchin already constructed Poisson modules in the cases of a smooth elliptic curve  by 
using Dolbeault type formulation of the Serre construction \cite{Hi_2011-2}.
\end{remark}
In order to show Proposition \ref{prop:Poisson module on surface},
we need local results of meromorphic extension of flat holomorphic connections on $\C^2$. Let $I$ be the ideal sheaf with support at the origin which is generated by $z_1, z_2$ over $\O_{\C^2}$, where $(z_1,z_2)$ be coordinates of $\C^2.$
Then every nontrivial extension of $I$ by $\O_{\C^2}$ is given by 
\bgn{equation}\label{extension k=N}
\xymatrix{
0\ar@{->}[r]&{\O_{\C^2}}\ar@{->}^i[r]&\O_{\C^2}\oplus \O_{\C^2}\ar@{->}[r]^j& I\ar@{->}[r]&0,
}
\end{equation}
where $i(1)=(\lam z_1, \lam z_2)$ and $j(f, g)=-z_2f+ z_1g$, for $\lam \neq 0\in \C.$
For simplicity, we consider the case of $\lam =1.$
We denote by $C$ the curve which is given by $\{z_1=0\}.$
Let $\nab^I$ be a logarithmic connection of the complement $I|_{\C^2\bsh C}$
$$
d-\frac{dz_1}{z_1}
$$
Then we have $\nab^I z_1=0.$
We denote by $\nab^{(0)}$ the trivial connection of $\O_{\C^2}$.
By using the splitting map $s$ such that $j\circ s=\id$ on the complement $\C^2\bsh C,$
we have an isomorphism $i\oplus s: \O_{\C^2\bsh C}\oplus I_{\C^2\bsh C}\to 
\O_{\C^2\bsh C}\oplus\O_{\C^2\bsh C}$. Then $\nab^{(0)}\oplus \nab^I$ induces a holomorphic connection $\nab$ on 
$\O_{\C^2\bsh C}\oplus\O_{\C^2\bsh C}$ by using the isomorphism $i\oplus s.$
Then we have 
\bgn{lemma}\label{extension lemma}
The holomorphic connection $\nab$ can be extended to be a meromorphic connection of 
$\O_{\C^2}\oplus\O_{\C^2}$ with a simple pole along $C=\{z_1=0\}.$
\end{lemma}
\bgn{proof}
Let $\{e_1, e_2\}$ be a basis $\O_{\C^2}\oplus\O_{\C^2}$ which is given by 
$e_1=(1,0)$ and $e_2=(0,1).$
Since $i\circ \nab^{(0)}=\nab\circ i$, 
 it follows from $\nab^{(0)}1=0$ that 
\bgn{align}\label{nab z1e1+z2e2)}
\nab(z_1 e_1+ z_2 e_2)=&(dz_1) e_1 +(dz_2) e_2 +z_1\nab e_1+z_2\nab e_2=0.
\end{align}
A splitting map $s$ is written as 
$s(z_1)=e_2+ f(z_1 e_1+ z_2 e_2)$ for $f\in \O_{\C^2}.$
since $\nab s(z_1)=\nab^I z_1=0$, it follows from (\ref{nab z1e1+z2e2)}) that 
\bgn{equation}\label{nab sz1=}
\nab s(z_1)=\nab e_2 + df (z_1 e_1+ z_2 e_2)=0.
\end{equation}
Then from (\ref{nab z1e1+z2e2)}) and
(\ref{nab sz1=}), we obtain 
\bgn{align}
\nab e_1=&-\frac{dz_1}{z_1}e_1-\frac{dz_2}{z_1}e_2+df(z_2 e_1+\frac{z_2^2}{z_1}e_2)\\
\nab e_2=&-df(z_1e_1+z_2e_2).
\end{align}
Thus the connection $\nab $ is extended to be a meromorphic connection on $\C^2$ with simple pole along $\{z_1=0\}$
\end{proof}
\bgn{lemma}\label{extension lemma of simple normal crossing}
Let $C'$ be a curve defined by $\{z_1z_2=0\}$ on $\C^2$ and $I$ the ideal sheaf generated by $z_1, z_2$ as before. 
We define a connection $\nab^{I'}$ by 
$$d-\frac{dz_1}{z_1}-\frac{dz_2}{z_2}.$$
Then by using a splitting map $s$ as in (\ref{extension k=N}), 
the connection $\nab^{(0)}\oplus \nab^{I'}$  yields a connection $\nab' $ on $\O_{C^2\bsh C'}\oplus \O_{\C^2\bsh C'}$. 
Then $\nab$ is extended to be a meromorphic connection on $\C^2$ with a simple pole along $C'=\{z_1z_2=0\}$
\end{lemma}
\bgn{proof}Our proof is the similar as before. 
For a splitting map $s$, we have $s(z_1 z_2) =z_2 e_2 + f(z_1e_1+z_2 e_2).$
Then $\nab' s(z_1z_2) =\nab^{I'}z_1z_2=0.$. Thus we have 
$$
z_2\nab' e_2+ dz_2 e_1 + (df)(z_1e_1+z_2e_2)=0
$$
Since $\nab' (z_1 e_1+z_2e_2)=0,$ we also have 
\bgn{align}
\nab' e_1=&-\frac{dz_1}{z_1}e_1-\frac{dz_2}{z_2}e_2-\frac{z_2}{z_1}\nab' e_2\\
=&-\frac{dz_1}{z_1}e_1-\frac{dz_2}{z_2}e_2-\frac{dz_2}{z_1}e_1-(df)(e_1+\frac{z_2}{z_1}e_2)
\end{align}
 Hence $\nab'$ is extended to be a meromorphic connection with a simple pole along $C'=\{z_1z_2=0\}.$
\end{proof}

\bgn{proof} [ Proof of Proposition \ref{prop:Poisson module on surface}]
Let $X=\cup_\a U_\a$ be an covering of affine open sets.
The zero set of $\b$ is given by 
$\{f_\a=0\}$ on $U_\a$, where $\{f_\a\}$ gives a section of the anti-canonical line bundle.
Then it turns out that 
$$\{d-\frac{df_\a}{f_\a}\}$$
yields a meromorphic flat connection of $K_X^{-1}\cong \O(N),$ which has a logarithmic pole along $C.$ Thus it follows from Proposition \ref{poisson and meromorphic connection} that 
$K_X^{-1}$ is a Poisson module. 
Let $A_\a:=\dstyle-\frac{df_\a}{f_\a}$ be a connection form of the meromorphic connection of $K^{-1}_X$. 
Since the complement $X\bs C$ is Stein,
the extension (\ref{extension}) restricted to the complement $X\bsh C$ gives the trivial extension :

\bgn{equation}\label{extension restrict }
\xymatrix{
0\ar@{->}[r]&{\O_X}|_{\ss X\bsh C}\ar@{->}[r]&E|_{X\bsh C}\ar@{->}[r]&{\O}_X(N)|_{X\bsh C}\ar@{->}[r]&0,
}
\end{equation}
Thus $ E|_{X\bsh C}\cong {\mathcal O}_X|_{X\bsh C}\oplus {\mathcal O}_X(N)|_{X\bsh C}$ has the holomorphic flat connection $\nab$ which is 
induced from the product of flat connections of ${\mathcal O}_X|_{X\bsh C}\oplus {\mathcal O}_X(N)|_{X\bsh C}.$
From Lemma \ref{extension lemma}
and Lemma \ref{extension lemma of simple normal crossing} , it turns out that every holomorphic flat connection on $X\bsh C$ can be extended as a meromorphic connection along a single pole along $C.$
From Proposition \ref{poisson and meromorphic connection}, $E$ becomes a Poisson module.

\end{proof}
Let 
$\psi=e^{-\sqrt{-1}\ome}$, where  $\ome$ is a K\"ahler structure.
Then we have the notion of the ordinary stability with respect to $\ome$ and the notion of
$\psi$-stability of generalized holomorphic vector bundles.
As in before, a Poisson module $(\E, \del_\b^\E)$ gives 
the generalized holomorphic vector bundle $(E, \ol\pa^E_{\J_\b})$, where 
$\E={\O }(E)$.

\bgn{proposition}
If $E$ as in (\ref{extension}) is stable as a holomorphic vector bundle, then 
$E$ is $\psi$-stable as Poisson module (a generalized holomorphic vector bundle).
\end{proposition}
\bgn{proof}
Since we have the slope inequality $\mu(\F)<\mu(\E)$ for all subsheaves $\F$, we also have the slope inequality for all 
Poisson subsheaves.
\end{proof}

We shall explain the $\psi$-stability and $\psi$-semistability in the case $X=\C P^2.$
For $k=3,$ we have the following:
\bgn{equation}\label{In the cases of k=2}
\xymatrix{
0\ar@{->}[r]&{\O_X}\ar@{->}[r]&E\ar@{->}[r]&{\O}_X(3)\otimes J_Y\ar@{->}[r]&0,
}
\end{equation}
Then it is known that $E$ is stable if and only if $Y=\{p_1, p_2, p_3\}$ is not contained in a line $l\subset \C P^2$. 
If $Y\subset l$, then $E$ is not stable as a holomorphic vector bundle,
that is, $E$ does not admit any Einstein-Hermitian metric in the ordinary sense.

Thus it turns out that the ordinary stability and semistability of $\E$ depend on a configuration of points $p_1, \cdots, p_m$. 
However the following shows that $\psi$-stability of Poisson modules is different from the ordinary stability:

\bgn{proposition}\label{psi-stable always}
We assume that $C$ is smooth.
A Poisson module $E$ as in (\ref{In the cases of k=2})
is always $\psi$-stable. Thus $E$ admits an Einstein-Hermitian metric of a generalized holomorphic vector bundle over generalized K\"ahler manifold $(M, \J_{J_\b}, \J_\psi).$
\end{proposition}


The following Figure 1 in the cases of three points, explains the difference between the ordinary stability and $\psi$-stability:


\begin{center}

\tikzset{every picture/.style={line width=0.75pt}} 

\begin{tikzpicture}[x=0.75pt,y=0.75pt,yscale=-1,xscale=1]

\draw  [color={rgb, 255:red, 18; green, 115; blue, 229 }  ,draw opacity=1 ][fill={rgb, 255:red, 255; green, 255; blue, 255 }  ,fill opacity=1 ] (138,164.42) .. controls (138,109.74) and (220.85,65.42) .. (323.04,65.42) .. controls (425.24,65.42) and (508.08,109.74) .. (508.08,164.42) .. controls (508.08,219.09) and (425.24,263.42) .. (323.04,263.42) .. controls (220.85,263.42) and (138,219.09) .. (138,164.42) -- cycle ;
\draw    (174.08,170.83) .. controls (377.08,-11.75) and (232.08,335.25) .. (451.08,149.83) ;

\draw [color={rgb, 255:red, 208; green, 2; blue, 27 }  ,draw opacity=1 ]   (179,118) -- (454,210.75) ;

\draw   (220.08,132.75) .. controls (220.08,131.65) and (220.74,130.75) .. (221.54,130.75) .. controls (222.35,130.75) and (223,131.65) .. (223,132.75) .. controls (223,133.85) and (222.35,134.75) .. (221.54,134.75) .. controls (220.74,134.75) and (220.08,133.85) .. (220.08,132.75) -- cycle ; \draw   (220.9,132.07) .. controls (220.9,131.96) and (220.97,131.87) .. (221.05,131.87) .. controls (221.13,131.87) and (221.19,131.96) .. (221.19,132.07) .. controls (221.19,132.18) and (221.13,132.27) .. (221.05,132.27) .. controls (220.97,132.27) and (220.9,132.18) .. (220.9,132.07) -- cycle ; \draw   (221.89,132.07) .. controls (221.89,131.96) and (221.96,131.87) .. (222.04,131.87) .. controls (222.12,131.87) and (222.18,131.96) .. (222.18,132.07) .. controls (222.18,132.18) and (222.12,132.27) .. (222.04,132.27) .. controls (221.96,132.27) and (221.89,132.18) .. (221.89,132.07) -- cycle ; \draw   (220.81,133.95) .. controls (221.3,133.42) and (221.78,133.42) .. (222.27,133.95) ;
\draw   (308,159.75) .. controls (308,160.3) and (307.78,160.75) .. (307.5,160.75) .. controls (307.22,160.75) and (307,160.3) .. (307,159.75) .. controls (307,159.2) and (307.22,158.75) .. (307.5,158.75) .. controls (307.78,158.75) and (308,159.2) .. (308,159.75) -- cycle ; \draw   (307.72,160.09) .. controls (307.72,160.15) and (307.7,160.19) .. (307.67,160.19) .. controls (307.64,160.19) and (307.62,160.15) .. (307.62,160.09) .. controls (307.62,160.03) and (307.64,159.99) .. (307.67,159.99) .. controls (307.7,159.99) and (307.72,160.03) .. (307.72,160.09) -- cycle ; \draw   (307.38,160.09) .. controls (307.38,160.15) and (307.36,160.19) .. (307.33,160.19) .. controls (307.3,160.19) and (307.28,160.15) .. (307.28,160.09) .. controls (307.28,160.03) and (307.3,159.99) .. (307.33,159.99) .. controls (307.36,159.99) and (307.38,160.03) .. (307.38,160.09) -- cycle ; \draw   (307.75,159.15) .. controls (307.58,159.42) and (307.42,159.42) .. (307.25,159.15) ;
\draw   (394.94,190.24) .. controls (394.97,189.42) and (395.22,188.79) .. (395.5,188.84) .. controls (395.77,188.88) and (395.96,189.59) .. (395.93,190.41) .. controls (395.9,191.23) and (395.64,191.86) .. (395.37,191.81) .. controls (395.1,191.77) and (394.9,191.06) .. (394.94,190.24) -- cycle ; \draw   (395.24,189.78) .. controls (395.24,189.7) and (395.26,189.64) .. (395.29,189.64) .. controls (395.32,189.65) and (395.34,189.72) .. (395.34,189.8) .. controls (395.33,189.88) and (395.31,189.94) .. (395.28,189.94) .. controls (395.25,189.93) and (395.23,189.86) .. (395.24,189.78) -- cycle ; \draw   (395.57,189.84) .. controls (395.58,189.76) and (395.6,189.69) .. (395.63,189.7) .. controls (395.66,189.7) and (395.68,189.77) .. (395.67,189.86) .. controls (395.67,189.94) and (395.64,190) .. (395.62,190) .. controls (395.59,189.99) and (395.57,189.92) .. (395.57,189.84) -- cycle ; \draw   (395.15,191.18) .. controls (395.33,190.81) and (395.5,190.84) .. (395.64,191.26) ;

\draw (223,108) node   {$p_{1}$};
\draw (316,139) node   {$p_{2}$};
\draw (392,173) node   {$p_{3}$};
\draw (236,195) node   {$ \begin{array}{l}
C:elliptic\ curve\\
\end{array}$};
\draw (142,63) node   {$Poisson\ module\ over\ CP^{2}$};
\draw (93,100) node   {$Y=\{p_{1} ,p_{2} ,p_{3}\}$};
\draw (91,227) node  [align=left] {E is not poly stable \\in the ordinary sense\\However};
\draw (95,275) node   {$\psi -stable$};
\draw (418,217) node [color={rgb, 255:red, 208; green, 2; blue, 27 }  ,opacity=1 ]  {$line$};
\end{tikzpicture}
\hspace{8cm} Figure 1.
\end{center}

\bgn{proof}[Proof of Proposition \ref{psi-stable always}]

Tensoring ${\mathcal O}(-l)$ with the exact sequence (\ref{In the cases of k=2}), we obtain 
\bgn{equation}\label{}
\xymatrix{
0\ar@{->}[r]&{\O_X}(-l)\ar@{->}[r]&E(-l)\ar@{->}[r]&{J_Y}(3-l)\ar@{->}[r]&0.
},
\end{equation}
where 
$E(-l):=E\otimes{\mathcal O}(-l)$ and 
$J_Y(3-l):={\mathcal O}_X(3-l)\otimes J_Y.$
This sequence is an exact sequence of Poisson modules, that is, meromorphic flat connections are preserved.
Let $\nab^{E(-l)}$ be the meromorphic flat connection of $E(-l)$.
We denote by $H^0_\b(X, {E}(-l))$ the 
space of holomorphic parallel sections of ${ E}(-l),$
that is, 
$\{s\in H^0(X, {E}(-l))\, |\, \nab^{E(-l)}s=0\, \}.$
We also denote by $H^0_\b(X, {J_Y}(3-l))$ 
the space of holomorphic parallel sections of ${J_Y}(3-l).$
We have $H^1(X,{\mathcal O}(-l))=0$ since $X=\C P^2.$
If $l\geq 1,$ we have $H^0(X, {\mathcal O}_X(-l))=0.$

Then the long exact sequence gives an exact sequence of holomorphic parallel sections:
\bgn{equation}
\xymatrix{
0\ar@{->}[r]& 0\ar@{->}[r]&
H^0_\b(X, E(-l))\ar@{->}[r]&H^0_\b(X, {J_Y}(3-l))\ar@{->}[r]&0
},
\end{equation}
Thus we have the isomorphism $H^0_\b(X, E(-l))\cong H^0_\b(X, {J_Y}(3-l)).$

If $l\geq 3$, then $H^0(X, J_Y(3-l))=0.$ 
In the case of $l=2$, every section of $H^0(X,J_Y(1))$ has zero along a line. 
Since $C$ is smooth, every section of $H^0(X,J_Y(1))$ has zero on the complement $X\bsh C.$ Since the flat connection of $J_Y(3-l)$ is holomorphic on the complement 
$X\bsh C,$ it implies that $H^0_\b(X, J_Y(3-l))=0.$

We assume that there exists a destabilizing object, that is, 
a Poisson subsheaf $i: {\mathcal O}(l)\to \E$ with 
$l\geq 2.$
Thus the inclusion $i$  gives a nonzero parallel section $s$ of $H^0_\b(X, E(-l))$. However $H^0_\b(X, E(-l))=0$ and then $s=0.$
This is a contradiction.
\end{proof}
On the other hand, if the anti-canonical divisor admits singularities, then a different aspect of
$\psi$-stability appears.

\bgn{proposition}
Let $\b$ be a Poisson structure on $X:=\C P^2$ whose zero set consists of three lines $l_1, l_2$ and $l_3$ in a general position
(See Figure 2: three points and three lines).
We denote by $\E$ a 
Poisson module which is given by the extension as in (\ref{In the cases of k=2}).
If the support $Y$ of the ideal $J_Y$ is included in  one of lines, then $E$ is not $\psi$-stable.
\end{proposition}
\bgn{proof}
We assume that $Y$ is included in a line $l_1.$
Let $X=\cup_\a U_\a$ be an open covering as before and
$f_\a^{(1)}$ a defining equation of the line $l_1$ on each $U_\a.$
Then as before, 
$$
d-\frac{df_\a^{(1)}}{f_\a^{(1)}}
$$
yields a meromorphic flat connection $\nab^{(1)}$ of ${\O}(1)$ 
and the section $s^{(1)}:=\{f_\a^{(1)}\}$ is a flat section of ${\O}(1),$
with respect to $\nab^{(1)}.$
Since the zero set of $s^{(1)}$ is $l_1$ and $Y\subset l_1,$
we see that 
$H^0_\b(X, J_Y(1))\neq 0.$
Since $H^0_\b(X, J_Y(1))\neq 0$, 
there is an inclusion $i: {\O(2)}\to E$ which preserves flat connections. 
Thus there is a Poisson submodule ${\O}(2)$ of $E.$
Hence $E$ is not $\psi$-stable since $\mu(E)=\frac32.$
$E$ is irreducible, since $c_2(E)=m\neq0.$
Thus $E$ is not $\psi$-polystable.
\end{proof}

The following Figure 2 explains that 
a Poisson module is not $\psi$-stable if an anticanonical divisor consists of threee lines and
three points are on a line of them.

\bgn{center}\label{Figure: three lines}
\tikzset{every picture/.style={line width=0.75pt}} 

\begin{tikzpicture}[x=0.75pt,y=0.75pt,yscale=-1,xscale=1]

\draw  [color={rgb, 255:red, 18; green, 115; blue, 229 }  ,draw opacity=1 ][fill={rgb, 255:red, 255; green, 255; blue, 255 }  ,fill opacity=1 ] (112,156.42) .. controls (112,101.74) and (194.85,57.42) .. (297.04,57.42) .. controls (399.24,57.42) and (482.08,101.74) .. (482.08,156.42) .. controls (482.08,211.09) and (399.24,255.42) .. (297.04,255.42) .. controls (194.85,255.42) and (112,211.09) .. (112,156.42) -- cycle ;
\draw [color={rgb, 255:red, 208; green, 2; blue, 27 }  ,draw opacity=1 ]   (171.08,188.02) -- (440.54,190.79) ;

\draw   (198.08,187.75) .. controls (198.08,186.65) and (198.74,185.75) .. (199.54,185.75) .. controls (200.35,185.75) and (201,186.65) .. (201,187.75) .. controls (201,188.85) and (200.35,189.75) .. (199.54,189.75) .. controls (198.74,189.75) and (198.08,188.85) .. (198.08,187.75) -- cycle ; \draw   (198.9,187.07) .. controls (198.9,186.96) and (198.97,186.87) .. (199.05,186.87) .. controls (199.13,186.87) and (199.19,186.96) .. (199.19,187.07) .. controls (199.19,187.18) and (199.13,187.27) .. (199.05,187.27) .. controls (198.97,187.27) and (198.9,187.18) .. (198.9,187.07) -- cycle ; \draw   (199.89,187.07) .. controls (199.89,186.96) and (199.96,186.87) .. (200.04,186.87) .. controls (200.12,186.87) and (200.18,186.96) .. (200.18,187.07) .. controls (200.18,187.18) and (200.12,187.27) .. (200.04,187.27) .. controls (199.96,187.27) and (199.89,187.18) .. (199.89,187.07) -- cycle ; \draw   (198.81,188.95) .. controls (199.3,188.42) and (199.78,188.42) .. (200.27,188.95) ;
\draw    (305.08,63.02) -- (221.08,223.02) ;

\draw    (251.08,76.02) -- (381.08,218.02) ;

\draw   (287.08,189.75) .. controls (287.08,188.65) and (287.74,187.75) .. (288.54,187.75) .. controls (289.35,187.75) and (290,188.65) .. (290,189.75) .. controls (290,190.85) and (289.35,191.75) .. (288.54,191.75) .. controls (287.74,191.75) and (287.08,190.85) .. (287.08,189.75) -- cycle ; \draw   (287.9,189.07) .. controls (287.9,188.96) and (287.97,188.87) .. (288.05,188.87) .. controls (288.13,188.87) and (288.19,188.96) .. (288.19,189.07) .. controls (288.19,189.18) and (288.13,189.27) .. (288.05,189.27) .. controls (287.97,189.27) and (287.9,189.18) .. (287.9,189.07) -- cycle ; \draw   (288.89,189.07) .. controls (288.89,188.96) and (288.96,188.87) .. (289.04,188.87) .. controls (289.12,188.87) and (289.18,188.96) .. (289.18,189.07) .. controls (289.18,189.18) and (289.12,189.27) .. (289.04,189.27) .. controls (288.96,189.27) and (288.89,189.18) .. (288.89,189.07) -- cycle ; \draw   (287.81,190.95) .. controls (288.3,190.42) and (288.78,190.42) .. (289.27,190.95) ;
\draw   (379.08,189.75) .. controls (379.08,188.65) and (379.74,187.75) .. (380.54,187.75) .. controls (381.35,187.75) and (382,188.65) .. (382,189.75) .. controls (382,190.85) and (381.35,191.75) .. (380.54,191.75) .. controls (379.74,191.75) and (379.08,190.85) .. (379.08,189.75) -- cycle ; \draw   (379.9,189.07) .. controls (379.9,188.96) and (379.97,188.87) .. (380.05,188.87) .. controls (380.13,188.87) and (380.19,188.96) .. (380.19,189.07) .. controls (380.19,189.18) and (380.13,189.27) .. (380.05,189.27) .. controls (379.97,189.27) and (379.9,189.18) .. (379.9,189.07) -- cycle ; \draw   (380.89,189.07) .. controls (380.89,188.96) and (380.96,188.87) .. (381.04,188.87) .. controls (381.12,188.87) and (381.18,188.96) .. (381.18,189.07) .. controls (381.18,189.18) and (381.12,189.27) .. (381.04,189.27) .. controls (380.96,189.27) and (380.89,189.18) .. (380.89,189.07) -- cycle ; \draw   (379.81,190.95) .. controls (380.3,190.42) and (380.78,190.42) .. (381.27,190.95) ;

\draw (202,153) node   {$p_{1}$};
\draw (290,155) node   {$p_{2}$};
\draw (382,171) node   {$p_{3}$};
\draw (122,43) node   {$Poisson\ module\ over\ CP^{2}$};
\draw (73,80) node   {$Y=\{p_{1} ,p_{2} ,p_{3}\}$};
\draw (75,255) node   {$E\ is\ not\ \psi -stable$};
\draw (367,98) node   {$three\ lines$};

\end{tikzpicture}
\end{center}
\hspace{6cm}Figure 2: three points and three lines

\section{Vanishing theorems of generalized holomorphic vector bundles}
\label{Vanishing theorems of generalized holomorphic vector bundles}

\subsection{Vanishing theorem in the case of $\psi=e^{-\sqrt{-1}\ome}$}
Let $(E,h)$ be an Hermitian vector bundle over a generalized K\"ahler manifold 
of symplectic type $(M, \J, \J_\psi)$. 
In this subsection, we assume that  $\psi=e^{\sqrt{-1}\ome}$ for simplicity.
Let $\D: \Gam(E)\to \Gam(E\otimes(\TT))$ be a generalized connection:
$$
\D=d+A+V=D^A+V,
$$
where $D^A=d+A$ is an ordinary unitary connection and $V$ is a section of $\End(E)\otimes T_M.$
We have the decomposition of $(\TT)$ with respect to $\J$, that is,  
$$(\TT)^\C=\L_\J\oplus\ol\L_\J$$
and we also have the another decomposition with respect to $\J_\psi$
$$(\TT)^\C=\L_\psi\oplus\ol\L_\psi.$$
Since $\J$ and $\J_\psi$ commute, we have the simultaneous decomposition:
$$(\TT)^\C=\L_\J^+\oplus \L_\J^-\oplus \ol\L_\J^+\oplus\ol\L_\J^-,$$
where $\L_\J=\L_\J^+\oplus \L_\J^-$, 
$\L_\psi=\L_\J^+\oplus \ol\L_\J^-$. 
Note that $\L_\psi =\{\, e\in \TT\, |\, e\cdot \psi =0\}.$
The following table explains our decomposition:
\bgn{center}
\bgn{tabular}{|l|r|r|}\hline
&$\L_\J$&$\ol\L_\J$\\\hline
$\ol\L_\psi$&$\L_\J^-$&$\ol\L_\J^+$\\\hline
$\L_\psi$&$\L_\J^+$&$\ol\L_\J^-$\\\hline
\end{tabular}
\end{center}
According to the decomposition, 
a generalized connection  $\D$ is decomposed into
\bgn{align}
\D=&\D^{1,0}+\D^{0,1},\qquad \D^{1,0}\in \L_\J\otimes \End(E),\,\,\, \D^{0,1}\in \ol{\L}_\J\otimes \End(E)\\
\D=&\D'+\D'',\qquad\,\,\,\,\,\,\, \D'\in \L_\psi\otimes\End(E),\,\,\,\D''\in \ol{\L}_\psi\otimes\End(E)\\
\D=&\D^{1,0}_++\D^{1,0}_-+\D^{0,1}_++\D^{0,1}_-,
\end{align}
where $\D^{1,0}_\pm\in \L_\pm\otimes\End(E)$ and $\D^{0,1}_\pm\in \ol{\L}_\pm\otimes\End(E).$
Then we have

$$\D'=\D^{1,0}_++\D^{0,1}_-, \quad\D''=\D^{1,0}_-+\D^{0,1}_+$$

We shall introduce the following ordinary connections $D'$ and $D''$:
\bgn{align}
D'=&D^A+\sqrt{-1}\ome(V),\\
D''=&D^A-\sqrt{-1}\ome(V)
\end{align}
where $\ome(V)\in T_M^*\otimes u(E)$  is given by the contraction between the real $2$-form $\ome$
and $V\in $u$(E)\otimes T_M.$
Then we have
\bgn{lemma}\label{lem:dh(s1, s2)}
$$dh(s_1, s_2)=h(D's_1, s_2)+ h(s_1, D''s_2)$$
$$dh(s_1, s_2)=h(D''s_1, s_2)+ h(s_1, D's_2)$$
\end{lemma}
\bgn{proof}
Since $D^A$ is a unitary connection, we have 
$$
dh(s_1, s_2) =h(D^As_1, s_2) + h(s_1, D^A s_2)
$$
Since $V$ is a Skew-Hermitian matrix valued vector field,
then the contraction $\sqrt{-1}\ome(V)$  is an Hermitian matrix-valued $1$-form. 
Thus we have 
$$
0=h(\sqrt{-1}\ome(V)s_1, s_2) +h(s_1, \,-\sqrt{-1}\ome(V))
$$
Thus we have $dh(s_1, s_2)=h(D's_1, s_2)+ h(s_1, D''s_2).$
\end{proof}
As in before, 
$D'$, $D''$ are extended as the operators $d^{D'}, d^{D''}$ of
$E\otimes \w^*T^*_M$, respectively.
Note that $D'$ and $D''$ are connections which are different from generalized connections $\D'$ and $\D''.$
The relations between them are given by 
\bgn{lemma}\label{D' D''}
For a section $s\in E,$ we have
$$(\D''s)\cdot\psi=d^{D''}(s\, \psi)$$
$$(\D's)\cdot\ol\psi=d^{D'}(s\,\ol\psi)$$
\end{lemma}
\bgn{proof}
Since $\D=d+A+V=\D'+\D''$ and $\psi=e^{-\sqrt{-1}\ome}$, we have 
\bgn{align}
\D'=&\frac12(d-\sqrt{-1}\ome^{-1}d+ A-\sqrt{-1}\ome^{-1}(A)+V+\sqrt{-1}\ome(V))\\
\D''=&\frac12(d+\sqrt{-1}\ome^{-1}d+ A+\sqrt{-1}\ome^{-1}(A)+V-\sqrt{-1}\ome(V))
\end{align}
Then we have 
$$
(\D's)\cdot\ol\psi= ds\cdot\ol\psi+(As)\cdot\ol\psi+\sqrt{-1}(\ome(V)s)\cdot\ol\psi=d^{D'}(s\ol\psi)
$$
$$
(\D''s)\cdot\psi= ds\cdot\psi+(As)\cdot\psi-\sqrt{-1}(\ome(V)s)\cdot\psi=d^{D''}(s\psi)
$$
\end{proof}
For $s_1\otimes \a_1, s_2\otimes \a_2\in E\otimes \w^\bullet T^*_M$, we define 
\bgn{equation}\label{eq:hlan s_1otimes a_1,  s_2otimes a_2 ran}
h\lan s_1\otimes \a_1, \,\, s_2\otimes \a_2\ran:=h(s_1, s_2)\a_1\w\sig(\ol\a_2)\in \w^\bullet T^*_M.
\end{equation}
Taking the projection to the top forms, i.e., ($2n$-forms), we define
$$
h\lan s_1\otimes \a_1, \,\, s_2\otimes \a_2\ran_s:=h(s_1, s_2)\(\a_1\w\sig(\ol\a_2)\)_{top}\in \w^{2n}T^*_M.
$$
Then we obtain 
\bgn{lemma}\label{dhlans1otimes a1}
$$
d h\lan s_1\otimes \a_1, \,\, s_2\otimes \a_2\ran =h\lan d^{D'}(s_1\otimes \a_1),\,\, s_2\otimes\a_2\ran
+(-1)^{|\a_1|+|\a_2|}h(s_1\otimes \a_1, \,\, d^{D''}(s_2\otimes \a_2)\ran,
$$
$$
d h\lan s_1\otimes \a_1, \,\, s_2\otimes \a_2\ran =h\lan d^{D''}(s_1\otimes \a_1),\,\, s_2\otimes\a_2\ran
+(-1)^{|\a_1|+|\a_2|}h(s_1\otimes \a_1, \,\, d^{D'}(s_2\otimes \a_2)\ran,
$$
where $|\a_i|$ denotes the degree of differential form $\a_i.$
\end{lemma}
\bgn{proof}
We have the following: 
$$
d\sig(\a)=
\bgn{cases}
+\sig(d\a), \qquad (|\a|=\text{even})\\
-\sig(d\a), \qquad (|\a|=\text{odd})
\end{cases}
$$
Then the result follows from Lemma \ref{lem:dh(s1, s2)}
and (\ref{eq:hlan s_1otimes a_1,  s_2otimes a_2 ran})
\end{proof}
Let $s\in\Gam(E).$
By the action of the exterior derivative $d$ on both differential forms 
$$h\lan d^{D'}(s\psi), \,\, s\psi\ran\qquad\text{\rm  and }\qquad h\lan d^{D''}(s\psi), \,\, s\psi\ran,$$
we obtain the following equalities of top forms from Lemma \ref{dhlans1otimes a1},
\bgn{align}\label{dhlandD'spsitop}
\(dh\lan d^{D'}(s\psi), \,\, s\psi\ran\)_{top}=&h\lan d^{D''}d^{D'}(s\psi), \,\, s\psi\ran_s-
h\lan d^{D'}(s\psi), \,\, d^{D'}(s\psi)\ran_s
\\
\(d\,h\lan d^{D''}(s\psi), \,\, s\psi\ran\)_{top}=&h\lan d^{D'}d^{D''}(s\psi), \,\, s\psi\ran_s-
h\lan d^{D''}(s\psi), \,\, d^{D''}(s\psi)\ran_s\notag
\end{align}
Since
$D'=d^A+\sqrt{-1}\ome(V)$, $D''=d^A-\sqrt{-1}\ome(V),$ we have
\bgn{align}
d^{D'}d^{D''} (s\psi)=&(d^A+\sqrt{-1}\ome(V))(d^A-\sqrt{-1}\ome(V))(s\psi)\\
=&(F_A+\ome(V)\cdot\ome(V))(s\psi)\\
&-\sqrt{-1}d^A(\ome(V)(s\psi))+\sqrt{-1}\ome(V)d^A(s\psi)\\
\\
d^{D''}d^{D'}(s\psi)=&(F_A+\ome(V)\cdot\ome(V))(s\psi)\\
&+\sqrt{-1}d^A(\ome(V)(s\psi))-\sqrt{-1}\ome(V)d^A(s\psi)
\end{align}

Hence we have
\bgn{lemma}\label{dD''dD'+dD''dD'}
$$(d^{D''}d^{D'}+d^{D'}d^{D''}) (s\psi)=2(F_A+\ome(V)\cdot\ome(V))(s\psi)$$
\end{lemma}

\bgn{lemma}\label{positivety lemma}
Let $s$ be a nonzero section of $E$ satisfying 
$\D^{0,1}s=0$. $(s\neq 0).$ 
We set $i^{+n}\lan \psi, \,\,\ol\psi\ran$ to be a volume form on $M.$
Then the following $2n$-forms are positive with respect to the volume form,
\bgn{align}
&i^{+n}h\lan d^{D'}(s\psi), \,\, d^{D'}(s\psi)\ran_s>0
,\qquad i^{+n}h\lan d^{D''}(s\psi), \,\, d^{D''}(s\psi)\ran_s>0
\end{align}
\end{lemma}
In other words, we have 
\bgn{align}
&\frac{ h\lan d^{D'}(s\psi), \,\, d^{D'}(s\psi)\ran_s}{\lan \psi, \,\,\ol\psi\ran_s}>0,\qquad
\frac{ h\lan d^{D''}(s\psi), \,\, d^{D''}(s\psi)\ran_s}{\lan \psi, \,\,\ol\psi\ran_s}>0
\end{align}
We need the following two lemmas for our proof of Lemma \ref{positivety lemma}.
Let 
$\{e_{-,i}^{1,0}\}_{i=1}^n$ be a unitary basis of $\L_\J^-$ which are locally defined, that is, 
$$\lan e^{1,0}_{-, i}, \, \ol{e^{1,0}_{-, j}}\ran_{\tt}=-\del_{i,j}
$$
\bgn{lemma}\label{lan e1,0 -icodpsiol}
$$
\lan e^{1,0}_{-,i}\cdot\psi, \,\, \ol{e^{1,0}_{-,j}}\cdot\ol\psi\ran_s=2\del_{i,j}\lan\psi, \,\,\ol\psi\ran_s
$$
\end{lemma}
\bgn{proof}
\bgn{align}
\lan e^{1,0}_{-,i}\cdot\psi, \,\, \ol{e^{1,0}_{-,j}}\cdot\ol\psi\ran_s=&
-\lan \ol{e^{1,0}_{-,j}}\cdot e^{1,0}_{-,i}\cdot\psi, \,\, \ol\psi\ran_s\\
&-\lan( \ol{e^{1,0}_{-,j}}\cdot e^{1,0}_{-,i}+{e^{1,0}_{-,j}}\cdot \ol{e^{1,0}_{-,i}})\cdot\psi, \,\,\ol\psi\ran_s\\
=&-2\lan e^{1,0}_{-, i}, \, \ol{e^{1,0}_{-, j}}\ran_{\tt}\lan \psi, \,\,\ol\psi\ran_s\\
=&2\del_{i,j}\lan \psi, \,\,\ol\psi\ran_s
\end{align}
\end{proof}

Let $\{e^{1,0}_{+,i}\}_{i=1}^n$ be a unitary basis of  $\L_\J^+$, that is, 
$$\lan e^{1,0}_{+, i}, \, \ol{e^{1,0}_{+, j}}\ran_{\tt}=+\del_{i,j}
$$

\bgn{lemma}\label{lan e1,0+icdotolpsi}
$$
\lan e^{1,0}_{+,i}\cdot\ol\psi, \,\, \ol{e^{1,0}_{+,j}}\cdot\psi\ran_s=-2\del_{i,j}\lan \ol\psi, \,\,\psi\ran_s
$$
\end{lemma}
\bgn{proof}
The result follows from the same way as before.
\end{proof}

\bgn{proof}[Proof of Lemma \ref{positivety lemma}]
First we shall show 
$$
\frac{ h\lan d^{D''}(s\psi), \,\, d^{D''}(s\psi)\ran_s}{\lan \psi, \,\,\ol\psi\ran_s}>0.
$$
From Lemma \ref{D' D''}, we obtain
\bgn{align}
h\lan d^{D''}(s\psi), \,\, d^{D''}(s\psi)\ran_s=&h\lan (\D''s)\cdot\psi, \,\, (\D''s)\cdot\psi\ran_s
\end{align}
Since $\D''=\D^{1,0}_-+\D^{0,1}_+$ and $\D^{0,1}s=0$, we have $\D^{0,1}_+s=0$.
Thus we obtain 
\bgn{align}
h\lan d^{D''}(s\psi), \,\, d^{D''}(s\psi)\ran_s=&h\lan (\D^{1,0}_-s)\cdot\psi, \,\, (\D^{1,0}_-s)\cdot\psi\ran_s
\end{align}

By using our unitary basis, $\D^{1,0}_- s$ is written by 

$$
(\D^{1,0}_-s)=\sum_i s_i e^{1,0}_{-,i}
$$
Applying Lemma \ref{lan e1,0 -icodpsiol}, we obtain
\bgn{align}
h\lan d^{D''}(s\psi), \,\, d^{D''}(s\psi)\ran_s=&\sum_{i,j}h\lan s_ie^{1,0}_{-,i}\cdot\psi, \,\, 
s_j e^{1,0}_{-,j}\cdot\psi\ran_s\\
=&\sum_{i,j}h(s_i, s_j)\lan e^{1,0}_{-,i}\cdot\psi, \,\, \ol{e^{1,0}_{-,j}}\cdot\ol\psi\ran_s\\
=&2\sum_i h(s_i, s_i)\lan \psi, \,\,\ol\psi\ran_s
\end{align}
Thus we obtain
$$
\frac{ h\lan d^{D''}(s\psi), \,\, d^{D''}(s\psi)\ran_s}{\lan \psi, \,\,\ol\psi\ran_s}=\sum_i h(s_i, s_i)
=\|\D's\|^2\geq 0
$$
Then it turns out that the equality holds if and only if $s=0.$

Secondly, we shall show 
$$\frac{ h\lan d^{D'}(s\psi), \,\, d^{D'}(s\psi)\ran_s}{\lan \psi, \,\,\ol\psi\ran_s}>0.$$
Taking the complex conjugate, we have
$$
\frac{ \ol{h\lan d^{D'}(s\psi), \,\, d^{D'}(s\psi)\ran_s}}{\ol{\lan \psi, \,\,\ol\psi\ran_s}}>0
$$
\bgn{lemma}
$$
\ol{h\lan d^{D'}(s\psi), \,\, d^{D'}(s\psi)\ran_s}=
-h\lan d^{D'}(s\ol\psi), \,\, d^{D'}(s\ol\psi)\ran_s
$$
\end{lemma}
\bgn{proof}
Let $\{\t_i\}$ be a symplectic basis of $T^*_M$ which satisfies 
$\ome=\sum_i \t^{2i-1}\w\t^{2i}$.
Then $D's$ is written as
$D's=\sum_i s_i\t^i$.
Since $ \lan \t^i\ol\psi, \t^i\psi\ran_s=0$, we have  
\bgn{align}
\ol{h\lan d^{D'}(s\psi), \,\, d^{D'}(s\psi)\ran_s}=&\sum_{i,j}\ol{h(s_i, s_j)}\lan \t^i\ol\psi, \t^j\psi\ran_s\\
=&\sum_{i\neq j}\ol{h(s_i, s_j)}\lan \t^i\ol\psi, \t^j\psi\ran_s
\quad \\
=-&\sum_{i,j}h(s_j, s_i)\lan \t^j\ol\psi, \t^i\psi\ran_s\\
=-&\sum_{i,j}h\lan s_j \t^j\cdot\ol\psi, \,\, s_i\t^i\cdot\ol\psi \ran_s\\
=-&h\lan d^{D'}(s\ol\psi), \,\, d^{D'}(s\ol\psi)\ran_s
\end{align}
\end{proof}
Then we obtain
\bgn{align}
h\lan d^{D'}(s\ol\psi), \,\, d^{D'}(s\ol\psi)\ran_s=&
h\lan (\D's)\cdot\ol\psi, \,\, (\D's)\cdot\ol\psi)\ran_s\\
=&h\lan (\D^{1,0}_+s)\cdot\ol\psi, \,\, (\D^{1,0}_+s)\cdot\ol\psi\ran_s
\end{align}

By using our local basis, we obtain
$$
(\D^{1,0}_+s)=\sum_i s_i e^{1,0}_{+,i}
$$
Then applying Lemma \ref{lan e1,0+icdotolpsi}, we have
\bgn{align}
\ol{h\lan d^{D'}(s\psi), \,\, d^{D'}(s\psi)\ran_s}=-&\sum_{i,j}h\lan s_ie^{1,0}_{+,i}\cdot\ol\psi, \,\, 
s_j e^{1,0}_{+,j}\cdot\ol\psi\ran_s\\
=-&\sum_{i,j}h(s_i, s_j)\lan e^{1,0}_{+,i}\cdot\ol\psi, \,\, \ol{e^{1,0}_{+,j}}\cdot\ol\psi\ran_s\\
=&2\sum_i h(s_i, s_i)\lan \ol\psi, \,\,\psi\ran_s
\end{align}
Thus we obtain
$$
\frac{ \ol{h\lan d^{D'}(s\psi), \,\, d^{D'}(s\psi)\ran_s}}{\ol{\lan \psi, \,\,\ol\psi\ran_s}}=\sum_ih(s_i, s_i)
=\|\D''s\|^2>0
$$
Hence we have the result
$$\frac{ h\lan d^{D'}(s\psi), \,\, d^{D'}(s\psi)\ran_s}{\lan \psi, \,\,\ol\psi\ran_s}>0$$
\end{proof}
Applying Lemma \ref{dD''dD'+dD''dD'} to (\ref{dhlandD'spsitop}), we obtain
\bgn{align}
d\(h\lan d^{D'}(s\psi), \,\, s\psi\ran+h\lan d^{D''}(s\psi), \,\, s\psi\ran\)_{[2n-1]}
=&h\lan( F_A+\ome(V)\cdot\ome(V))(s\psi), \,\, s\psi\ran_s\\
-&h\lan d^{D'}(s\psi), \,\, d^{D'}(s\psi)\ran_s\\
-&h\lan d^{D''}(s\psi), \,\, d^{D''}(s\psi)\ran_s
\end{align}
Then we have
\bgn{theorem}
$$
0=\int_M i^{+n}h\lan( F_A+\ome(V)\cdot\ome(V))(s\psi), \,\, s\psi\ran_s-\int_M(\|D's\|^2+\|D''s\|^2)\vol_M
$$
\end{theorem}

A generalized Hermitian connection $\D$  satisfies the Einstein-Hermitian condition if and only if
$$
F_A+\ome(V)\cdot\ome(V)=-i\lam \ome\id, 
$$
where $\lam$ is Einstein constant.
This is equivalent to 
$$
\sqrt{-1}\Lam_\ome(F_A+\ome(V)\w\ome(V))=\lam\,\id
$$
Thus we have
$$
i^n\lan ( F_A+\ome(V)\cdot\ome(V))(s\psi), \,\, s\psi\ran_s=ni^n\lam h(s,s)\lan \psi, \,\,\ol\psi\ran_s
$$
Since $\dstyle{\psi=e^{\frac{\ome}{i}}}$, we have 
\bgn{align}
i^n\lan \frac{\lam \ome}{i}\w\psi, \,\,\ol\psi\ran_s=\lam i^n \frac{\ome}{i}\w
 (\frac{\ome}{i})^{n-1}\frac1{(n-1)!}
\end{align}
and we have
$$\lan\psi, \,\,\ol\psi\ran_s=\frac1{n!}(\frac{\ome}{i})^n$$
Thus we have the following vanishing theorem:

\bgn{theorem}\label{vanishing theorem}
Let $\D$ be a generalized Hermitian connection which satisfies the Einstein-Hermitian condition for negative Einstein constant $\lam<0$. 
Then a generalized holomorphic section $s$ of $E$ which satisfies
$\D^{0,1}s=0$ is zero.
If $\lam =0,$ then a generalized holomorphic section $s$ is parallel.
\end{theorem}
\bgn{proof}We have the following;
$$\int_M i^{+n}h\lan( F_A+\ome(V)\cdot\ome(V))(s\psi), \,\, s\psi\ran_s=
\lam n\|s\|^2_{L^2}\vol(M)=\|D's\|^2_{L^2}+\|D''s\|^2_{L^2}$$
If $\lam<0$, then $s=0$.
If $\lam=0$, then $D's=D''s=0.$ Then $\D s=0.$
\end{proof}
\subsection{In the general cases of $\psi=e^{b-\sqrt{-1}\ome}$}
Let $(E, h)$ be an Hermitian vector bundle over $(M, \J, \J_\psi)$,
where $\psi=e^{b-\sqrt{-1}\ome}$.
Then by the action of the $d$-closed form $b$, we have a generalized K\"ahler structure $(\J_b, \J_{\psi'})$, where 
 $\J_b=e^{^b}\J e^b$ and $\J_{\psi'}=e^{-\sqrt{-1}\ome}$
Let $\D^\A=\D^{1,0}+\D^{0,1}$ be a generalized Hermitian connection satisfying 
the followings: 
 $$
 D^{0,1}\circ \D^{0,1}=0
 $$
$$
\K_\A(\psi)=\lam \id_E,
$$
where $\K_\A(\psi):=\pi^{Herm}_{U^{0,-n}}\F_\A(\psi).$
Since we have 
$$
\F_{\Ad_{e^b}\A}(\psi)=e^b\F_\A({e^{-b}\psi}),
$$
it follows that 
$\K_{\Ad_{e^b}\A}(\psi)=\K_\A(\psi).$
Thus  $\D_b := e^b\circ \D^\A\circ e^{-b}= \D^{\Ad_{e^b}\A}=\D^{1,0}_b+\D^{0,1}_b$ satisfies the Einstein-Hermitian condition with same constant $\lam.$
If $s\in\Gam(E)$ is a section satisfying $\D^{0,1}s=0$, 
Then $\D^{1,0}_b s=e^{-b}\D^{1,0}e^b s=0.$
Hence $s\in\Gam(E)$ is a section with $\D^{0,1}_b s=0$ and $\D_b$ is a generalized connection over $(\J_b, \J_{\psi'}).$
Then applying Theorem \ref{vanishing theorem}, we obtain $s=0.$
Thus we have 
\bgn{theorem}\label{th: vanishing}
Let $(E, h)$ be an Hermitian vector bundle over a compact generalized K\"ahler manifold $(M, \J, \J_\psi)$
of symplectic type. We assume that a generalized Hermitian connection $\D^\A:=\D^{1,0}+\D^{0,1}$ satisfies 
$\D^{0,1}\circ\D^{0,1}=0$ and the Einstein-Hermitian condition 
$$
\K_\A(\psi)=\lam \id_E
$$
If $\lam$ is negative, then every section $s\in \Gam(E)$ satisfying $\D^{0,1}s=0$ is a zero section. 
If $\lam=0,$ then every generalized holomorphic section is parallel.
\end{theorem}
\bgn{proposition}
Let $(E, h)$ be an Hermitian vector bundle over a compact generalized K\"ahler manifold $(M, \J, \J_\psi)$
of symplectic type. We assume that a generalized Hermitian connection $\D^\A:=\D^{1,0}+\D^{0,1}$ satisfies 
$\D^{0,1}\circ\D^{0,1}=0$ and the Einstein-Hermitian condition 
$$
\K_\A(\psi)=\lam \id_E
$$
Then the dual bundle $(\otimes^pE)\otimes (\otimes^q E^*)$ also satisfies the Einstein-Hermitian condition with factor $(p-q)\lam$. In particular, $\w^p E$ satisfies the Einstein-Hermitian condition with factor $p\lam$
\end{proposition}
\bgn{proof}
The result directly follows from the standard connection theory.
\end{proof}
We define the slope of $E$ by the ratio
$$
\mu(E):= \frac{\lan c_1(M)\cdot\psi,\,\ol\psi\ran_s}{\lan \psi, \,\,\ol\psi\ran_s}
$$
\bgn{theorem}
Let $(E, h)$ be an Hermitian vector bundle over a compact generalized K\"ahler manifold $(M, \J, \J_\psi)$
of symplectic type. We assume that a generalized Hermitian connection $\D^\A:=\D^{1,0}+\D^{0,1}$ satisfies 
$\D^{0,1}\circ\D^{0,1}=0$ and the Einstein-Hermitian condition 
$$
\K_\A(\psi)=\lam \id_E
$$
Let $E_1$ be a generalized holomorphic sub bundle of $E$. Then we have the following inequality 
$$
\mu(E_1)\leq \mu(E)
$$
\end{theorem}
\bgn{proof}
Let $E_1$ be a generalized holomorphic subbundle of $E$ of rank $p$. We denote by $j$ the inclusion $j: E_1\to E.$
Then taking $p$th skew-Symmetric power of both sides, we have the map $\w^p: \w^pE_1 \to \w^p E. $
Thus we have a non zero section $s$ of $\w^p\otimes (\w^p E_1)^*$ which is generalized holomorphic. 
Every line bundle admits a generalized Einstein-Hermitian metric and $\w^p E$ also satisfies the Einstein-Hermitian condition. Thus $\w^p\otimes (\w^p E_1)^*$ also satisfies the Einstein-Hermitian condition with factor 
$-p\mu(E_1)+p\mu(E).$
Hence it follows from Theorem \ref{th: vanishing} that we have 
$$
0\leq -p\mu(E_1)+p\mu(E)
$$
Thus we have the inequality 
$$
\mu(E_1)\leq \mu(E).
$$
\end{proof}

\bgn{theorem}
Let $(E, h)$ be an Einstein-Hermitian generalized holomorphic vector bundle over a compact generali zed K\"ahler manifold $(M, \J, \J_\psi).$ 
Then $E$ is $\psi$-poly-stable and $E$ is a direct sum 
$$
(E, h) =\bigoplus_i(E_i, h_i)
$$
where each $(E_i, h_i)$ is $\psi$-stable generalized holomorphic vector bundle with the same factor.
\end{theorem}
\bgn{proof}The result follows from Theorem \ref{th: vanishing}.
\end{proof}



\medskip
\noindent
E-mail address: goto@math.sci.osaka-u.ac.jp\\
\noindent
Department of Mathematics, Graduate School of Science,\\
\noindent Osaka University Toyonaka, Osaka 560-0043, JAPAN 
\end{document}